\documentclass{article}

\usepackage{hyperref}
\usepackage{subcaption}

\usepackage{color}

\usepackage{amssymb, amsmath}
\usepackage[mathscr]{eucal}
\usepackage{graphics,graphicx}
\usepackage{float}
\newcommand{\mycomp}[1]{{#1}^{\:\!^{\scriptstyle \mathsf{c}}}}

\newcommand{\C}{{\mathcal{C}}}

\newcommand{\T}{{\mathcal T}}

\newcommand{\B}{{\mathbb{B}}}

\newcommand{\Cstar}{C^{\,\!\scalebox{0.7}{*}}}
\newcommand{\Cstarhat}{{\widehat C}^{\,\!\scalebox{0.7}{*}}}

\newcommand{\Sstar}{S^{\,\!\scalebox{0.7}{*}}}
\newcommand{\Sstarsmall}{S^{\,\!\scalebox{0.6}{*}}}

\newcommand{\D}{{\mathcal{D}}}
\renewcommand{\H}{{\mathcal{H}}}

\newcommand{\R}{{\mathbb{R}}}

\newcommand{\cA}{{\cal A}}

\newcommand{\wh}{\widehat}

\newcommand{\argmin}{\mathop{\rm argmin}}

\newtheorem{theorem}{Theorem}[section]
\newtheorem{lemma}{Lemma}[section]
\newtheorem{proposition}{Proposition}[section]

\newtheorem{corollary}{Corollary}[section]

\numberwithin{equation}{section}

\begin{document}
\title{Multiscale geometric feature extraction for high-dimensional and non-Euclidean data with applications}
\author{$^1$Gabriel Chandler and $^2$Wolfgang Polonik\\ $^1$Department of Mathematics, Pomona
College, Claremont, CA 91711\\ $^2$Department of Statistics,
University of California, Davis, CA 95616\\[5pt]\small E-mail:
gabriel.chandler@pomona.edu \;and\; polonik@wald.ucdavis.edu}
\date{}
\date{\small \today}
\maketitle 

\vspace*{-1cm}

\begin{abstract}
A method for extracting multiscale geometric features from a data cloud is proposed and analyzed. Based on geometric considerations, we map each pair of data points into a real-valued feature function defined on the unit interval.  Further statistical analysis is then based on the collection of feature functions. The potential of the method is illustrated by different applications, including classification and anomaly detection. 

For continuous data in Euclidean space, our feature functions contain information about the underlying density at a given base point (small scale features), and also about the data depth of the base point (large scale feature). Theoretical investigations show that the method combats the curse of dimensionality, and it also shows some adaptiveness towards sparsity. Connections to other concepts, such as random set theory, localized depth measures and nonlinear multidimensional scaling, are also explored.  
\end{abstract}

\allowdisplaybreaks
\section{Introduction} Extracting qualitative features from a multivariate data cloud is an important task that receives a lot of attention in the literature. A popular approach for the construction of such methods is based on the kernel-trick (e.g. see Shawe-Taylor and Cristianini, 2004) that maps the observed $n$ data points into a usually infinite dimensional feature space. In this work, we propose a new feature extraction method that, similar to the kernel trick, also maps points to an infinite dimensional feature space. However, rather than constructing a feature function for each data point, we construct a {\em matrix of feature functions}: A feature function for each pair of data points. By construction, these feature functions encode certain information about the geometry and shape of the point cloud. This geometric aspect underlying our construction infuses intuition to the resulting methodology that enhances interpretability. Moreover, our feature functions are real-valued functions of a one-dimensional variable, even for multivariate or infinite-dimensional data. Consequently, these functions can be plotted, resulting in a visualization tool. The contributions of this work can be summarized as follows:     \\[5pt]
%
(i) Present a methodology for extracting information about the shape of the underlying distribution based on an iid sample of size $n$ from a Hilbert space, and accomplish this in a multiscale fashion;\\[3pt]
%
(ii) Provide supporting visualization tools; \\[3pt]
%
(iii) Introduce the novel idea of a {\em distribution} of local depth values for a given point;\\[3pt]
%
(iv) Provide supporting theory for the case of Euclidean data, which in particular shows that our method exhibits some inherent adaptivity to sparsity;\\[3pt]
%
(v) Illustrate versatility of the approach through various applications, including anomaly detection and classification;\\[3pt]
%
(vi) Implicitly, some fundamental challenges in non-parametric statistics are addressed: How to choose smoothing/tuning parameters? What is the `right' scale in high dimensions? How to find informative regions or directions in high dimensions?\\[5pt]
%
To heuristically describe the underlying novel idea, assume that the data lie in $\R^d.$  For a given {\em anchor point} $x \in \R^d$, we construct a {\em distribution} of depth values of $x$, and the corresponding quantile function will be the feature function associated with $x$. Defining an anchor point in terms of two data points results in a matrix of feature functions.

The distribution of depths is obtained by randomly selecting subsets of $\R^d$ containing $x$, then finding the depths of $x$ within the subsets.  The notion of depth used here is specified below. The consideration of distributions of depths was motivated by the work of Ting et al. (2013), who consider random split points in the context of what they called `mass estimation'.

By construction, our depth distribution is such that small quantiles (small scales) provide information about the value of the density, while large quantiles (large scales) contain information about the global depth (cf. Lemma~\ref{multiscale}). Intermediate scales might perhaps be most important for feature extraction, particularly in high dimensions. Finding depths within subsets is a way to localize depth. Various local depth measures have been proposed (e.g. Agostinelli and Ramanazzi 2011, Paindaveine and van Bever 2012, Dutta et al. 2016, Kotik and Hublinka 2017, Agostinelli 2018, Serfling 2019), and their multiscale nature has also been discussed. See also section~\ref{local-depth}.
It should also be pointed out that we do not consider the quantile level as a tuning parameter, but we consider the entire quantile function. This is similar to the concepts underlying the mode tree (Minotte and Scott, 1993), the siZer (Chaudhuri and Marron, 1999), or to  Betti-curves or other representations of persistence diagrams considered in topological data analysis (TDA), e.g. Bubenik and Kim (2007), and Bubenik (2015). These methods also do not choose a tuning (smoothing) parameter, but consider all possible values of the parameter. However, for both small and large values of the tuning parameter, the resulting quantities do not contain useful information as they are degenerate. This is in stark contrast to the methodology presented here, where all the values of the tuning parameter give informative quantities. Also, the methods just mentioned summarize the entire data cloud in essentially one function (or in a few, as in the case of Betti curves), while our approach constructs functions for all (pairs of) individual data points. 

In contrast to the kernel trick, our construction is not related to an RKHS. Instead, our construction is more directly guided by the objective to extract geometric information. The computational burden might be somewhat higher than for RKHS-based methodologies, but given the currently available computing power, it is manageable. The algorithmic complexity of our procedure is quadratic in the sample size and linear in the dimension, meaning that a high dimension is not the major computational challenge. It turns out that, to a certain extent, the same applies to the statistical behavior. 

The remainder of the paper is organized as follows. The construction of population versions of our feature functions in $\R^d$ are described in section~\ref{pop}, and corresponding empirical versions are given in section~\ref{feature-constr}. Section~\ref{object} discusses generalizations to Hilbert space data.  Several applications of our approach are presented in section~\ref{applications}. Section~\ref{relations} discusses  relations to other concepts, while the choice of tuning parameters is discussed in section~\ref{tuning}. Theoretical analyses of our feature functions are presented in section~\ref{theory}. This includes robustness to the curse of dimensionality, some adaptation to sparsity, and asymptotic distribution theory, both pointwise and as a process. All the proofs are presented in section 9.

\section{Distribution of depths and depth quantile functions}
\label{pop}
The basic idea underlying our approach is to define a distribution of depths of a given point $x \in \R^d, d \ge 2,$ 
by randomly selecting subsets of ${\mathbb R}^d$ containing $x$ and finding the depth of $x$ within this subset.  In an attempt to balance complexity and computational cost, the usual challenge in high-dimensional situations, our approach uses (right spherical) {\em cones} as subsets. 

In what follows, we only consider the case $d \ge 2$. Let $F$ be a distribution on $\R^d$ with Lebesgue density $f$. In slight abuse of notation, we also denote by $F(t)$ with $t \in \R^d$ the corresponding distribution function, and below we follow the same convention for distributions denoted by other symbols. For $x\in \R^d$ and a vector $u$ on the unit sphere, a line $\ell = \ell_{x,u} \subset \R^d$ is given by $\ell_{x,u} = \{ y = y(t) = x + tu,\, t \in\R\}$. Notice that $\ell_{x,u} = \ell_{x,-u}$. 

Then, for $(x,u) \in \R^d \times S^{d-1}$ and $s \in \R,$ we denote by $C_{x,u}(s)$ the (right spherical) cone with fixed opening angle $\alpha \in [0,\pi)$ that contains $x,$ has $\ell_{x,u}$ as its axis of symmetry, and whose tip $y(s) = x + su$ has distance $|s|$ to $x$. Note that this includes cones opening to both sides, where the orientation of the cone $C_{x,u}(s)$ depends on the sign of $s$, where, by convention, the orientation of $C_{x,u}(0)$ is chosen to be the same as for $s > 0.$ For $s \ne 0$, we also have $C_{x,u}(s) = C_{x,-u}(-s)$. The vector $x$ is serving as our anchor point. In practice, we will have to specify a rule for how to choose an anchor point, and this rule might depend on the specific application considered (see sections~\ref{feature-constr} and \ref{tuning}). 

Recall that $\langle z - x, u \rangle$ is the signed distance between $x \in \R^d$ and the projection of $z \in \R^d$ onto the line $\ell_{x,u}.$ With $X \sim F$, we define the depth of $x$ within $C_{x,u}(s)$ as
\begin{align}\label{def-local-depth}
d_{x,u}(s) = \min \big(P\big[\big\{X \in\, &C_{x,u}(s)\big\} \cap \big\{\langle X - x,u \rangle \le 0\big\}\big],\nonumber\\
&P\big[\big\{X \in C_{x,u}(s)\big\} \cap \big\{\langle X - x,u \rangle \ge 0\big\}\big]\,\big).
\end{align}
We can interpret $d_{x,u}(s)$ as the Tukey depth of $x$ with respect to the (in general improper) one-dimensional distribution $F_{x,u,s}$ on $\ell_{x,u}$, defined by $F_{x,u,s}(t) = P(X \in C_{x,u}(s), \langle X - x, u \rangle \le t).$ As $x \in\R^d$, `Tukey depth' of $x$ might be associated with half-space depth, which is not what we mean. Instead $d_{x,u}(s)$ is the (one-dimenional) Tukey depth of the origin of the one-dimensional  and in general improper distribution with cdf $F_{x,u,s}(t),$ i.e., we have $d_{x,u}(s) =  \min\big\{F_{x,u,s}(0), F_{x,u,s}(\infty) - F_{x,u,s}(0)\big\},$ where $F_{x,u,s}(\infty) = \sup_{t\in \R}F_{x,u,s}(t)$ denotes the total mass of $F_{x,u,s}.$ A different way of writing this is
\begin{align}\label{alternative-def}
d_{x,u}(s) = \min\big(F(A_{x,u}(s)), F(B_{x,u}(s))\big),
\end{align}
where
\begin{align}\label{A-B-sets-def}
A_{x,u}(s) = C_{x,u}(s) \cap H_{x,u}^-\quad\text{and}\quad B_{x,u}(s) = C_{x,u}(s) \cap H_{x,u}^+,
\end{align}
and $H_{x,u}^\pm$ are the two closed half-spaces defined by the hyperplane $H_{x,u} = \{z\in {\mathbb R}^d:\, \langle z,u\rangle = \langle x,u\rangle\}$ with normal direction $u$, passing through $x,$ and $H_{x,u}^-$ contains the tip $y(s) = x + su.$ 
%
%
%
\begin{figure}[h]
\centerline{\includegraphics[height=2.3in, width = 3.1in]{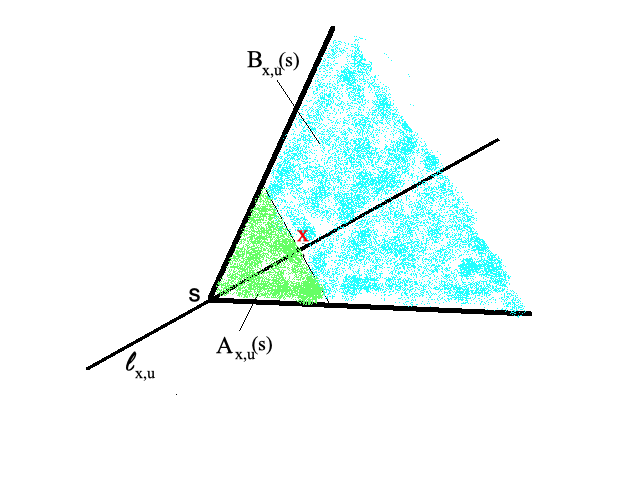}}
\vspace*{-1cm}
\caption{\small{$A$ and $B$-sets are illustrated. The minimum probability content of the two sets equals the depth of $x$ within the indicated cone.}}\label{A-B-sets}
\end{figure}
In other words, the hyperplane $H_{x,u}$ splits the cone $C_{x,u}(s)$ into two closed subsets $A_{x,u}(s)$ and $B_{x,u}(s)$ with $A_{x,u}(s) \cap B_{x,u}(s) = H_{x,u}  \cap C_{x,u}(s).$  By definition, $A_{x,u}(s)$ contains the tip $y(s)$, i.e. $A_{x,u}(s)$ is itself a cone with base $H_{x,u}  \cap C_{x,u}(s)$, and $B_{x,u}(s) = C_{x,u}(s) \setminus A_{x,u}(s)$ is a frustum. 
%
%
Figure \ref{A-B-sets} illustrates this representation for a given line and a given anchor point.

Notice that even for $|s| \to \infty$, we are not obtaining Tukey half-space depth of $x$ because, for a fixed $u$, we only consider parallel half-planes orthogonal to $u$. If we, however, consider all directions for a given $x$, then we can obtain (Tukey) half-space depth as $\inf_{u \in S^{d-1}}\max_{s \in \R}d_{x,u}(s).$  Note that $d_{x,-u}(s) = d_{x,u}(-s)$, and thus these  functions depend on the direction of $u$, not just on the line $\ell_{x,u}.$ This should not cause any confusion, however. 

The depth quantile functions defined next are the crucial objects. Empirical versions of them will be used as feature functions (see below). Assume that $s$ is randomly chosen according to a distribution on $\R$, meaning that we randomly choose the tip of the cone (or its distance to the anchor point $x$ along with the orientation of the cone) on a given axis of symmetry. Formally, given a distribution $G$ on $\R$, let $S \sim G,$  and define $q_{x,u}(\delta)$ as the quantile function of the resulting depth distribution of $x$, i.e.
\begin{align}\label{def-depthity-mult}
q_{x,u}(\delta) = \inf\{t \ge 0: G(s: d_{x,u}(s) \le t) \ge \delta\}. 
\end{align}
%
%
The distribution $G$ can be considered to be a tuning parameter. Often, we will assume $G$ to be a continuous distribution; see section~\ref{tuning} for some discussion on the choice of $G$. 

The following two lemmas describe some properties of the depth function $d_{x,u}(s)$ and the corresponding depth quantile function $q_{x,u}(\delta)$. We denote by $I_{x,u}(t)$ the sublevel sets of $d_{x,u}(s)$, i.e.
$$I_{x,u}(t) = \{s \in\R: d_{x,u}(s) \le t\}.$$
\vspace*{-0.4cm}

\begin{lemma}\label{quant}   Let $(x,u) \in \R^d \times S^{d-1}$ be fixed. \\[3pt]
(a) The function $s \to d_{x,u}(s), s \in \R$ is non-increasing for $s \le 0$ and non-decreasing for $s \ge 0.$ Thus, its sublevel sets $I_{x,u}(t), t \ge 0,$ are either empty or intervals containing $0$.\\[3pt]
(b) Let $F$ be continuous with positive Lebesgue density on $\R^d$, and let $m_{x,u} = \sup_s d_{x,u}(s).$ The function $s \to d_{x,u}(s), s \in \R$ is continuous with minimum value 0, and it is {\em strictly} decreasing on $\{s\le 0: d_{x,u}(s) < m_{x,u}\}$, and strictly increasing on $\{s\ge 0: d_{x,u}(s) < m_{x,u}\}$. The intervals $I_{x,u}(t), t \ge 0,$ are thus all closed.\\[3pt]
(c) (rigid transformations and scale changes) Let $\theta = (\mu,\sigma,O)$ with $\mu \in {\mathbb R}^d, \sigma > 0$ and $O \in \R^{d \times d}$ orthogonal. If $F_{\theta}(t) = F\big(\tfrac{O^{-1} (t - \mu)}{\sigma}\big), t \in \R^d$, and $d_{x,u}^{\,\theta}(s)$ denotes the depth function corresponding to $F_{\theta}$,  while $q_{x,u}(\delta)$ is the one for $F$, then
\begin{align}
    d^{\,\theta}_{x,u}(s) = d_{\sigma O x + \mu,Ou}\big(\sigma s\big),\;\;s \in \R.
\end{align}
\end{lemma}

Let $F_{x,u}(t)$ denote the cdf of the one-dimensional distribution obtained by projecting all the mass onto the line $\ell_{x,u}$. 

\begin{lemma} 
\label{multiscale}
In parts (a)-(c) assume that both $F$ and $G$ have positive Lebesgue densities $f$ and $g,$ respectively. Let $(x,u) \in \R^d \times S^{d-1}$ be fixed. \\[3pt] 
(a)  With $I_{x,u}(q_{x,u}(\delta)) = [s^l_{x,u}(\delta), s^r_{x,u}(\delta)],\,\delta \in (0,1)$, we have 
\begin{align}\label{depth-quant-cont}
q_{x,u}(\delta) =  d_{x,u}(s^l_{x,u}(\delta)) = d_{x,u}(s^r_{x,u}(\delta)).
\end{align}
Moreover, $q_{x,u}(\delta)$ is strictly increasing with $q_{x,u}(0) = 0,$ and if $g$ is symmetric about zero then $q_{x,u}(\delta) = q_{x,-u}(\delta)$ for all $\delta \in [0,1].$ \\[3pt]
%
(b)\;(large scales) \;$q_{x,u}(1) = \min\big(F_{x,u}(x), 1 - F_{x,u}(x)\big)$ is the global Tukey depth of $x$ for the distribution with cdf $F_{x,u}(t)$. \nonumber \\[3pt]
(c)\; (small scales)\; $\lim_{\delta \to 0} \frac{q_{x,u}(\delta)}{\delta^d} = c_d \frac{f(x)}{g^d(0)}$, where $c_d > 0$ is a known constant.\\[3pt]
 (d)\; (invariance) Let $\theta = (\mu,\sigma,O)$ (as in Lemma~\ref{quant}(d)), and let $G_\theta(t) = G\big( \tfrac{t}{\sigma}\big)$. If $q^{\theta}_{x,u}(\delta)$ denotes the depth quantile function defined with  $G_\theta$ and $F_{\theta}$, while $q_{x,u}(\delta)$ is the one for $F$, then, for any $\theta$,
 $$q^{\theta}_{\sigma Ox + \mu,Ou}(\delta) = q_{x,u}(\delta).$$
\end{lemma}

Assertion (c) makes precise the above statement that, for small values of $\delta,$ the depth quantile function $q_{x,u}(\delta)$ contains information about the density at $x$. Note that this localization is achieved even though we are not using local neighborhoods of $x$. Instead, localization follows from the use of Tukey depth and the choice of cones to define the distribution of depths. (It is the acute angle of the cone that is important here; cf. proof of assertion (c).) Property (d) indicates a certain type of scale-invariance, provided the change in scale in $F$ is reflected in $G$ also. (See discussion at the end of section~\ref{feature-constr} for more on this.)\\[5pt]
{\em Interpretation of depth quantile functions.} By definition, $d_{x,u}(s) = F(C^*(s))$ with $C^*(s) \in \{A_{x,u}(s), B_{x,u}(s)\}$. Thus, if we were to consider only negative or only positive values of $s$, the function $d_{x,u}(s)$ can be interpreted as a generalized distribution function (only evaluated on a line) of a sub-probability measure. (Notice, though, that even for fixed $(x,u)$ the sets $C^*(s)$ might not necessarily be nested.) The total mass of this measure equals the `Tukey depth of the base point $x$ among the projections of all the data on the line.' This distribution changes with the direction $u$ and with the base point $x$. Choosing the tip randomly and considering the quantile functions allows both the combining of the negative and positive values of $s$, as well as bringing the various `distribution functions' onto the same scale. This, in turn, makes averaging more reasonable, for instance (see section~\ref{classification}). Also note that the quantile functions are not quantiles of the `distribution functions' $d_{x,u}(s)$, but they are quantiles of the distribution of the random variable $d_{x,u}(S)$ where $S \sim G$. As a result, the depth quantile function $q_{x,u}(\delta)$ is a reparametrization of $d_{x,u}(s)$ (see \ref{depth-quant-cont}). A more geometric construction of $q_{x,u}(\delta)$ is as follows. Consider the level sets of $d_{x,u}(s)$ at level $t$, which are the intervals $I_{x,u}(t).$ Fix $\delta > 0$, and imagine increasing $t$ until the probability content of $I_{x,u}(t)$ under $G$ equals $\delta.$ Then this value $t = t_\delta$, i.e. the value of the depth function at the limits of the interval $I_{x,u}(t_\delta),$ equals $q_{x,u}(\delta)$. (This is (\ref{depth-quant-cont}).) As for comparing two quantile functions: $q_{x_1,u_1}(\delta) > q_{x_2,u_2}(\delta)$ indicates a higher mass concentration about $x_1$ than about $x_2$ (at least when `looking' in directions $u_1$ and $u_2$, respectively), in the sense that there exist intervals about the respective anchor points with the same $G$-mass $\delta$ (the level sets $I_{x_1,u_1}(t_\delta)$ and $I_{x_2,u_2}(t_\delta)$, respectively), such that the sets $C^*$ corresponding to the endpoints of $I_{x_1,u_1}(t_\delta)$ carry a higher mass content than the ones corresponding to $I_{x_2,u_2}(t_\delta)$. See also the discussion on average depth quantile functions in section~\ref{random-sets}.

\section{Construction of the feature functions \texorpdfstring{$\widehat q_{ij}(\delta)$}{a} for data in \texorpdfstring{$\R^d$}{a}} 
\label{feature-constr}

\begin{figure}[h]
\begin{subfigure}{0.4\textwidth}
\hspace*{-0.3cm}
\includegraphics[height=1.5in, width= 2.2in]{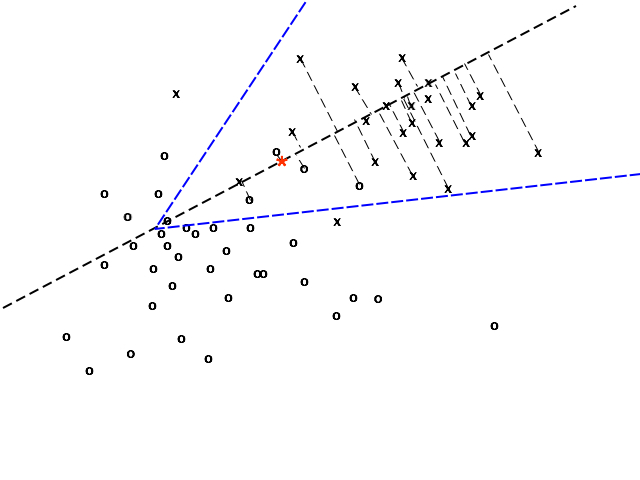}
\end{subfigure}
\hspace*{0.7cm}
\begin{subfigure}{0.4\textwidth}
\includegraphics[height=1.5in, width= 2.2in]{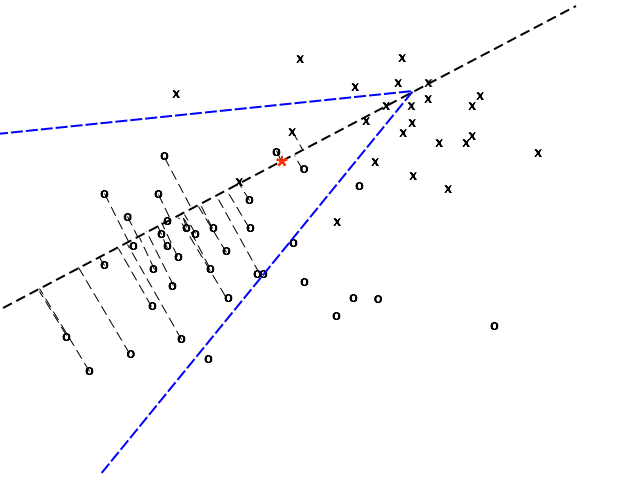}
\end{subfigure}
\caption{\small Illustrates the construction of $\widehat d_{ij}(s)$ in dimension $d=2$ for a 2-class classification problem (different characters correspond to different classes).  The two panels correspond to one positive and one negative value of $s$. The dashed line is the line $\ell_{ij}$; the anchor point $m_{ij}$ is shown in red; the cones $C_{ij}(s)$ are shown in blue; projection of the data inside the cone onto $\ell_{ij}$ are indicated.}
\label{projection}
\end{figure}

Given iid data $X_1,\ldots,X_n$ from $F$, an empirical counterpart of $q_{x,u}(\delta)$ is obtained by simply replacing the distribution $F$ by the empirical distribution $F_n$ resulting in 
\begin{align}\label{emp-quant-simple}
\widehat{q}_{x,u}(\delta) = \inf \{t \in \R: G(s: \wh d_{x,u}(s) \le t) \ge \delta)\},
\end{align}
where $\widehat{d}_{x,u}(s) = \min(F_n(A_{x,u}(s)), F_n(B_{x,u}(s)).$ The question now is: how to pick both $x$ and $u$. We propose to choose them depending on the data as follows. 
For  each data pair $X_i,X_j$, $i\neq j$, let $u_{ij} = \tfrac{X_i - X_j}{\|X_i - X_j\|}$ be the direction given by $X_i$ and $X_j$, and let  $m_{ij} =\frac{X_i + X_j}{2}$ be the anchor point. The anchor point and the direction define the line $\ell_{ij} = \{y\in {\mathbb R}^d: y = y(t) =  m_{ij} + t u_{ij}, t \in {\mathbb R}\}$ passing through $X_i$ and $X_j.$  This gives
%
%
%
$$\wh d_{ij}(s) = \min\big(F_n(A_{ij}(s)), F_n(B_{ij}(s))\big),$$
where $A_{ij}(s)$ and $B_{ij}(s)$ is the sets $A_{x,u}(s), B_{x,u}(s)$  with  $x = m_{ij}$ and $u = u_{ij}.$ The corresponding empirical depth quantile function $\wh q_{ij}(\delta)$ is 
\begin{align}\label{emp-quantile}
\wh q_{ij}(\delta) = \inf \big\{t: G\big(\{s: \wh d_{ij}(s) \le t\}\big) \ge \delta  \big\}.
\end{align}
The functions $\wh q_{ij}(\delta)$ are our feature functions. 
An alternative way of thinking about  constructing $\wh d_{ij}(s)$ (see Fig.~\ref{projection}) is as follows. 
Project all the data inside $C_{ij}(s) = C_{m_{ij},u_{ij}}(s)$ onto $\ell_{ij},$ and let $\widehat d_{ij}(s)$ be the `(one-dimensional) Tukey depth of $m_{ij}$ among these projections'. 
%
%
%
Similar to the populations version, this means the following: Each projection of a point $X_k \in C_{ij}(s)$ can be written as $y(t_k),$ where  $t_k  = \langle X_k - m_{ij}, u_{ij}\rangle \in \R$ is the projection score, so that, in particular, $m_{ij} = y(0).$ The value $\wh d_{ij}(s)$ then is the (scaled by the total sample size $n$) Tukey depth of $0$ among the $t_k$'s.
%
%
%

In addition to the distribution $G$ on $\ell_{ij}$ of the tip, the angle $\alpha$ and the rule to choose the anchor point $m_{ij} \in \ell_{ij}$ are also tuning parameters. Our choice of the anchor point as the midpoint of $X_i,X_j$ was originally motivated by the application of the proposed methodology to classification (see section~\ref{classification}), though consideration of anchor points other than the actual observations appears to have some benefits, at least in high dimensions. See section~\ref{tuning} for some discussion of the choice of these tuning parameters.

Observe that, in analogy to the population version discussed in section~\ref{pop},  the empirical depth functions $\wh d_{ij}(s)$ and $\wh d_{ji}(s)$ are mirrored versions of each other. If $G$ is symmetric about zero, then they lead to the same feature functions, i.e. in this case $\wh q_{ij}(\delta) = \wh q_{ji}(\delta)$. 




\section{Depth quantile functions for object data, the kernel-trick, $Z_1$-$Z_2$-plots, and non-linear multi-dimensional scaling}
\label{object}

\subsection{Object data and the kernel trick} 
The approach discussed above can also be applied to object data in a Hilbert space. The simple, but key, underlying observation here is that determining our depth quantile functions only requires the ability to define cones, to calculate which data points fall into the cones, and to find distances between points. All of this, of course, hinges on the existence of an inner-product.  So, if our underlying data are objects $O_1,\ldots,O_n \in {\cal O},$ lying in a Hilbert space $H$ equipped with an inner-product $\langle\cdot,\cdot \rangle_H,$ then we can simply apply our methodology outlined above to the object data with this inner-product. 

Interestingly, applying our methodology to functional data in $L^2$ will map functional data into a different set of functional data, and it remains to be investigated how such an approach compares to working with the functional data directly. 

Another possible application of our methodology is in combination with the kernel trick. The kernel trick consists of mapping  data $O_i$ into an RKHS $H$ via a feature map $\Psi$, say, such that the dot-product $\langle \cdot,\cdot\rangle_H$ satisfies $\langle \Psi(O_i), \Psi(O_j)\rangle_H = K(O_i,O_j)$ for some kernel $K: {\cal O} \times {\cal O} \to \R$. After applying this map, we can simply apply our  methodology outlined above to the transformed data $\Psi(O_i),\,i = 1,\ldots,n$ with the corresponding dot-product in the RKHS. Kernel methods are available for a variety of objects, including trees, graphs, matrices, strings, tensors, functions, persistence diagrams, etc. (e.g. see Genton 2001, Shawe-Taylor and Cristianini, 2004, Cuturi 2010). Since the elements in the feature space $\Psi(O_i)$ lie in a Hilbert space, we can apply our methodology outlined above. This can be used to investigate and to compare the geometry of the data in feature spaces. 

\subsection{Relation to non-linear multidimensional scaling: $Z_1$-$Z_2$-plots} 
Given a pair of data, say, $O_i,O_j$ with $\ell = \ell_{ij}$ passing through $O_i$ and $O_j$, and the anchor point $m_{ij} = \frac{O_i+ O_j}{2},$ all we need to determine the depth functions $\wh d_{ij}(s)$, or the quantile functions $\wh q_{ij}(\delta)$, are the set of two-dimensional points $\{(Z_{1k}^{ij},Z_{2k}^{ij}),\,k=1,\ldots,n\},$ where
\begin{align*}
Z^{ij}_{1k} = \Big\langle O_k - m_{ij}, \frac{O_i - O_j}{\|O_i-O_j\|} \Big\rangle \quad\text{\rm and}\quad Z^{ij}_{2k} = \sqrt{\|O_k - m_{ij}\|^2 - |Z_{1k}|^2\,},
\end{align*}
where $Z^{ij}_{1k}$ is the (signed) distance between $m_{ij}$ and the projection of $O_k$ onto $\ell_{ij}$, and $Z^{ij}_{2k}$ is the distance of $O_k$ to $\ell_{ij}.$ Checking whether $O_k$ lies in $C_{ij}(s)$ with $s < 0$, say, simply means checking whether $Z_{1k}^{ij} \ge s,$ and if yes, determining whether $ Z_{2k}^{ij} \ge \tan(\alpha)(Z_{1k}^{ij} -s).$ 
In other words, the pairs $\big(Z^{ij}_{1k},Z^{ij}_{2k}\big) \in \R^2,\,k = 1,\ldots,n,$ fully determine $\wh d_{ij}(s),$ and thus $\wh q_{ij}(\delta)$. 

A plot of the pairs $(Z^{ij}_{1k},Z^{ij}_{2k})$, $k\neq i,j$, is what we call a $Z_1$-$Z_2$-plot. 
It is some kind of non-linear projection onto two dimensions, which is in the spirit of multidimensional scaling. However, in contrast to multidimensional scaling, we do not just have one plot, but a total of $n(n-1)/2$ such plots. Each (unordered) pair $(O_i,O_j)$ determines one $Z_1$-$Z_2$-plot, which provides a visual impression of the geometry of the data cloud if viewed in direction $\ell_{ij}$. 
%
%
One can think of the function $\wh q_{ij}(\delta)$ as a summary of the $Z_1$-$Z_2$-plot. 
Another reason for introducing these plots is as follows. Suppose we want to perform binary classification of high-dimensional data. The collection of $Z_1$-$Z_2$-plots can be used to visually inspect whether good linear binary classification is possible. To this end let $\ell_{ij}$ be a line connecting two points $(O_i,O_j)$ and consider the corresponding $Z_1$-$Z_2$-coordinates. Observe that a hyperplane orthogonal to the line $\ell_{ij}$ corresponds to a vertical line in the $Z_1$-$Z_2$-plot. Thus, if we can find a pair of points $O_i,O_j$ such that the corresponding two-dimensional $Z_1$-$Z_2$-coordinates allow good classification with a vertical line, then there exits a good binary linear classifier in the original space.

As an illustration of a possible use of the $Z_1$-$Z_2$-plots, we consider an infinite dimensional situation. Here the underlying space is the RKHS generated by the popular radial basis kernel $K_{\sigma}(u,v)=\mbox{exp}(-\|\frac{u-v}{\sigma}\|^2)$.  The objects $O_j$ are the corresponding feature functions created by applying the kernel-trick to two classes (Virginica vs. Versicolor) of Fisher's famous Iris data, and we visualize the geometry of the cloud of feature functions by using the RKHS geometry as described above. 

We are interested in exploring the effect on the geometry of the feature functions induced by the choice of the tuning parameter $\sigma.$ Fixing a pair $(O_i,O_j)$ and varying $\sigma$ results in a class of $Z_1$-$Z_2$ plots $\{(Z_{1k}^{ij}(\sigma),Z_{2k}^{ij}(\sigma)), k = 1,\ldots,n\}, \sigma > 0.$ Alternatively, one can also think of them as a collection of curves $\sigma \to $ $ (Z_{1k}^{ij}(\sigma),Z_{2k}^{ij}(\sigma)),$ $k = 1,\ldots,n.$

\begin{figure}[h]
\centerline{\includegraphics[width=8cm,height=5.5cm]{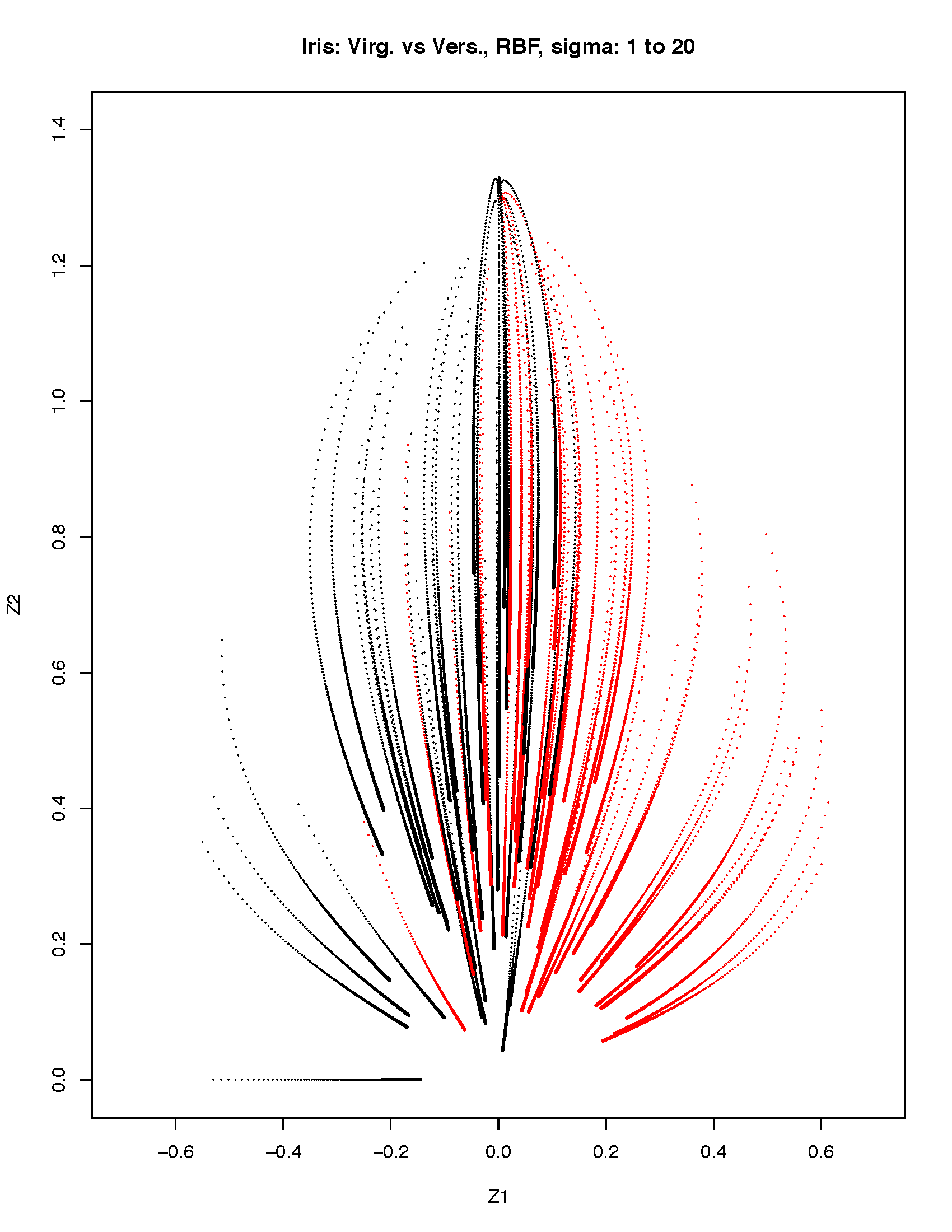}}
\caption{{\small Overlaid $Z_1$-$Z_2$ plots, visualizing point cloud geometry from RBF kernel with varying $\sigma$ parameter value for two classes of the Iris data (color coded). Larger values of $\sigma$ correspond to larger dot sizes. Each curve corresponds to an individual observation as its $Z_1$-$Z_2$-coordinates  change as a function of $\sigma$.}   
}
\label{fig:changesigma}
\end{figure}

Figure \ref{fig:changesigma} shows such a collection of curves by varying $\sigma$ from 1 to 20. It turns out that, in our application, the curves hardly overlap, so that they can be tracked visually. The two colors (black, red) correspond to the two classes. The peculiar shape of the curves can partly be explained as follows. When $\sigma\approx 0$ (top part of the plot), all points in the feature space are essentially orthogonal as $K(u,v)\approx 0$ for $u\neq v$, and all points in the $Z_1$-$Z_2$-plot are approximately $(0,\sqrt{3/2}).$ The latter can be seen by observing that the norm of all the feature functions equals 1 (they all live on the unit sphere in the RKHS), and  projecting orthogonal unit vectors onto the line $\ell_{ij}$ connecting two of them gives $m_{ij}$, the midpoint between the two points determining the line. The distance to this midpoint is $\sqrt{3/2}.$ Since the projections equal the midpoint, their distance obviously equals zero. This then results in the point $(0,\sqrt{3/2})$ in the $Z_1$-$Z_2$-coordinates. As $\sigma\to\infty$, the points converge to (0,0) for all the distances in the RKHS converge to zero.



\section{Applications}\label{applications}
\footnote[1]{R code for computing the depth quantile functions is available at \href{https://github.com/GabeChandler/depthity/}{https://github.com/GabeChandler/depthity/}}Here we illustrate possible types of applications of the depth quantile functions $\wh q_{ij}(\delta)$. The goal of this exposition is to establish that these objects contain valuable information regarding the location of individual observations within a point cloud as well as the associated geometry of the point clouds. How best to use this information likely depends on the particular task at hand, and constitutes future work beyond the scope of the present paper. 

For all of the tasks discussed below, we summarize the $\binom{n}{2}$ depth quantile functions by considering averages of the form $\wh q_i(\delta) = n^{-1}\sum_{j\ne i}\wh q_{ij}(\delta)$, or with the sum restricted to certain subsets of the data (see section~\ref{classification}).  As a result, we have a new representation of our data that can be used to perform a variety of statistical tasks, such as EDA, classification, and anomaly detection.  Other statistical tasks may require combining these functions in a different way, or adapting the choice of the tuning parameters (cf section~\ref{tuning}), which the flexibility of this idea allows.  In all of the following, we use as our base distribution $G$ a uniform over a support large enough to ensure that at the boundary of the support, all the data lie inside the cone for all $u$ and all $x$ under consideration. We use $\alpha=\pi/2$ unless otherwise noted. 


\subsection{Topological data analysis, or discovering holes}
\begin{figure}[h]
    \centering
    \begin{subfigure}[t]{0.32\textwidth}
        \centering
        \includegraphics[height=1.4in, width=\linewidth]{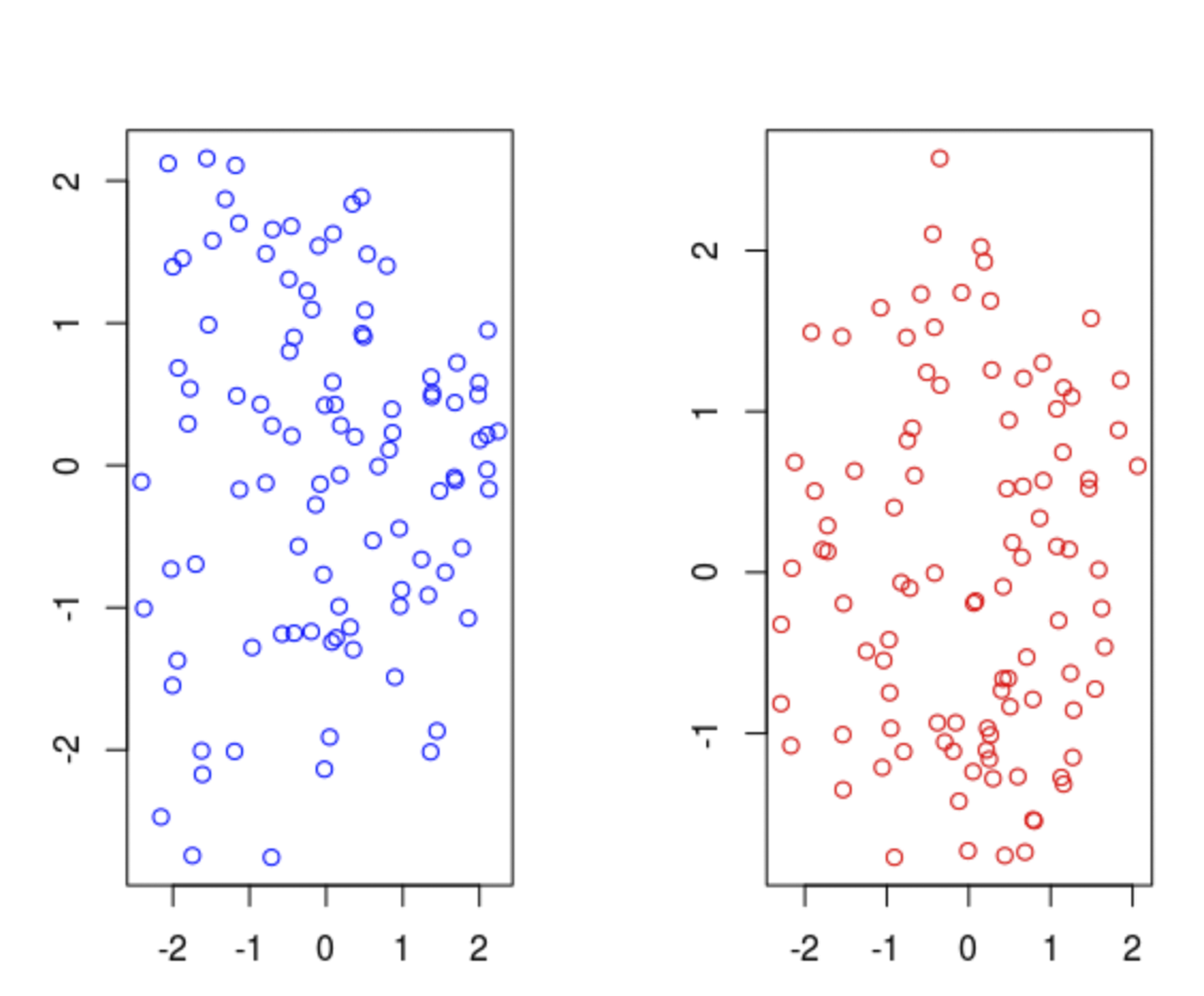} 
        \caption{MDS for 8-annulus (left) and 8-ball (right).} \label{fig:mds8ball}
    \end{subfigure}
    \hfill
        \centering
    \begin{subfigure}[t]{.32\textwidth}
        \includegraphics[height=1.4in, width=\linewidth]{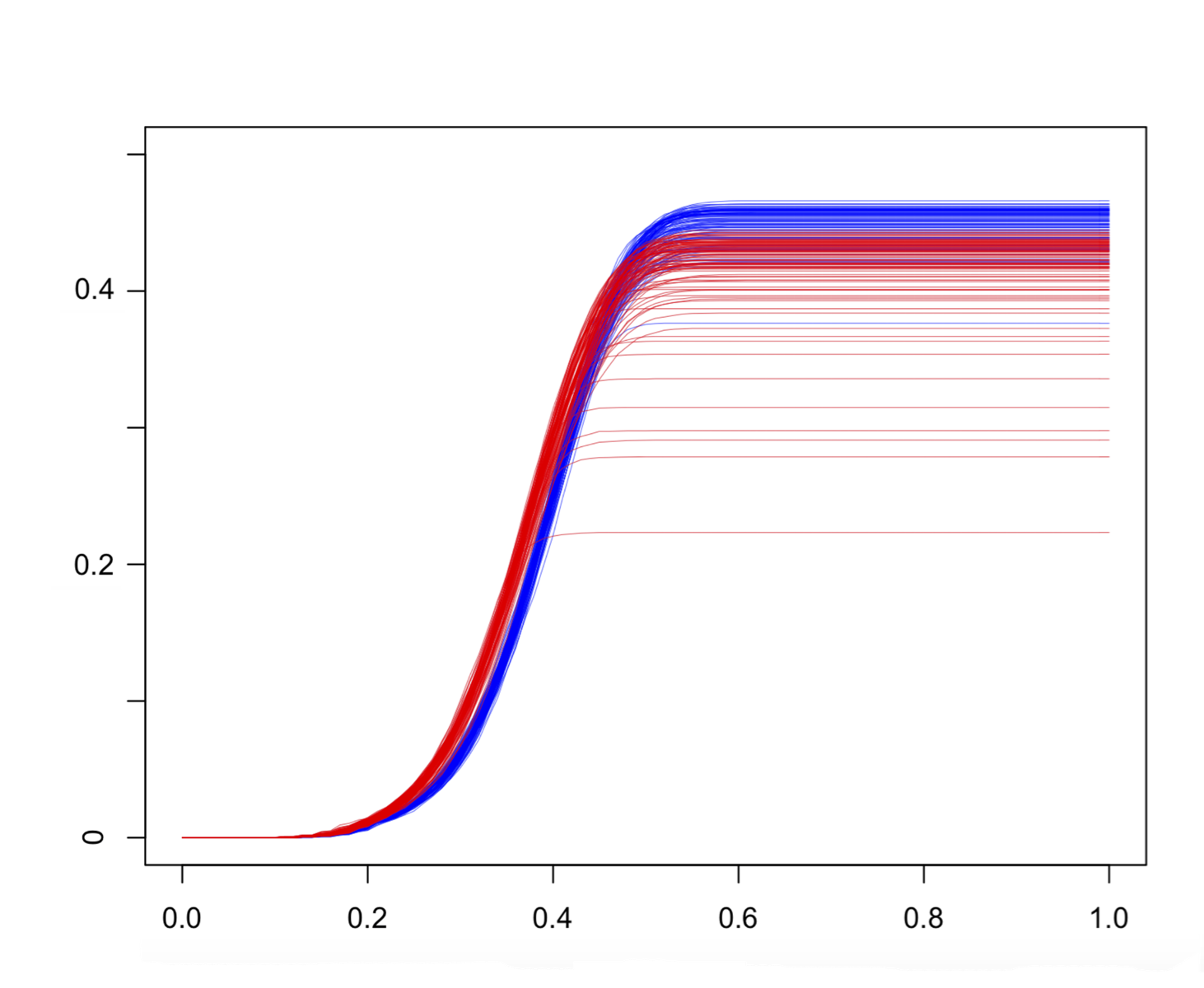} 
        \caption{Averaged Depth Quantile Functions. (Annulus: Blue, Ball: Red) } \label{fig:8balldepthquant}
    \end{subfigure}
\hfill
    \begin{subfigure}[t]{0.32\textwidth}
        \centering
        \includegraphics[height=1.4in, width=\linewidth]{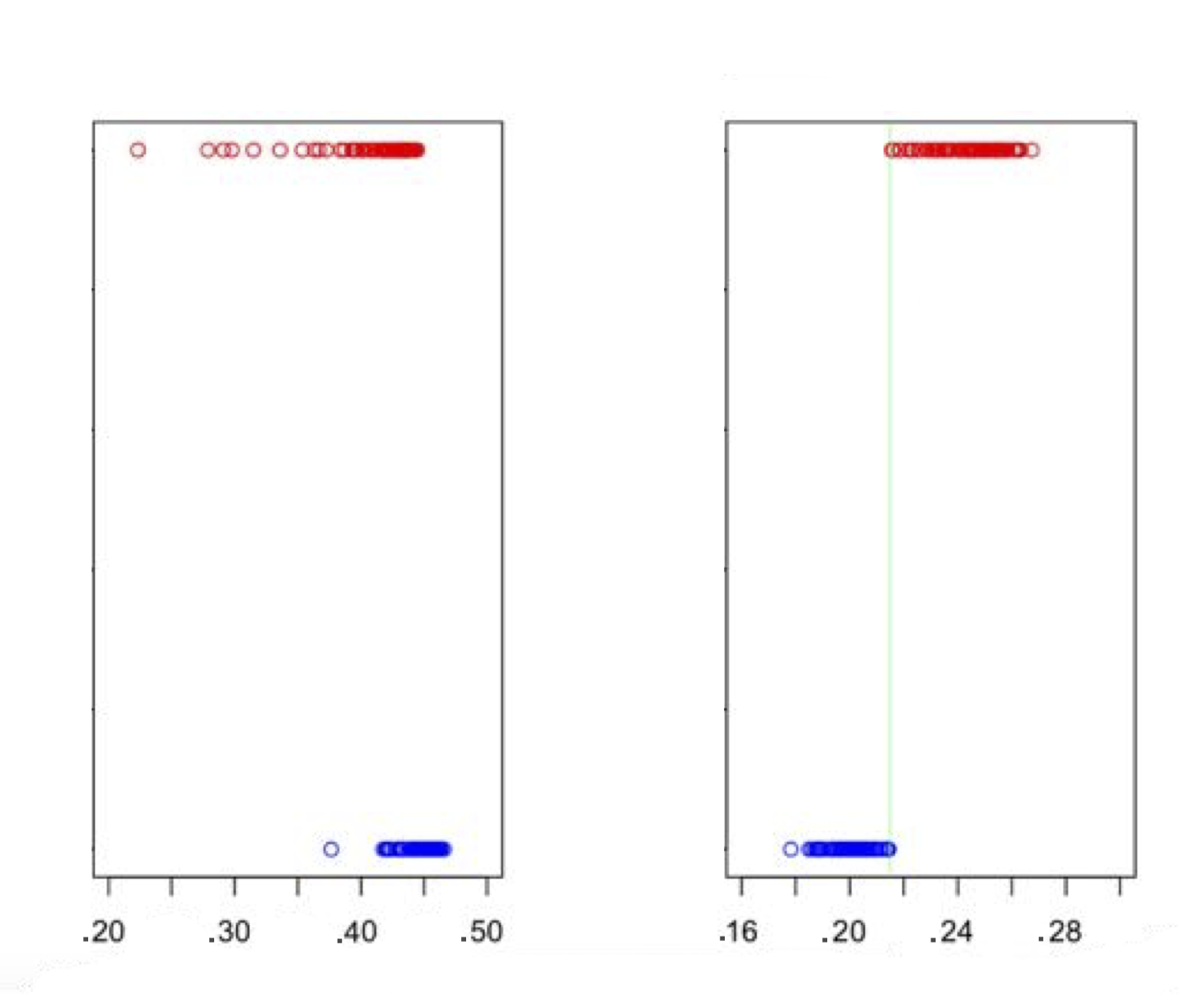} 
        \caption{Global depth performance, i.e. $\hat{q}_i(1)$ (left) vs $\hat{q}_i(.39)$ (right), with separating line (green).} \label{fig:8ball_global}
    \end{subfigure}
   \caption{\small Comparison of observation level information for distinguishing between 8-ball and 8-annulus.\\[-18pt]}
   \label{fig:8ball_indy_level}
\end{figure}

Consider observing a point cloud in $\R^8$ where the data is either generated uniformly in a ball of fixed radius or is generated uniformly in this ball with a smaller middle ball removed.  We call these two supports the 8-ball and 8-annulus respectively.  Specifically, our larger ball has radius 1.5 and the inner ball has radius 1.25, accounting for 25 percent of the total volume.  We wonder whether, given a sample of size of $n=100$ from each model, it is possible to detect if a hole is present.  Figure  \ref{fig:mds8ball} shows that a 2-dimensional multidimensional scaling (MDS) with distances based on a weighted $K$-nearest neighbors graph ($K=5$ shown here) seems unable to recognize the difference between the two supports.  However, the depth quantile functions show quite different behavior between the two situations, and considering all intermediate values provides more information than just considering, for instance, the right hand boundary.  In fact, as shown in figure~\ref{fig:8ball_indy_level}c, the intermediate  value $\delta = .39$ gives perfect separation of the two classes, unlike other values, such as the boundary value $\delta = 1$.  While this suggests that a single functional observation is sufficient to discriminate between the two classes, unsurprisingly, using all the functions together yields even stronger information.  

\subsection{Classification}
\label{classification}
Here, we assume our data is of the form $(X_i, Y_i)\in \R^d \times \{0,1\}$.  This setting provides the original motivation for considering midpoints $m_{ij}$ as our anchor points, as in the case of modal clustering, midpoints for within class comparisons would tend to live inside of a cluster (high density, low depth) while between class comparisons would result in midpoints between the clusters (low density, high depth). As we further demonstrate, the way that the depth quantile function transitions between these two extremes provides additional useful information for classification, even outside of the case of convex clusters (cf. section~\ref{simul}).

To perform binary classification, we consider the following heuristic. We summarize each observation via {\em two} functions: $$\wh{q}^{k}_i(\delta)=\frac{1}{\#\{y_j=k\}-{\bf 1}(y_i=k)}\sum_{j\neq i, y_j=k}\wh{q}_{ij}(\delta),\qquad k=0,1.$$  Functional principal component analysis (fPCA) is then used to associate each of these function with a four dimensional loading vector.  Doing this for all functions associates an 8-dimensional vector with each observation, with the first 4 values corresponding to $\hat{q}_i^0(\delta)$.  An `out-of-the-box' support vector machine is used to perform classification on these vectors, with results given below.  Extensions to $m$-class problems can be accomplished similarly by constructing $4m$-vectors associated with every observation. 
\subsubsection{Simulated Data}\label{simul}
\begin{figure}[h] 
    \centering
    \begin{subfigure}[t]{0.32\textwidth}
        \centering
    \includegraphics[width=\linewidth]{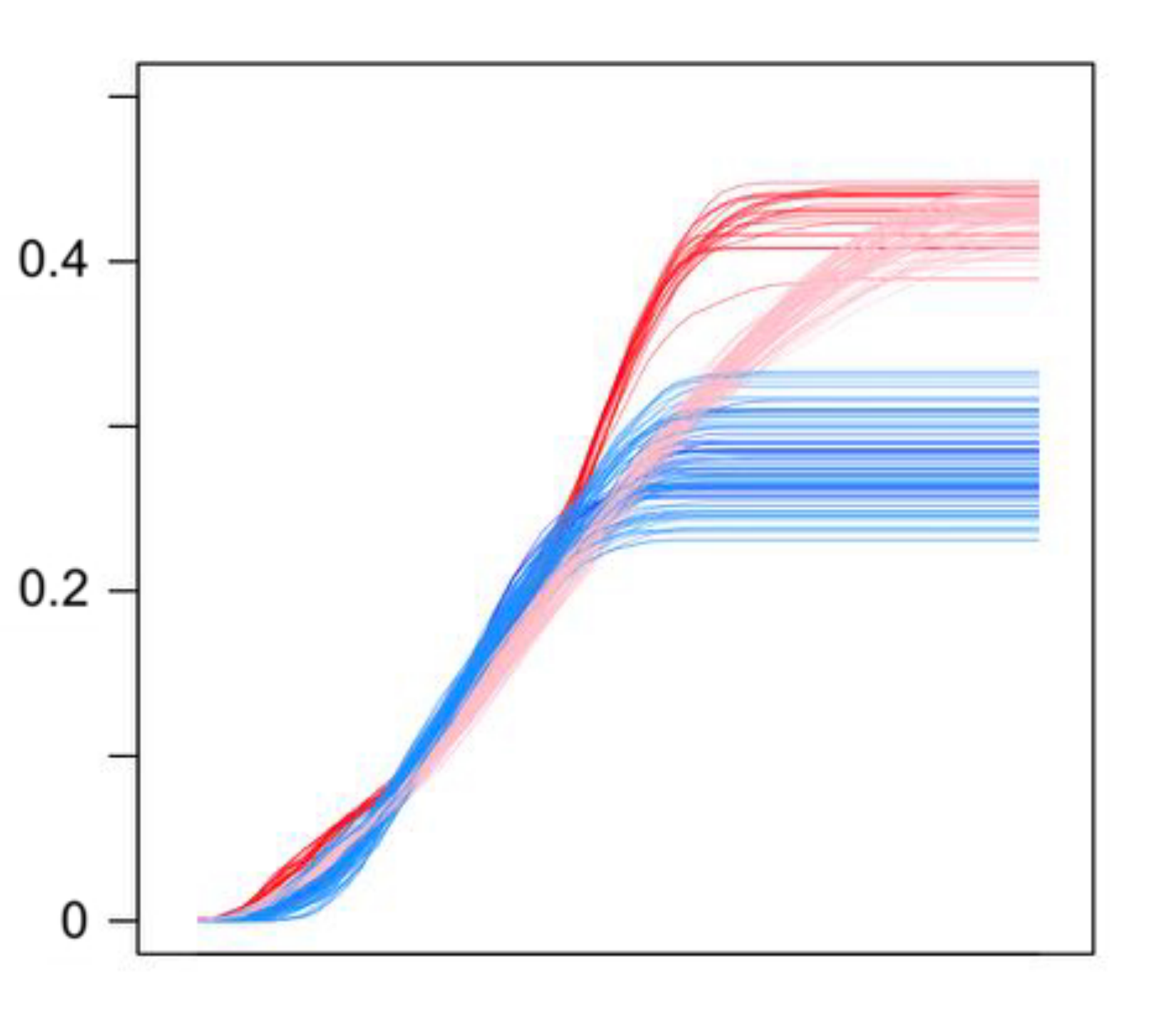} 
        \caption{\small $\hat{q}_i^k$ based on $x\in \boldsymbol{\R^2}$} \label{fig:cir_data}
    \end{subfigure}
    \hfill
        \centering
    \begin{subfigure}[t]{.32\textwidth}
        \includegraphics[height = 1.4in, width=\linewidth]{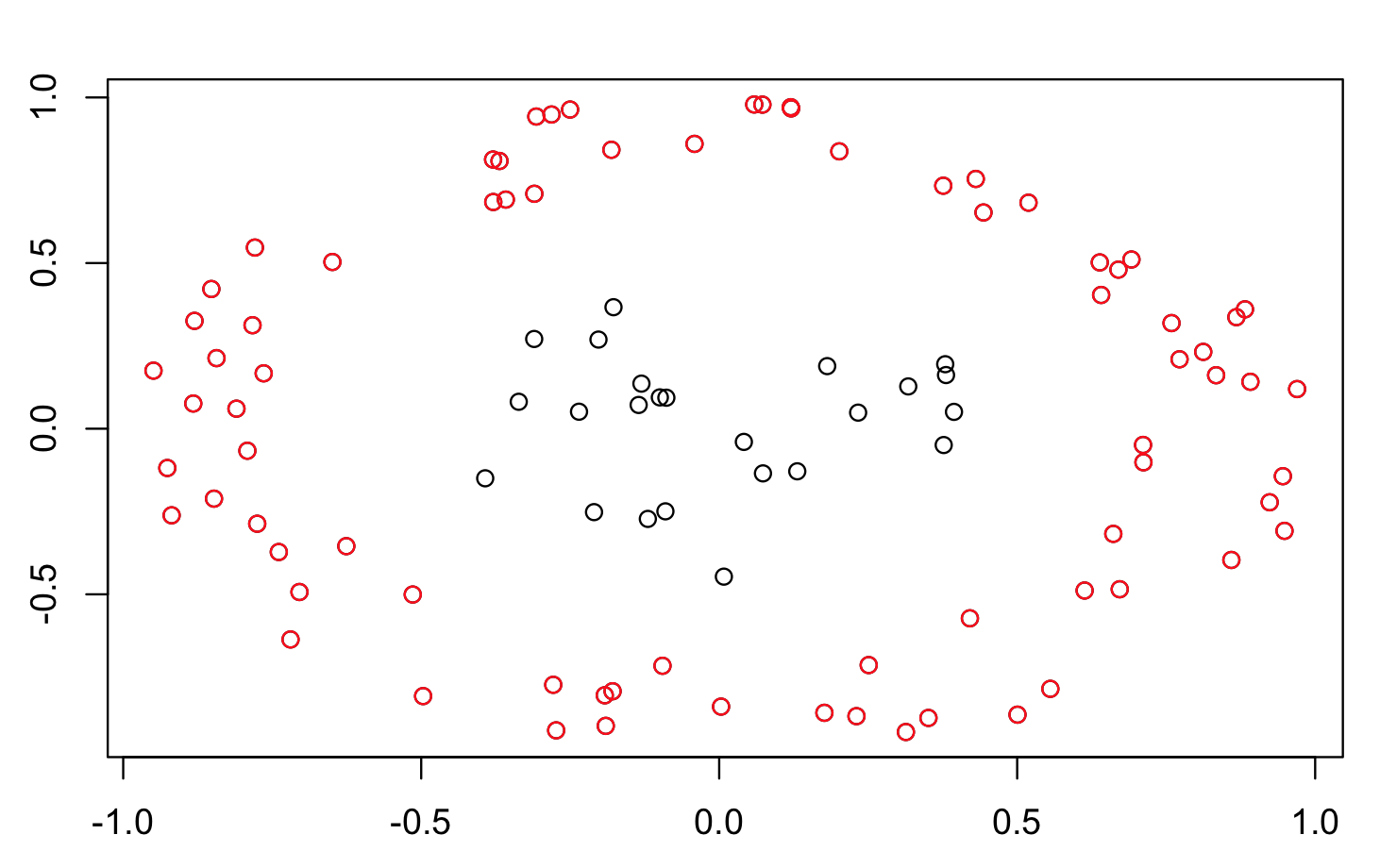} 
        \caption{\small Data} \label{fig:cir_lin}
    \end{subfigure}
\hfill
    \begin{subfigure}[t]{0.32\textwidth}
        \centering
        \includegraphics[width=\linewidth]{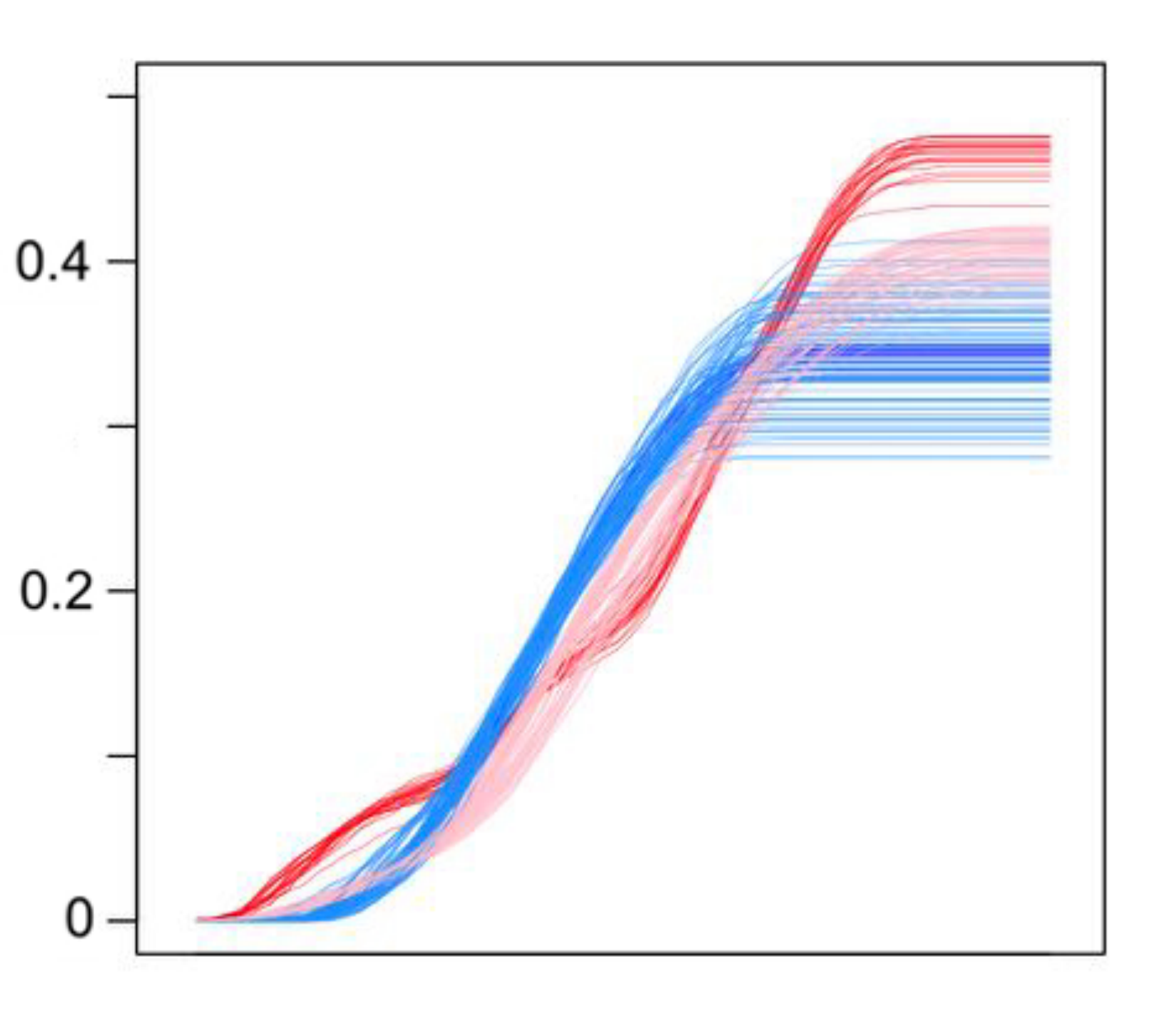} 
        \caption{\small $\hat{q}_i^k$ based on $\phi(x)\in \boldsymbol{\R^3}$} \label{fig:cir_kern}
    \end{subfigure}
    \caption{\small (a)  The averaged quantile functions $\hat{q}_i^k(\delta)$ (blue-between class comparisons, red-within the inside group, and pink-within the outside group, (b) a simulated data set of size $n=100$, and (c)  $\hat{q}_i^k(\delta)$ for data in the feature space.}
    \label{fig:nonconvex}
\end{figure}

The first classification data set we consider suggests that convex modal regions are not necessary for successfully extracting information relevant for classification with our method.  Figure  \ref{fig:nonconvex} considers a 2-dimensional data set in which the supports of the two classes are a disc and a concentric annulus, such that a third annulus exists between the two with density zero.  

Comparisons between classes tend to yield midpoints in this zero density region, and thus valuable information seems to be revealed by small $\delta$ values.  The 
mapping $\phi((x_1,x_2))=(x_1,x_2,x_1^2+x_2^2)$ makes the data linearly separable though causes the support to live on a elliptic paraboloid and thus all midpoints live in a zero density region.  Interestingly, it is the intermediate values of $\delta$ that now reveal useful information (see figure \ref{fig:cir_kern}).

\subsubsection{Real Data}
\begin{figure}[h]%
   \begin{subfigure}[c]{0.31\textwidth}
        \centering
        \includegraphics[width=\linewidth]{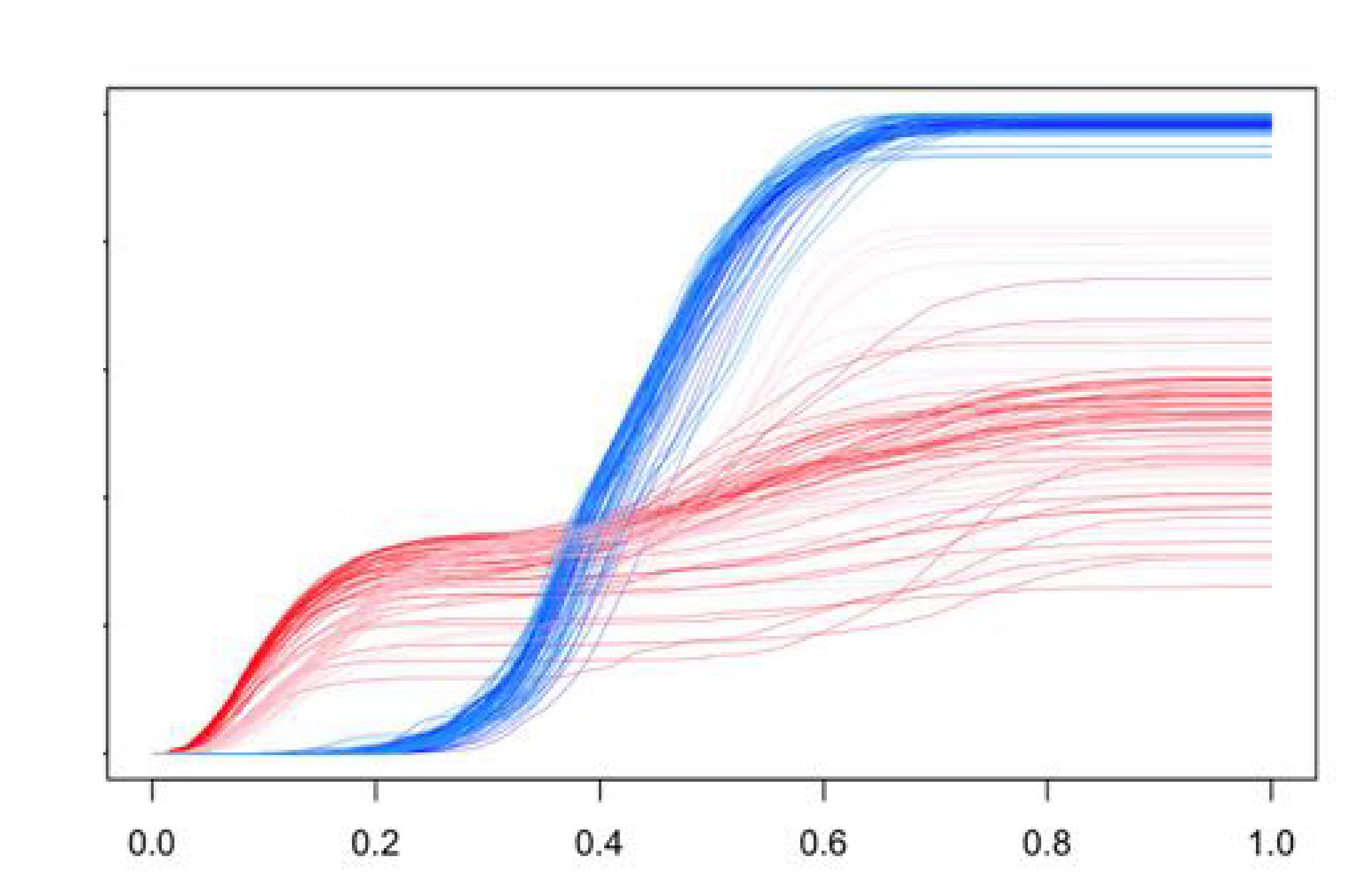} 
        \caption{Iris Data: Setosa vs. Versicolor} \label{fig:iris_set_ver}
    \end{subfigure}
    \hfill
         \begin{subfigure}[c]{0.32\textwidth}
        \centering
        \includegraphics[width=\linewidth]{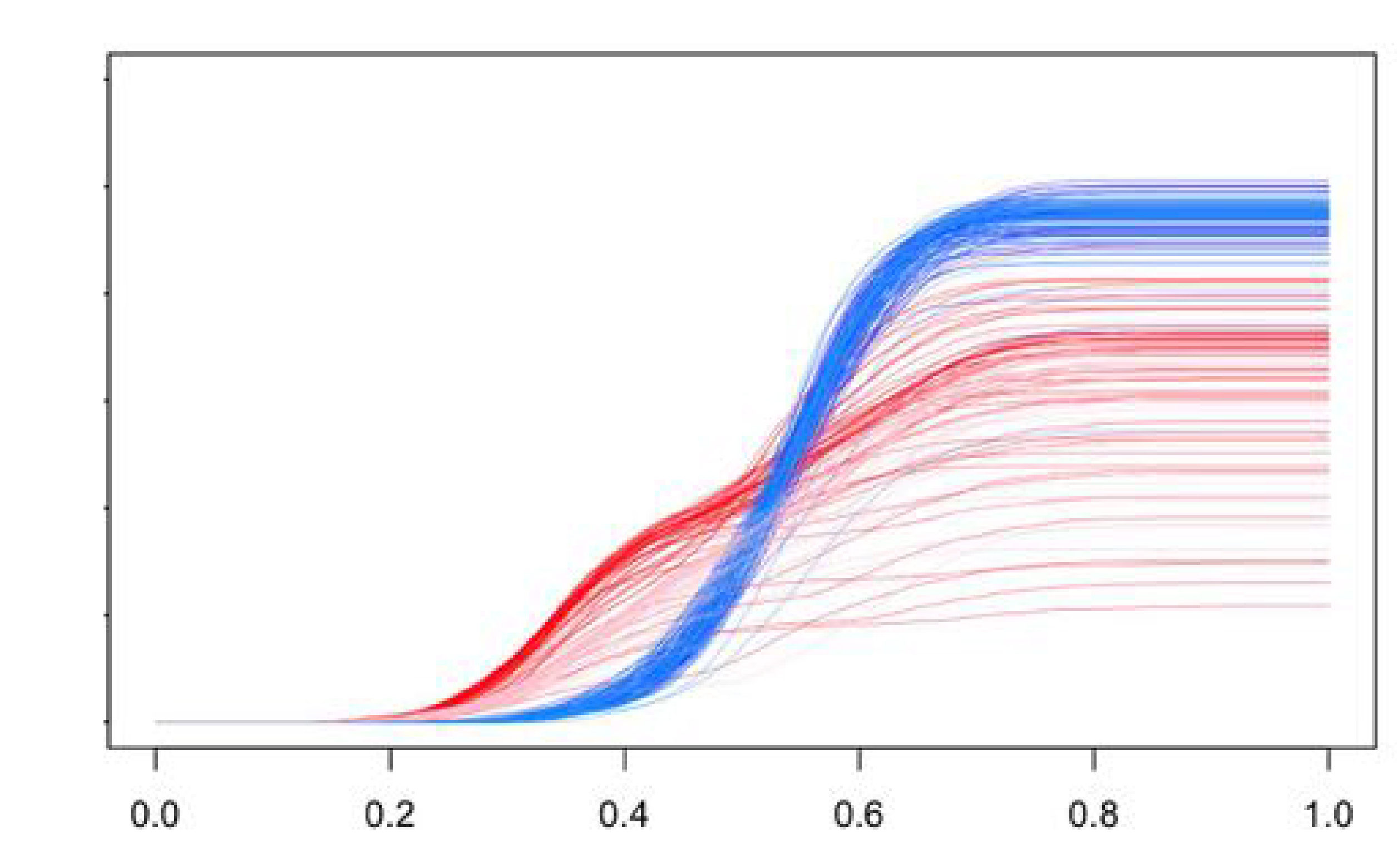} 
        \caption{Wine Data: Class 1 vs Class 3} \label{fig:win12}
    \end{subfigure}
        \hfill
         \begin{subfigure}[c]{0.32\textwidth}
        \centering
        \includegraphics[width=\linewidth]{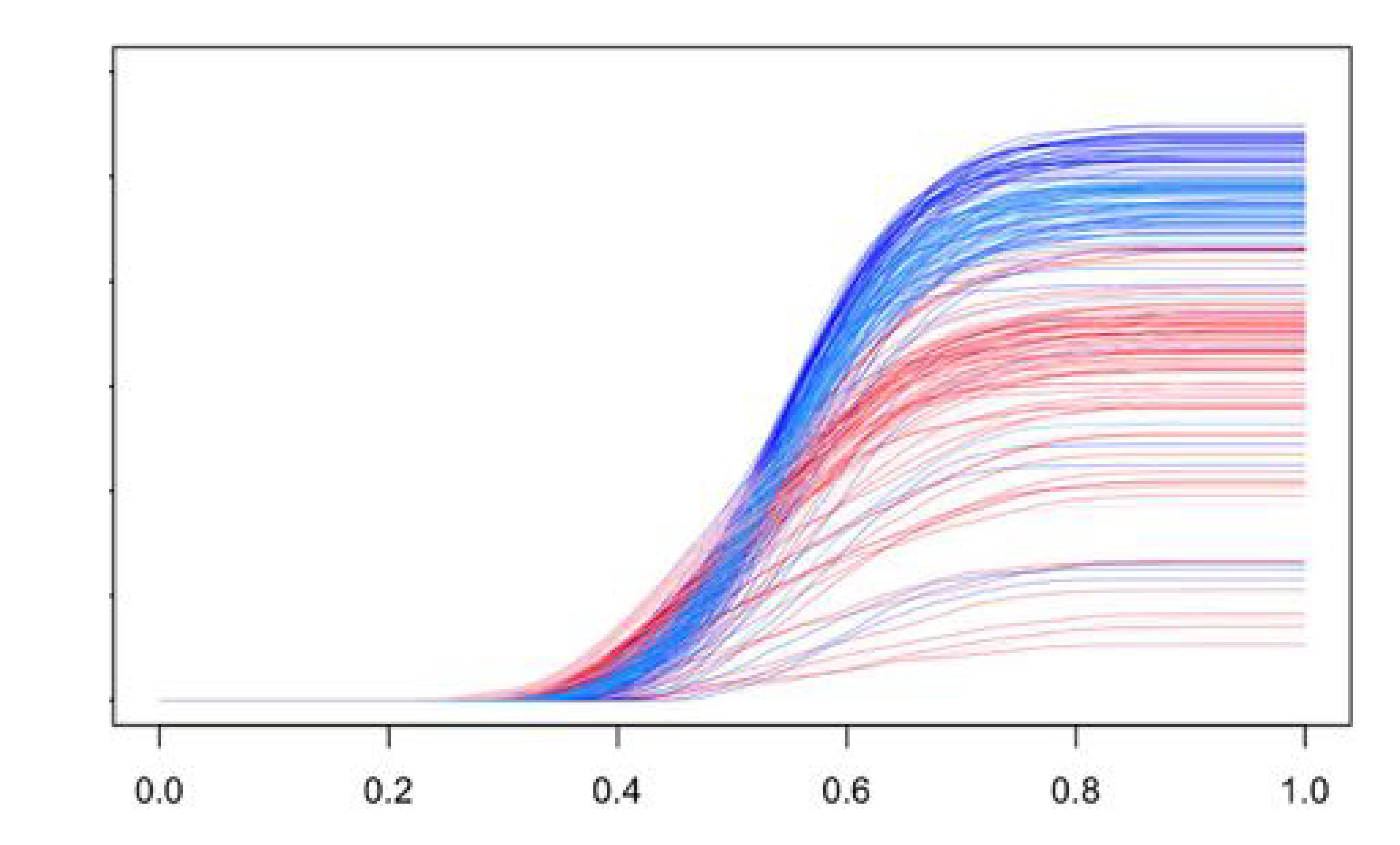} 
        \caption{Wine Data: Class 2 vs Class 3} \label{fig:win13}
    \end{subfigure}
\vspace*{-0.2cm}
   \caption{ \small $\hat{q}_{i}^k(\delta)$ for Iris and Wine data sets.  (Note: since $\hat{q}_{ij}(1)\leq .5$ by construction,  the scale of the vertical axis always ranges from 0 to slightly less than .5, and thus we omit it for the remainder of the plots for cleanliness of presentation.)\\[-8pt]}%
  \label{fig:iris}%
\end{figure}

We next illustrate the performance of our classification routine on several well known data sets.  For each data set, all variables are first $z-$scaled (all dimensions have mean 0 and variance 1).   First is Fisher's Iris data, consisting of 3 classes of $n_i=50$ observations each, with $x_i\in\boldsymbol{\R^4}$ (note that we actually perform the classification task in double the dimension in this example).  We evaluate the leave-one-out performance of the method, and find a correct classification rate of 97.33\%. For comparison, a $k$-nearest neighbors ($k$NN) classifier optimized with $k $=13 achieves a rate of 96.67\% on the $z$-transformed data.  Figure \ref{fig:iris_set_ver} demonstrates that not only are the two points clouds situated apart from each other, but they seem to have fairly different geometry as well.
 
The second real data example we consider is the {\it Wine} data set, available from the UCI Machine Learning database.  There are 178 observations over 3 classes in 13 dimensions.  Our method yields a leave-one-out correct classification rate of 96\%. A $k$NN ($k$=25) achieves a rate of 97.7\%, though other classifiers have shown even better performance, see Aeberhard et al. (1992). Figures \ref{fig:win12} and \ref{fig:win13} shows the functions for the easiest pairwise comparison (1 misclassification) and the hardest (3 misclassifications), respectively.  

\begin{figure}[h]
    \centering
        \includegraphics[width=\linewidth, height=1.6in]{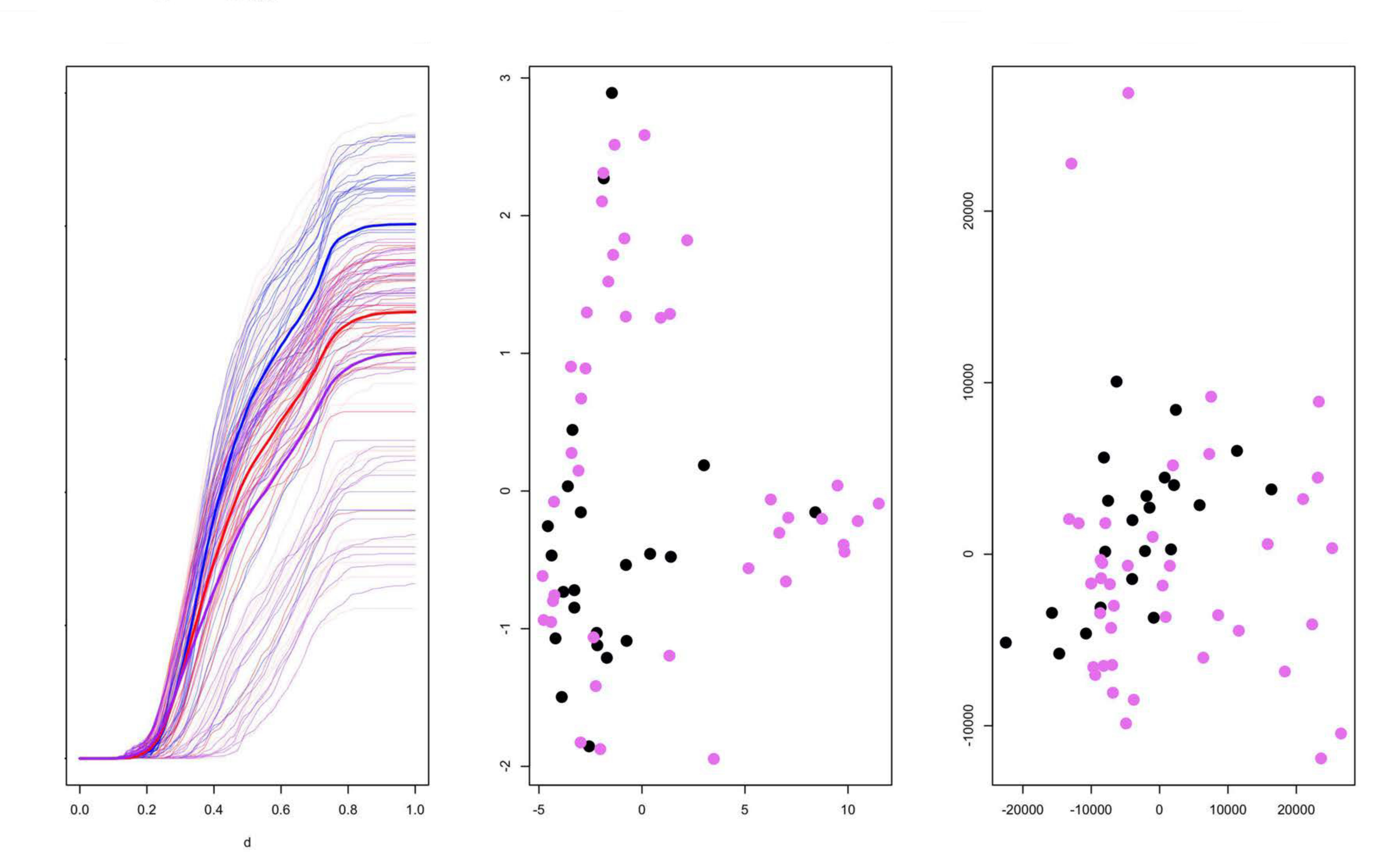} 
         \caption{\small Averaged depth quantile functions (left), and the first two dimensions of PCA based on averaged depth quantile functions (center) and raw data (right). Left: Purple and blue are within group comparisons for cancer and no-cancer, respectively. Red and pink: Between group comparisons. Group level averages are in bold (red and pink averages coincide by construction).}\label{fig:cancerdata}
\end{figure}

%
We next consider a gene expression data set used for colon cancer detection (Alon et al., 1999), consisting of 62 observations on 2000 variables.  We use an opening angle of 100 degrees for our cones.  Plots of the depth quantile functions (figure \ref{fig:cancerdata}) seem to suggest that the observations corresponding to the presence of cancer live in a higher density region (early initial dominance of the purple curve) near the boundary of the point cloud (small depth values, despite nearly 2/3 of the observations corresponding to this class).  Averages over classes are added to better show the information contained by the collection of functions. A visualization of the first two values of a principal components analysis (PCA) based on a fPCA yields more structure for the classification problem than that based on the distance matrix of the raw data. 
The method returns a correct classification rate of nearly 84\% (82\% if angle of $\pi/2$ is used). $k$NN ($k=13$) correctly classifies only 74.6\%.  A better comparison is to a linear SVM (88.7\% with cross-validated $C=.002)$, which demonstrated the best performance in the simulation study of Dutta et al. (2016) under a different evaluation scheme.

\begin{figure} [h]
    \centering
    \begin{subfigure}[t]{0.49\textwidth}
        \centering
        \includegraphics[height = 1.6in, width=\linewidth]{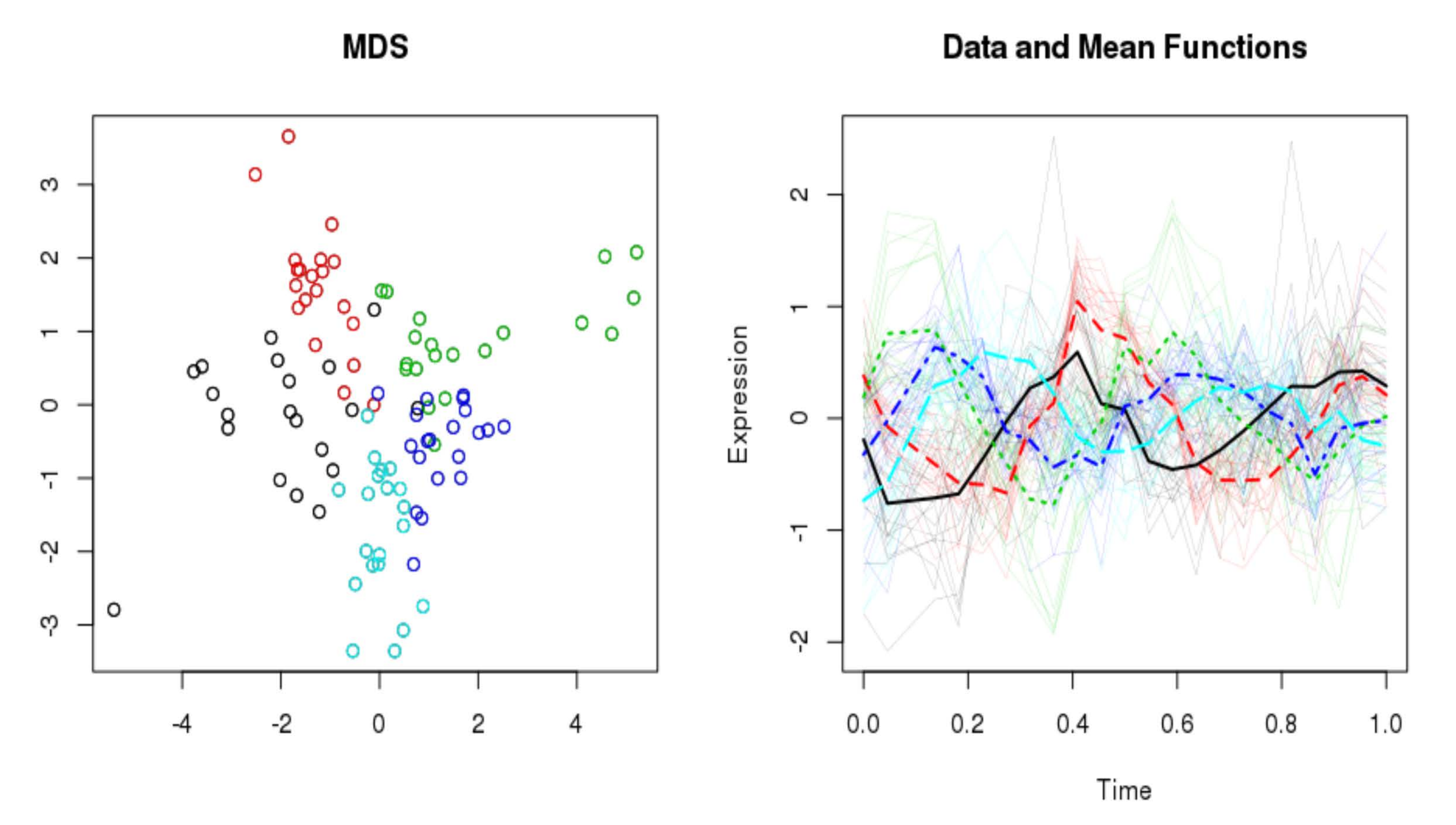} 
         \caption{\small MDS and raw data of Spellman data; right panel: group averages in bold}\label{fig:mds_spellman}
    \end{subfigure}
    \hfill
        \centering
    \begin{subfigure}[t]{.49\textwidth}
        \includegraphics[height = 1.6in, width=\linewidth]{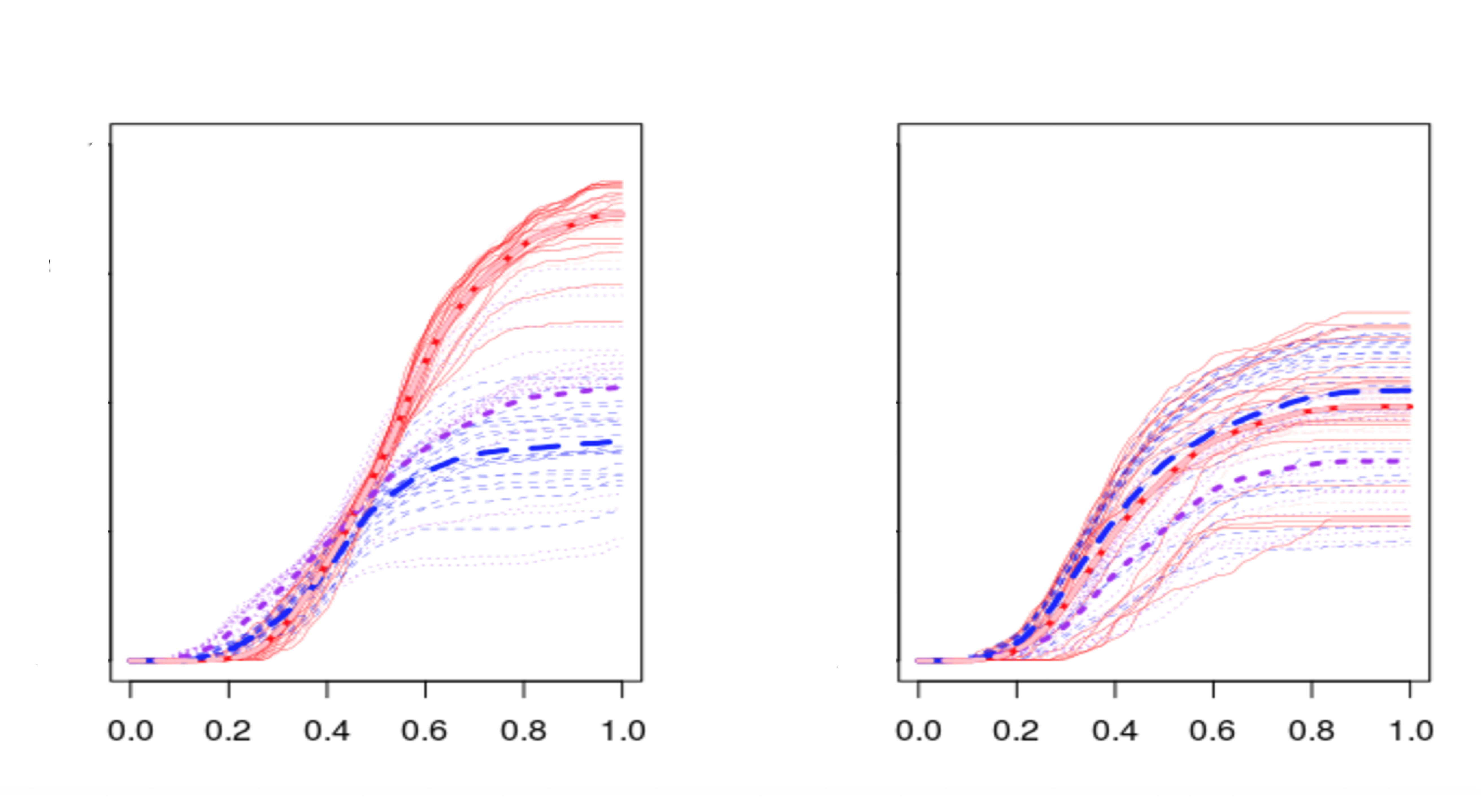} 
         \caption{\small Averaged depth quantile functions of easy and hard classification problems - blue, purple are within class comparisons, thicker lines are group averages.}\label{fig:dq_spellman}
    \end{subfigure}
    \caption{\small From MDS, G1 and M (red vs cyan in MDS plot) are easy to discriminate between, while S and G2 (green vs blue) are more difficult.\\[-24pt]}\label{fig:spellman_full}
\end{figure}
Finally, as an infinite dimensional example, we study Spellman's yeast cell expression data (Spellman et al., 1998), a continuous valued time series. We subsample the data so that only the first $n_i$=20 observations are available for each of 5 classes.  Figure \ref{fig:spellman_full} shows the raw data, as well as a PCA plot of the data. The second panel considers depth quantile functions for two comparisons, one suggested by PCA to be easy (G1 vs M), and a second suggested to be difficult (S vs G2).  The G1 vs M comparison yields a correct classification rate of 97.5\%, while S vs G2 is only 75\%. For comparison, we consider the task of discriminating between G1 and non-G1 phases, as explored in Leng and M\"{u}ller (2006).  We find a classification rate of 94\% compared to theirs of 90\%, though not directly comparable due to different subsetting of the available data.

The method detailed here contains many choices and associated tuning parameters beyond the construction of the depth quantile functions themselves: in particular the number of loadings retained from the fPCA, and the parameters associated with the SVM method (defaults were used above) as well as the averaging of the individual depth quantile functions themselves.  How best to use the information is left as future work.  However, it is clear that interesting information is available, as the results of this ad hoc method tend to be generally competitive with existing methods, while additionally allowing a visual representation of the data.  This visual aspect is particularly beneficial in the next task.        
\subsection{Anomaly Detection}
We conjecture that anchor points associated with outlying observations will tend to live in regions that have different geometry than inlying points.  We consider two data sets to illustrate how depth quantile functions are useful in identifying these anomalous points, strongly aided by the visual nature of this representation of the data.  

\begin{figure}[h]
\centering 
\begin{subfigure}[t]{0.4\textwidth}
\hspace*{0.4cm}\includegraphics[width=6cm,height=4cm]{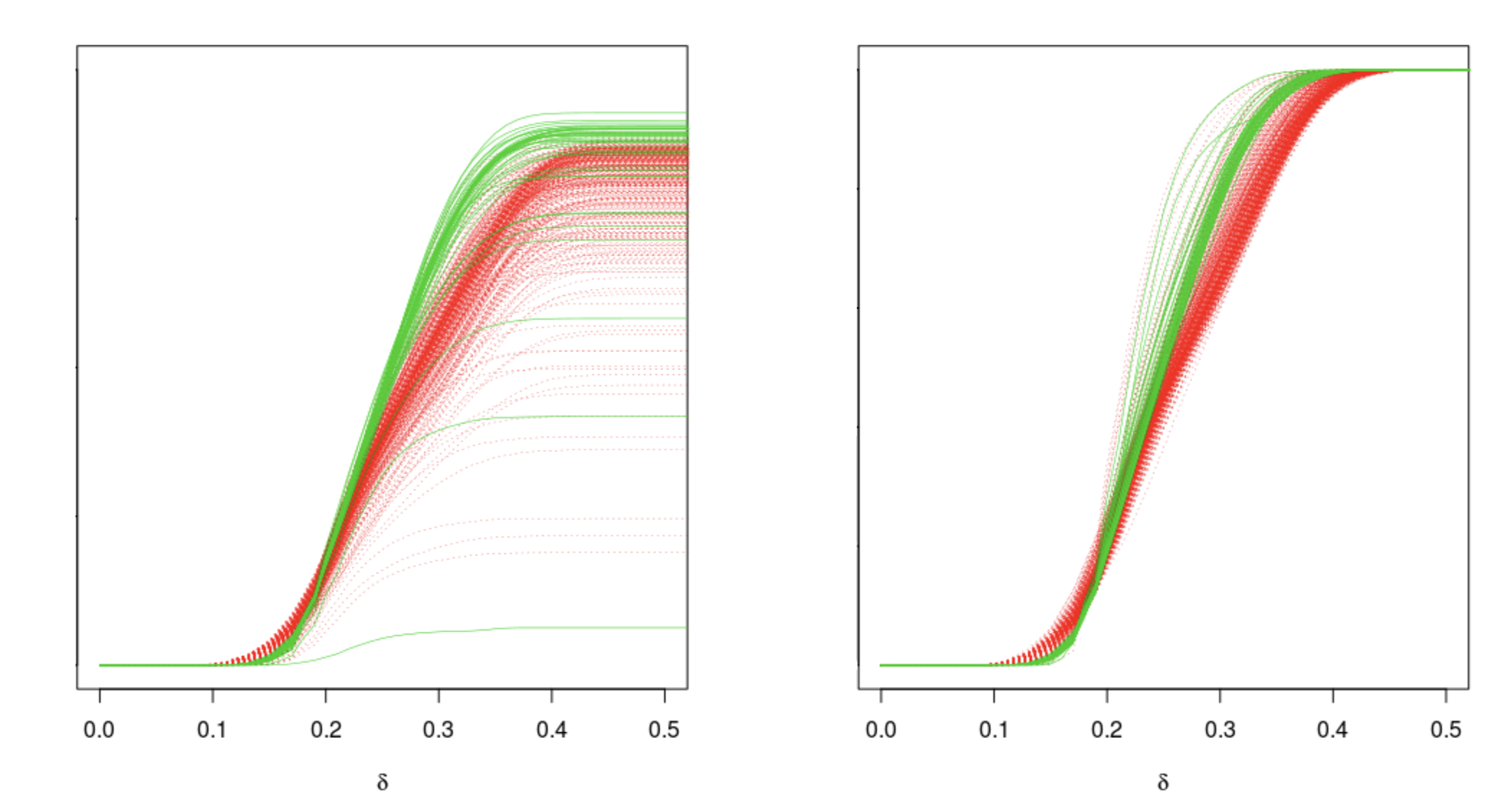}
\end{subfigure}
\begin{subfigure}[t]{0.4\textwidth}
\hspace*{0.2cm}\includegraphics[width=6cm,height=4cm]{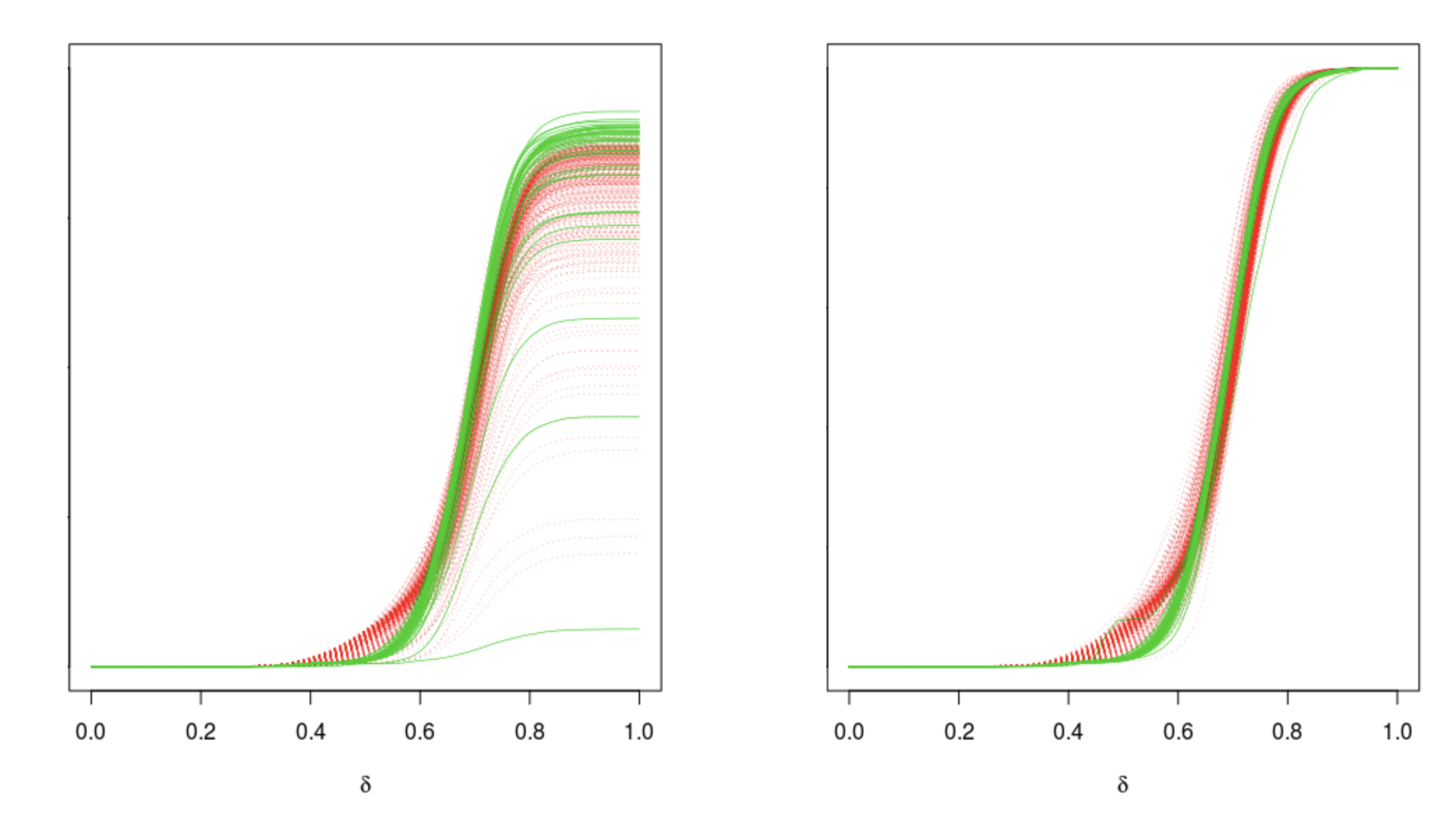}
\end{subfigure}
\caption{\small The left two figures show the averaged depth quantile function and the normalized version for the {\it mfeat} data set with opening angle of $\alpha=150^\circ$.  Green curves correspond to outlying (`0') patterns, red are inlying.  The right two curves show the same as the left, but with an opening angle of $\alpha=90^\circ$.}
\label{fig:mfeat}
\end{figure}

The first data set, the ``multiple features'' data set ({\it mfeat}) considered in Pham (2018) and available on the UCI machine learning repository, consists of $n=440$ observations living in $d=649$ dimensions.  Each dimension constitutes a different feature of handwritten numerals. This data set consists of 200 patterns each of the numbers `6' and `9' and 40 randomly selected `0' patterns.  To illustrate the effect of cone angle, we present the depth quantile functions using both an angle of $150^\circ$  and $90^\circ$ (figure \ref{fig:mfeat}). The similarity of the information between the two angles suggest that this is not necessarily a crucial ingredient.  Note that the location of the increase in the functions happens at different quantiles  depending on angle ($G$ is the same for both).  

Using only a single value of these functions with $\alpha=150^\circ$ on the normalized functions (so only looking at shape), specifically $\wh{q}_i(.17)/\wh{q}_i(1)$, yields an area under the response operator curve (ROC AUC) of .98, besting the results of Pham (2018) of .95, though not directly comparable due to the random subsetting of the data.  

\begin{figure}[h]
\centerline{\includegraphics[width=12cm,height=5cm]{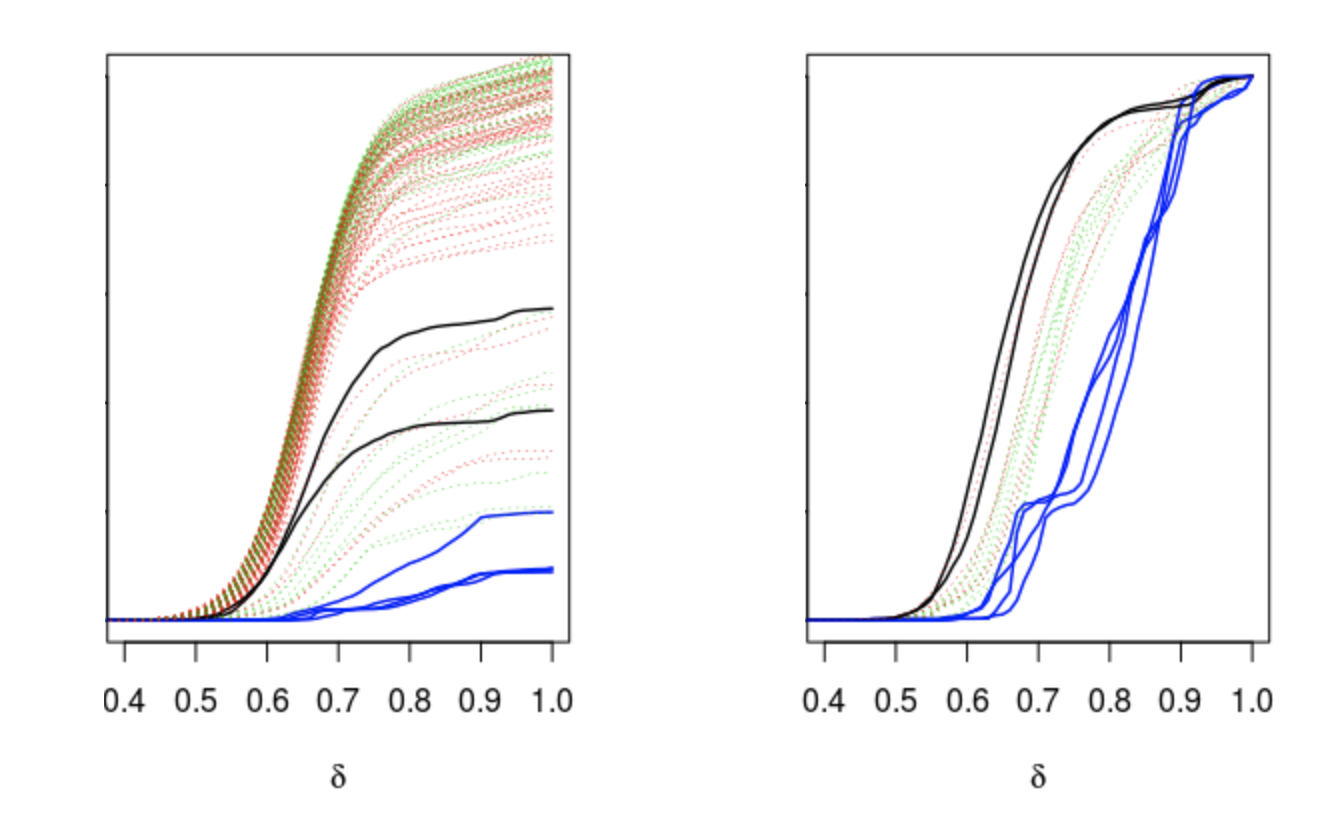}}
\vspace*{-0.2cm}

\caption{{\small  Average depth quantile functions for the lymphography data (left) and normalized version of these functions for the 19 `interesting' observations (right). Bold solid lines represent the 6 total observations from classes 1 and 4 (black and blue respectively).  Dashed lines are for inlying points, color coded for class.}
}
\label{fig:lymph}
\end{figure}
Next, we consider the {\it Lymphography} data set, also available at the UCI repository.  The data set consists of 4 classes over 148 observations in 18 dimensions. Classes 1 and 4 only constitute 6 instances and are treated as outliers. 

Figure \ref{fig:lymph} shows the analysis of this data using an angle of $\pi/2$.  The depth quantile functions suggest there are 19 observations (all 6 outliers and 13 inlying points) that warrant further inspection.  On inspection of the normalized functions, the outlying functions distinguish themselves from all but one of the inlying functions.  Interestingly, the variables in this data set are all numeric assignments of categorical variables, though clearly exploration of the geometry of the induced point cloud is still valuable. This seems to be consistent with Campos et al. (2016), who found an ROC AUC of 1 is attainable. 

\section{Other conceptual relations}
\label{relations}

\subsection{Interpreting average depth quantiles via random set theory}
\label{random-sets}

In an attempt to obtain further insight about the average depth quantiles, we use random set theory to provide a different point of view to the construction of our feature functions. For simplicity we fix $(ij) = (12).$ 
%
%

Recall that $A_{12}(s) = A_{m_{12},u_{12}}(s), B_{12}(s)= B_{m_{12},u_{12}}(s)$ with $m_{12} = \frac{X_1 + X_2}{2}$ and $u_{12} = \frac{X_1 - X_2}{\|X_1 - X_2\|}.$  Let 
%
$d_{12} (s) = \min\big(F(A_{12} (s)), F(B_{12}(s))\big) = F(C^*_{12}(s)),$
%
with $C^*_{12}(s)\in \{A_{12}(s), B_{12}(s)\}$. The joint distribution of $X_1, X_2$ induces the distribution of the random set $C^*_{12}(s),$ and thus of $d_{12}(s),$ the probability content of the set. 
We can write
$$d_{12}(s)= P_Z(Z \in C^*_{12}(s)|X_1,X_2) ,$$
where $Z\sim F$ independent of $X_1,X_2.$ The corresponding random depth quantile function is
%
$q_{12}(\delta) =  d_{12}(s_{12}(\delta)),$ 
%
where $s_{12}(\delta) = s_{12}^l(\delta)$ is as in Lemma~\ref{multiscale} but with $(x,u)$ replaced by $(m_{12},u_{12})$ (equivalently, we could use $s_{12}^r(\delta)$ instead of $s_{12}^l(\delta)$). Mimicking the construction of the average empirical feature functions used in the classification context (see section~\ref{classification}), we obtain by conditioning on $X_1$ and taking expectation over $X_2$ that
\begin{align*}
E_{X_2}\big[q_{1,2}(\delta)\big| X_1 = x_1\big] 
&= E_{X_2}P_Z\big( Z \in C^*_{12}(s_{12}(\delta)) \big) |X_1 = x_1\big] 
= E_Z \Psi_{x_1, \delta}(Z),
\end{align*}
with  $\Psi_{x_1,\delta}(z) = P_{X_2}\big(z \in C^*_{12}(s_{12}(\delta))\big| X_1 = x_1\big)$ 
the hitting function of the random set $C^*_{12}(s_{12}(\delta))$ conditional on $X_1 = x_1$. This random function in $\delta$ depends on the distributions of both $X_1$ and $X_2.$ In the context of binary classification (with a symmetric $G$), this gives three cases corresponding to three different comparison functions $\delta \to E_Z\Psi_{X_1,\delta}(Z)$: one between and two within class comparisons. Note that our empirical average depth quantile functions $\wh q_j^k(\delta)$ with $y_j, y_k \in \{0,1\}$ (see section~\ref{classification}) can be considered as samples from the (estimated) distributions of $E_Z \Psi_{X_1, \delta}(Z)$ for different combinations of class labels of $X_1,X_2.$ In order for our classification procedure to work, these distributions of the random functions need to differ. Notice that a similar interpretation holds for $\wh q_{12}(\delta)$  by replacing $P_Z$ by the empirical distribution of the sample $X_3, X_4,\ldots,X_n$, and the distribution of $X_2$ by the empirical distribution of the sample. 


\subsection{Depth and hitting functions in the literature}\label{conceptual}

Hitting functions have appeared it the context of data depth before, but in a different fashion. For instance, type A depth functions are defined in Zuo and Serfling (2000) as being of the form $d(x,F) = Eh(x;X_1,\ldots,X_n)$ for some function $h$, where $X_i \sim F, i = 1,\ldots,n.$ When $h(x;X_1,\ldots,X_n) = {\bf 1}(x \in S(X_1,\ldots,X_n)),$ for some (random) set $S = S(X_1,\ldots,X_n),$ we obtain the depth function as a hitting function $P(x \in S(X_1,\ldots,X_n)).$ The perhaps best known such depth measure is the simplicial depth introduced by Liu (1988, 1990), where $S(X_1,\ldots,X_n) = {\rm conv}[X_1,\ldots,X_{d+1}]$ with ${\rm conv}[x_1,\ldots,x_{d+1}]$ denoting the convex hull of the points $x_k, k = 1,\ldots,d+1.$  
Another instance is the depth based on pairwise connecting spheres (see Elmore et al., 2006) for data on a sphere, and recently extended to data on Riemannian manifolds in Fraiman et al. (2019). Here, given a Riemannian manifold ${\cal M}$ with geodesic distance $d_g$, a point $x \in {\cal M}$ is assigned the depth $d(x) = P(x \in B_{X_1,X_2})$, where $X_1,X_2$ are iid, and $B_{x_1,x_2}$ denotes the ball with radius $d_g(x_1,x_2),$ centered at the midpoint between $x_1$ and $x_2.$   The work of Ting et al. (2013) is based on a similar premise. Considering data in $\R^d,$ the authors propose to form binary partitions by using random, axis parallel splits. For a given point $x \in \R^d,$ the average mass of the rectangles containing $x,$ taken over a number of binary partitions is considered. It also defines some type of depth function, but the average of the probability contents over a random set is considered, rather than counting the number of boxes containing $x$. As the discussion above shows, in our work, expected hitting functions come into play when averaging depth quantile functions, and we in fact consider {\em expected} hitting functions for an entire class random sets, parametrized by $\delta$. 

\subsection{Relations to local depth}\label{local-depth} Notions of local depth by construction give rise to  `multiscale' approaches (e.g. Serfling, 2019). Given any such notion of local depth corresponding a scale parameter $\tau > 0,$ and given a point $x,$ one can consider a depth function $\tau \to d(x,\tau)$. If the scale parameter were chosen according to a distribution $G$, then this would also allow construction of a depth quantile function corresponding to these local depths similar to what is discussed above. However, one would obtain only one depth quantile function for a given $x.$ In our set-up, however, the notion of depth also depends on the direction, leading to more refined information. For instance, if one were to choose a data point $X_i$ as an anchor point (cf. see discussion on choice of base point given below), then our construction would result in $n-1$ quantile functions, each describing a different distribution of depths for $X_i.$ It might be interesting to study differences in the distributions of these $n-1$ functions over varying $X_i$'s.  Also, not all the existing notions of local depth lead to localized information (density) for small scales, as, for instance, the multivariate local simplical depth, or the local halfspace depth in Agostinelli and Romanazzi (2011) and Agostinelli (2018), respectively. 

\section{Discussion of choice of tuning parameters}
\label{tuning}

The proposed procedure has three ingredients that can be considered as tuning parameters. The most obvious one is the angle $\alpha$ of the cones, then there is the distribution $G$ under which to choose the distance of the tip of the cone to the base point, and last, but not least, there is the choice of the base point. We briefly discuss some aspects of the choice of these ingredients. Moreover, one might also consider our choice of (spherical) cones in the construction of our feature functions to be a `tuning parameter', and we give some discussion on possible alternative choices below.

{\it Choice of base point:} In this paper, we only consider anchor points of the form $m_{ij} = \frac{X_i + X_j}{2}$ in defining $\wh q_{ij}(\delta)$. It may seem surprising that we choose to describe an observation not by the geometry associated with the observation itself but rather with its associated midpoints. A motivation of this choice in the classification case has been given in the application section, though we use it for all statistical tasks.  The rational is that using our depth-based approach with high dimensional point clouds necessitates consideration of values other than the actual observations.  It is well known (see R\'{e}nyi and Sulanke, 1963) that in high dimension, many points will live on or near the convex hull.  For such points $x$, one of the two sets $A_{x,u}(s)$ or $B_{x,u}(s)$ will be empty for many if not all values of $s$, leading to degenerate, non-informative depth quantile functions. Empirically, we have found this to be the case, with many estimated depth quantile functions being identically zero when considering observations as the anchor points.  Considering anchor points that live inside the convex hull, which midpoints often are, helps in alleviating this problem. While other choices might be reasonable as well, we have not explored them here. 

{\it Choice of opening angle $\alpha$.} As a general guideline, a reasonable choice of $\alpha$ will be such that not all the feature functions will increase very steeply at about the same value of $\delta,$ because this will result in weak discrimination. For $G$ fixed, a too large value of $\alpha$ will result in a steep increase early (small $\delta)$, while for too small values of $\alpha$ this will happen late. In practice we observed that the shape of the feature functions is not too sensitive to the choice of $\alpha$, adjusting the support of $G$ accordingly. In our simulation studies the choice of an angle of $90^\circ$ seems to work well for a large range of dimensions.

Still, one would expect that the opening angle of the cones should depend on the dimension $d$. In fact, as $d$ increases, the opening angle needs to tend to $\pi$ to avoid the angle being `too small'. This is due to the well-known fact that for a fixed opening angle (and tip in 0, say), the portion of the Euclidean space contained in the cone will tend to zero as the dimension increases to infinity. 
For non-Euclidean data the choice of $\alpha$ appears more challenging, because of some lack of geometric intuition. 

As suggested by one of the referees, the angle could also be chosen randomly, similar to the tip $s$. This would mean that one has to choose a two dimensional distribution on $(s,\alpha).$ The resulting depth quantile function would still be a one-dimensional function (because the depth is a one-dimensional quantity), but computationally this would be more challenging, for the depth functions would now be a function of two variables, meaning that its level sets are no longer (one-dimensional) intervals. This will complicate the computation. We will  explore this further in future research. 

{\it Choice of the distribution $G$.} Choosing $G$ essentially means choosing a non-linear transformation to the horizontal axis, amounting to emphasizing different regions of the depth quantile plots. This is somewhat convoluted to the choice of the opening angle $\alpha$ (see above). Choosing $G$ to be unimodal about $0,$ for instance, would put higher weight on points close to the base-points, and would thus lead to depth quantile functions that are flatter for small delta than for a uniform $G$, for instance. We recommend to choose $G$ to be symmetric about zero because of the symmetry property of the quantile functions given in Lemma~\ref{multiscale} part (b). In our applications, we always chose $G$ to be the uniform distribution on an interval chosen such that all the projections used to determine the empirical depth quantiles fall into the support of $G$. The motivation for this choice is similar to the choice of a non-informative prior. 

{\em Choice of geometric shapes.} Using spherical cones in our procedure is a choice of convenience. Choosing elliptical cones (with the ellipse perhaps chosen adaptively) appears to be a first option. However, choosing an ellipse that works well for all directions $u$ does not appear to be straight-forward. Our current choice of spherically symmetric cones does imply no preference for a direction. Geometric objects with different shapes, such as parabolas or hyperbolas, can also be used.  See for instance Kotic and Hlubinka (2017) for a discussion on using such objects in the context of weighted depth. These choices would be similar to cones in that they also provide some sense of direction. In how much the choice of the shapes will effect the performance of the methodology needs to be explored.

\section{Properties of depth functions and depth quantile functions}\label{theory}

Here we study various properties of the appropriately normalized quantities $\widehat{d}_{x,u}(s), \wh q_{x,u}(\delta),$ $\wh d_{ij}(s)$, and the averaged feature functions $\frac{1}{n}\sum_{j=1, j \ne i}^n \wh q_{ij}(\delta)$ considered as processes in $s$ and $\delta$, respectively. On the one hand, the presented results are meant to highlight some interesting properties of these quantities related to the curse of dimensionality and adaptivity, and on the other hand they provide the basis for the formal construction of inference procedures. We do not present such constructions here. The main emphasis of this work is to introduce the basic approach and to illustrate its versatility. 
%
%
%
%
%
 
\subsection{Concentration and large sample behavior} 
\label{large sample}

Recall that $\ell_{x,u} = \{ y \in \R^d: y = x + su,\,s \in \R \} $.  For $C,D \subset \R^d$ measurable, define the pseudometric $d_F(C,D) = F(C \Delta D)$, where $C \Delta D = (C\setminus D) \cup (D \setminus C).$  
We let 
$$\C_{x_0,u_0} = \big\{ C_{x,u_0}(s),\, x \in \ell_{x_0,u_0}, s \in \R\big\}$$
be the set of all cones (with opening angle $\alpha$) with $\ell_{x_0,u_0}$ as their axis of symmetry. With $\H_{u}$ the class of all halfspaces with normal vector $u$, $\D_{x_0,u_0} := \C_{x_0,u_0} \cap \H_{u_0}$ equals the set $\{A_{x,u_0}(s), B_{x,u_0}(s),\, x \in \ell_{x_0,u_0}, s \in \R\}.$  Both $\C_{x_0,u_0}$ and $\H_{u_0}$ are VC-classes with VC-dimension $2$ (both are unions of two nested classes of sets), and thus $\D_{x_0,u_0}$ is a VC-class with VC-dimension bounded by $4$. The importance here is that the VC-dimension of $\D_{x_0,u_0}$  does not vary with the dimension, which is the reason for the estimators of the depth quantile functions not suffering from the curse of dimensionality, as is shown in the following result. Furthermore, if we let
$${\textstyle \D = \bigcup_{(x,u) \in \R \times S^{d-1}}\D_{x,u}.}$$
then, with $\C = \bigcup_{x\in \R,u \in S^{d-1}}\C_{x,u}$ and $\H = \bigcup_{u \in S^{d-1}} \H_u,$ we have $\D = \{C \cap H,\,C \in \C, H \in \H\}.$ $\D$ is a VC-class, since both $\C$ and $\H$ are VC-classes. However its VC-dimension depends on the dimension $d$. In fact, it is $O(d)$.

\subsubsection{Results for empirical depth functions} 
\begin{proposition}\label{stoch-bound}
For every given line $\ell_{x_0,u_0},$ there exists a constant $c > 0$ not depending on the dimension $d$ and the opening angle $\alpha$ such that, for all $M > 0,$ 
\begin{align}\label{depthity-bound}
P\big[ \sup_{x \in \ell_{x_0,u_0}, s \in \R}\big|\sqrt{n}(\wh d_{x,u_0}(s) - d_{x,u_0}(s))\big| > M \big]  \le c M^3 \exp\big\{-\tfrac{M^2}{2}\big\}.
\end{align}
\end{proposition}

Now recall that our methodology depends on the $\binom{n}{2}$ functions $\wh d_{ij}(s), i < j.$ We show that, with $A_{ij}(s) = A_{m_{ij}, u_{ij}}(s)$ and $B_{ij}(s) = B_{m_{ij}, u_{ij}}(s),$ these functions are uniformly close to the following population based functions
$$d_{ij}(s) = \min\big(F\big(A_{ij}(s)\big), F\big(B_{ij}(s)\big)\big).$$
While $d_{ij}(s)$ depends on the population distribution, it is still a random quantity, but only through the choice of the line and the base point. We now have the following result:
\begin{proposition}\label{uniform-bound} With $d = d_n \in {\mathbb N}$ and $\alpha = \alpha_n \in [0,\pi/2)$, both nonrandom sequences, we have, 
\begin{align}\label{uniform-depthity-bound}
\max_{i,j = 1,\ldots,n\atop i \ne j} \sup_{s \in \R}\big|\wh d_{ij}(s) - d_{ij}(s)\big| = O_P\big(\sqrt{\tfrac{\min(d,\log n)}{n}}\;  \big)\qquad \text{as }\,  n \to \infty.
\end{align}
Moreover, if the distribution $F$ is concentrated on a $d^{\,\prime}$-dimensional hyperplane in $\R^d$, then (\ref{uniform-depthity-bound}) holds with $d$ replaced by $d^{\,\prime}.$
\end{proposition}
{\sc Remark.} The `min' appearing in the rate tells us that $\big(\frac{\log n}{n}\big)^{1/2}$ is a dimension independent uniform bound for the rate of the deviations of the functions $\wh d_{ij}(s)$ from their population based counterparts $d_{ij}(s).$ The asserted adaptivity to some kind of `sparsity' is due to some geometry: Consider a cone in $\R^d$ with tip at zero and with axis of symmetry being one of the coordinate axes. The intersection of such a cone with a linear subspace of dimension $d^\prime$ is a cone in $\R^{d^\prime}$.\\[-8pt]

\noindent
Our next result concerns the asymptotic distribution of $\sqrt{n}\big(\wh d_{x,u}(s) - d_{x,u}(s)\big)$ as a process in $s \in \R$. The pointwise behavior of the process depends on whether the minimizing set in the definition of $d_{x,u}(s)$ is unique or not (see Corollary~\ref{pointwise}). Therefore, we are facing a challenge in small neighborhoods of points with non-unique minimizers, because the empirical minimizer, i.e. the sets minimizing $\wh d_{x,u}(s)$, might not stabilize in such neighborhoods. Moreover, suppose that $s_0$ is a point with $d_{x,u}(s_0) = F(A_{x,u}(s_0)) = F(B_{x,u}(s_0)),$ and consider $s$ in a small neighborhood of $s_0$. If the minimizers corresponding to $s < s_0$ are $A_{x,u}(s),$ and for $s^\prime > s_0$ they switch to $B_{x,u}(s^\prime)$, then the map $s \to {\bf 1}_{\Cstar_{x,u}(s)}$ is not continuous at $s_0$ (w.r.t. the $L_1(F)$-norm), because the sets $A_{x,u}(s)$ and $B_{x,u}(s^\prime)$ are not nested. 

The following definition is to exclude such values $s_0$. For $(x,u) \in \R^d \times S^{d-1},$ denote $\Delta_{x,u}(s) = \big| F(A_{x,u}(s)) - F(B_{x,u}(s)) \big|,$ and for $\eta > 0$ let
\begin{align}\label{Sxell}
S_{x,u}(\eta) = \big\{s \in \R: \Delta_{x,u}(s) \ge \eta\big\}.
\end{align}
Note that $\Delta_{x,u}(s) = 0$ means that the minimizer $\Cstar_{x,u}(s)$ in the definition of $d_{x,u}(s)$ is not unique (both $A_{x,u}(s)$ and $B_{x,u}(s)$ are minimizers). 
We also introduce the notation
$$\T_{x,u}(s) = \argmin_{C \in \{A_{x,u}(s), B_{x,u}(s)\}} F(C).$$
Observe that this set of minimizers $\T_{x,u}(s)$ either consists of exactly one, or of exactly two elements. The latter holds if and only if $\Delta_{x,u}(s) = 0.$ 
\begin{theorem} \label{depth}
For $(x,u) \in \R \times S^{d-1}$ set $\Sstar_{x,u}(\eta) = S_{x,u}(\eta) \cup \{s: \Delta_{x,u}(s) = 0\}$, and let $\{\eta_n\}$ be a sequence of real numbers with $\eta_n \sqrt{n}\to \infty$. Then, on an appropriate probability space, there exists a sequence of $F$-bridges $\{\B^{(n)}_{F}(C), C \in \D\}$, i.e. tight, mean zero Gaussian processes with  covariance ${\rm Cov}(\B^{(n)}_{F}(C),\B^{(n)}_{F}(D)) = F(C\cap D) - F(C)F(D),$ and almost surely continuous sample paths w.r.t. $d_F$, such that, for any given line $\ell_{x,u}$, we have with $\B^{(n)}_{x,u}(s) = \min\{\B^{(n)}_F(C),\,C \in {\cal T}_{x,u}(s)\}$ that
\begin{align*}
\sup_{s \in \Sstarsmall_{x,u}(\eta_n)}\big|\sqrt{n}(\wh d_{x,u}(s) - d_{x,u}(s)) - \B^{(n)}_{x,u}(s)\big| = o_P(1),\qquad \text{as } n \to \infty,
\end{align*}
where the $o_P(1)$-term 
does not depend on $d, \alpha$ and $u$. Moreover, we have 
\begin{align*}
\sup_{(x,u) \in \R \times S^{d-1}} \sup_{s \in \Sstarsmall_{x,u}(\eta_n)}\big|\sqrt{n}(\wh d_{x,u}(s) - d_{x,u}(s)) - \B^{(n)}_{x,u}(s)\big| = o_P(1),\;\; \text{as } n \to \infty,
\end{align*}
where the $o_P(1)$-term does not depend on $\alpha$, but it does depend on $d$.
\end{theorem}
{\sc Remark.} In the somewhat related context of minimum volume sets (or generalized quantile processes), a result of similar type has been obtained in Einmahl and Mason (1992). There the approximating process (to the generalized quantile process) is a maximum of $F$-bridges taken over the, in general non-unique, generalized quantiles. Also, Mass\'e (2004) showed convergence of the multivariate Tukey depth process to a similar type of limit process. 

\begin{corollary} \label{pointwise}
For $(x,u) \in \R^d \times S^{d-1},$ set $\sigma_{x,u}^2(s) = d_{x,u}(s) (1 - d_{x,u}(s)), s \in \R$. \\[5pt]
{\bf (i)} For $(x,u,s) \in \R^d \times S^{d-1} \times \R$ with $\Delta_{x,u}(s) \ne 0$, we have
$$\sqrt{n}\big(\wh d_{x,u}(s) - d_{x,u}(s)\big) \to {\cal N}(0, \sigma_{x,u}^2(s))\qquad\text{in distribution as $n \to \infty.$}$$
{\bf (ii)} For $(x,u,s) \in \R^d \times S^{d-1} \times \R$ with $\Delta_{x,u}(s) = 0$, we have
$$\sqrt{n}\big(\wh d_{x,u}(s) - d_{x,u}(s)\big) \to \min(N_1,N_2)\qquad\text{in distribution as $n \to \infty,$}$$
where $N_1,N_2$ are jointly normal with mean zero, variance $\sigma_{x,u}^2(s)$, and covariance Cov$(N_1,N_2) = - d^{\,2}_{x,u}(s).$
\end{corollary}

\subsubsection{Results for empirical depth quantile functions} First we give some finite sample concentration results.
\begin{proposition}\label{consistency-quant} For every fixed line $\ell_{x_0,u_0} \subset \R^d,$ there exists a constant $c > 0$, not depending on $d$ and $\alpha$, such that, for all $M > 0,$ 
\begin{align}\label{depthity-quant-bound}
P\big[\sqrt{n}\,\sup_{x \in \ell_{x_0,u_0},\, \delta \in [0,1]}\big| \wh q_{x,u_0}(\delta) - q_{x,u_0}(\delta) \big| > M\big] \le c M^3 e^{-M^2/2}.
\end{align}
\end{proposition}

Our next result addresses the asymptotic distribution of the depth quantile functions. We need the following assumptions. \\[5pt]
{\bf (A1)} {\em The distributions $F$ and $G$ have bounded, continuously differentiable Lebesgue densities $f > 0$ on $\R^d$, and $g > 0$ on $\R$, respectively.} \\[5pt]
{\bf (A2)} {\em (tail behavior of $f$) There exist constants $C, \epsilon > 0$ such that }
\begin{align*}
{\bf (i)}&\;\;\sup_{z \in \R^d} \|z\|^{d-1 + \epsilon}f(z) < C < \infty\;\,\text{\em and}\\
{\bf (ii)}&\;\;\sup_{z \in \R^d} \|z\|^{d-1 + \epsilon}\|{\rm grad}f(z)\| < C < \infty.
\end{align*}
{\bf (A3)} {\em For each $K > 0,$ there exists $\epsilon_K > 0$ such that 
$$\inf_{\|x\| \le K,\, u \in S^{d-1}}\inf_{s: \Delta_{x,u}(s)=0} \big|\tfrac{d}{ds}F(A_{x,u}(s)) - \tfrac{d}{ds}F(B_{x,u}(s))\big| > \epsilon_K,$$
where {\bf (A1)} and {\bf (A2)} imply that these derivatives exist (see Lemma~\ref{derivative} in the supplemental material).\\}

{\sc Discussion of {\bf (A3)}.} This assumption is used to show that, for each $(x,u)$ and $\eta$, the set  $\mycomp{S}_{x,u}(\eta)$ consists of a (finite) union of `small' intervals, and that  the same is true for is the set $\mycomp{D}_{x,u}(\eta,\epsilon)$; cf. Lemma~\ref{short-intervals}. Intuitively, assumption {\bf (A3)} can be expected to hold in many situations. The reason for this intuition is that $\mycomp{S}_{x,u}(\eta)$ can only be large, if  $F(A_{x,u}(s))$ and $F(B_{x,u}(s))$ have the same derivatives (w.r.t. $s$) in a neighborhood of a point $s_0$ with $F(A_{x,u}(s_0)) = F(B_{x,u}(s_0)).$ These derivatives are surface integrals (with respect to $f$) over the `lateral' surface areas of the cone $A_{x,u}(s)$ and the frustum $B_{x,u}(s)$, respectively (see proof of Lemma~\ref{derivative}). Given the geometry of these sets, these derivatives being equal over a (small) interval, puts some strong requirements on the geometry of the density $f.$ For more standard density models, these requirements do not seem to hold, at least intuitively, while explicit calculations appear challenging, even for a normal model, say.  Also, suppose that while $F(A_{x,u}(s_0)) = F(B_{x,u}(s_0)),$ the derivative of these two functions at $s_0$ is different, let's say, $\frac{d}{ds} F(A_{x,u}(s_0)) > \frac{d}{ds} F(B_{x,u}(s_0)).$  Then, in order for $F(A_{x,u}(s_1)) = F(B_{x,u}(s_1))$ to hold for some $s_1 > s_0$, the inequality between the derivatives has to change along the way. Again, this puts some requirements on the geometry of the underlying density. Intuitively, these requirements are more likely to hold, if $f$ has several modes, allowing for a less regular  behavior of these derivatives, by which we mean that they have a large number of local extrema. This leads to the intuition that $\mycomp{S}_{x,u}(\eta)$ (and $D_{x,u}(\delta)$) can be expected to consist of a finite number of small intervals assuming that the number of modes of $f$ is finite. 

Recall the notation $\{s: d_{x,u}(s) \le t\} = [s^l_{x,u}(t), s^r_{x,u}(t)]$, and that for $t = q_{x,u}(\delta),\, 0 \le \delta < 1$, we use the short-hand notation
$$\{s: d_{x,u}(s) \le q_{x,u}(\delta)\} = [s^l_{x,u}(\delta), s^r_{x,u}(\delta)].$$ 
Also recall the definition of the set $S_{x,u}(c)$ given in (\ref{Sxell}), and that $\Cstar_{x,u}(s)\in \{A_{x,u}(s),B_{x,u}(s)\}$
denotes the minimizer in the definition of $d_{x,u}(s)$. We now have the following result. 

\begin{theorem}\label{as-distr-quant} Suppose that assumptions {\bf (A1)} - {\bf (A3)} hold. For $C, \epsilon, \eta > 0$, let $D_{x,u}(\eta,\epsilon) = \{\delta \in [\epsilon,1-\epsilon]: s^r_{x,u}(\delta), s^l_{x,u}(\delta) \in S_{x,u}(2\eta)\}.$ 
For each $K > 0$, and for real numbers $\eta_n > 0$ with $\eta_n \to 0$ and $\eta_n \sqrt{n} \to \infty$, we have, as $n \to \infty$, \\[-15pt]
\begin{align*}
\sup_{u \in S^{d-1} \atop \|x\| \le K} \sup_{\delta \in D_{x,u}(\eta_n,\epsilon)}&\Big| \sqrt{n}\big(\wh q_{x,u}(\delta) - q_{x,u}(\delta) \big) \\[-5pt]
&\hspace*{-0.5cm}- \big[\alpha_{x,u}(\delta) \nu_n(s^l_{x,u}(\delta))- (1-\alpha_{x,u}(\delta))\,\nu_n(s^r_{x,u}(\delta))\big]\Big| = o_P\big(1\,\big),
\end{align*}
where   $\alpha_{x,u}(\delta) =  \big(1 +\frac{g(s_{x,u}^l(\delta))}{g(s_{x,u}^r(\delta))}\cdot \big|\frac{d^\prime_{x,u}(s_{x,u}^r(\delta))}{d^\prime_{x,u}(s_{x,u}^l(\delta))}\big| \big)^{-1},$ and $\nu_n(s^{l(r)}_{x,u}(\delta)) = \sqrt{n}(F_n - F)(\Cstar_{x,u}(s^{l(r)}_{x,u}(\delta))).$ 
Moreover, on an appropriate probability space, there exists a sequence of mean zero $F$-bridges $\{\B^{(n)}_{F}(C), C \in \D\}$ 
, such that with $\B^{(n)}_{x,u}(s) = \B_F^{(n)}(\Cstar_{x,u}(s))$, we have, as $n \to \infty$,
\begin{align*}
\sup_{u \in S^{d-1}\atop \|x\| \le K}\sup_{\delta \in D_{x,u}(\eta_n,\epsilon)}&\Big| \sqrt{n}\big(\wh q_{x,u}(\delta) - q_{x,u}(\delta) \big) \\
&\hspace*{-0.8cm}- \big[\alpha_{x,u}(\delta) \B_{x,u}^{(n)}(s^l_{x,u}(\delta)) - (1-\alpha_{x,u}(\delta))\,\B_{x,u}^{(n)}(s^r_{x,u}(\delta))\big]\Big| = o_P\big(1\,\big).
\end{align*}
\end{theorem}

{\sc Discussion:} (i) As $\sqrt{n}\big(\wh q_{x,u}(\delta) - q_{x,u}(\delta) \big)$ is a quantile-type process, and the approximating process is a (mixture of two) empirical processes given by $\nu_n(s^{\,l}_{x,u}(\delta))$ and $\nu_n(s^{\,l}_{x,u}(\delta))$, respectively, the first assertion can be seen as a generalized Bahadur-Kiefer representation. \\[3pt]
(ii) The approximating Gaussian processes $\B^{(n)}_{x,u}(s^{\,l (r)}_{x,u}(\delta))$ depend on $\Cstar_{x,u}(s)$, which is either $A_{x,u}(s)$ or $B_{x,u}(s).$ In other words, either $\B_{x,u}^{(n)}(s^{\,l(r)}_{x,u}(\delta)) =  \B_F^{(n)}(A_{x,u}(s^{\,l(r)}_{x,u}(\delta)))$ or $\B_{x,u}^{(n)}(s^{\,l(r)}_{x,u}(\delta)) =  \B_F^{(n)}(B_{x,u}(s^{\,l(r)}_{x,u}(\delta))).$ Both the processes $\{\B_F^{(n)}(A_{x,u}(s^{\,l(r)}_{x,u}(\delta))), \delta \in [0,1]\}$ are time-transformed one-dimensional Brownian bridges, because both $\{A_{x,u}(s), s \in \R\}$ and $\{B_{x,u}(s), s \in \R\}$ consists of nested classes of sets, and because the functions $\delta \to s^{\,l(r)}_{x,u}(\delta)$ are monotonic. As $\delta$ varies, we can think of $\B_{x,u}^{(n)}(s^{l(r)}_{x,u}(\delta))$ as switching between these two processes, depending on which of the two sets $A_{x,u}(s)$ or $B_{x,u}(s)$ is the corresponding minimizer of the depth functions. The approximation works for $\delta \in D_{x,u}(\eta_n,\epsilon)$, because we know that the empirical minimizer and the theoretical minimizers coincide with probability tending to one. By our assumption, the set $D_{x,u}(\eta_n,\epsilon)$ is a finite union of intervals, separated by small neighborhoods about points $s_0$ with $\Delta_{x,u}(s_0) = 0$, so that there is no switching between approximating processes within a subinterval.\\[3pt]
(iii) The mixing proportion $\alpha _{x,u }(\delta)$ is determined by the two quantities\\ $g(s^r_{x,u }(\delta))/d_x^{\,\prime}(s^r_{x,u }(\delta))$ and $g(s^l_{x,u }(\delta))/d_x^{\,\prime}(s^l_{x,u }(\delta))$, which are the derivatives of $G(s \ge 0: d_{x,u}(s) \le t)$ and $G(s < 0: d_{x,u}(s) \le t)$, respectively,  evaluated at $t = q_{x,u}(\delta)$. The restrictions on $x$ and $\delta$ given in the statements of the theorem make sure that these ratios are uniformly bounded away from zero and infinity, and that they are uniformly continuous maps as $x$ and $\delta$ vary on the restricted sets. This is needed to assure that the asserted uniformity of the asymptotic statements.


\begin{corollary} Under the assumptions of Theorem~\ref{as-distr-quant}, we have for each $(x,u) \in \R^d \times S^{d-1}$ and $\delta \in [0,1)$ with $\Delta(s_{x,u}^r(\delta)) \ne 0$ and $\Delta(s_{x,u}^l(\delta)) \ne 0$, we have, as $n \to \infty$,
$$\sqrt{n}\big(\wh q_{x,u}(\delta) - q_{x,u}(\delta) \big) \to \alpha_{x,u}(\delta) N_{x,u}^l(\delta) + \big((1-\alpha_{x,u}(\delta)\big) N^r_{x,u}(\delta)\;\;\text{in distribution,}$$
where $N^l_{x,u}(\delta),N^r_{x,u}(\delta)$ are mean zero, jointly normal random variables with variances $F\big(\Cstar_{x,u}(s^l_{x,u}(\delta)) \big)$ $\big[1- F\big(\Cstar_{x,u}(s^l_{x,u}(\delta))\big) \big],$ and $F\big(\Cstar_{x,u}(s^r_{x,u}(\delta))\big) \big[1- F\big(\Cstar_{x,u}(s^r_{x,u}(\delta))\big) \big],$ respectively. Their convariance is given by ${\rm Cov}(N^l_{x,u}(\delta),$ \scalebox{0.96}{$N^r_{x,u}(\delta)) = F\big(\Cstar_{x,u}(s^l_{x,u}(\delta)) \cap \Cstar_{x,u}(s^r_{x,u}(\delta))\big) - F\big(\Cstar_{x,u}(s^l_{x,u}(\delta))\big) F\big(\Cstar_{x,u}(s^r_{x,u}(\delta))\big).$}
\end{corollary}
 
\subsubsection{Asymptotic normality of averaged depth quantile functions.} Here we study the asymptotic behavior of the averaged depth quantile functions used in the applications. Recall that $m_{ij} = \frac{X_i + X_j}{2}$ and $u_{ij} = \frac{X_i - X_j}{\|X_i - X_j\|}.$ For $K > 0$, let
%
$\wh q_{i,K}(\delta) = \frac{1}{n}\sum_{j = 1 \atop j\ne i}^n \wh q_{ij}(\delta){\bf 1}(\|m_{ij} \| \le K),\;i = 1,\ldots, n,$ 
%
and let the corresponding population-based quantities be
%
$q_{i,K}(\delta) = \frac{1}{n}\sum_{j = 1 \atop j\ne i}^n q_{ij}(\delta){\bf 1}(\|m_{ij}\| \le K),\;i = 1,\ldots, n.$
%
Consider the process 
$$\sqrt{n} \big(\wh q_{i,K}(\delta) - {\rm E}\big(q_{i,K}(\delta)\big| X_i\big)\,\big),\,\delta \in [0,1],$$ 
conditional on $X_i$. We will see that this process can be approximated by a U-process, which will then gives us weak convergence to a Gaussian process. Recall the definition of $D_{x,u}(\eta,\epsilon)$ given in Theorem~\ref{as-distr-quant}, and let
%
$D_{ij}(\eta,\epsilon) = D_{m_{ij},u_{ij}}(\eta,\epsilon).$ 
%
We need to following additional assumption:\\[10pt]
{\bf (A4)}  {\em Let $\Psi_{\delta,u}(x) = F\big(A_{x,u}(s_{x,u}(\delta))\big) - F\big(B_{x,u}(s_{x,u}(\delta))\big)$. The directional derivatives $\nabla_u\Psi_{\delta,u}(x)$ at $x$ in direction $u \in S^{d-1}$ exist and are equi-Lipschitz continuous for $u \in S^{d-1}, \delta \in [\epsilon, 1 - \epsilon], \epsilon > 0.$ Moreover, for some $\epsilon_K > 0,$ }
$$\inf_{\|x\| \le K \atop u \in S^{d-1}}\inf_{\delta: \Delta_{x,u}(s_{x,u}(\delta))=0} \Big|\nabla_u\Psi_{\delta,u}(x)\Big| > \epsilon_K.$$

{\sc Discussion of {\bf (A4)}.} Validity of assumption {\bf (A4)} depends on the size and the shape of $\mycomp{D}_{x,u}(\eta,\epsilon)$ {\em as a function in $(x,u)$}, which in turn is another implicit assumption of the shape of the underlying density $f$. Suppose we fix a line $\ell$ (in particular $u$ is fixed), and we consider the directional derivative of $F(A_{x,u}(s))$ and $F(B_{x,u}(s))$ for fixed $s$ w.r.t. $x$ in direction $u$. For $ s > 0$, the directional derivative of $F(B_{x,u}(s))$ w.r.t. $x$ in direction $u$ is given by the negative of the entire surface integral over the entire surface of the frustum $B_{x,u}(s)$. The corresponding directional derivative of $F(A_{x,u}(s))$ equals the difference of the surface integral w.r.t. $f$ over the `lateral' surface and the base of $A_{x,u}(s)$. The fact that this derivative is a difference makes it even more difficult to investigate than the derivative in {\bf (A3)}. Moreover, in {\bf (A4)}, the value $s$ also is a function of $(x,u)$, which adds to the complexity. Nevertheless, heuristic arguments similar to the ones given in the discussion of {\bf (A3)} leads to the intuition that assumption {\bf (A4)} holds in many cases.

\begin{theorem}\label{as-norm-ave-quant}
Suppose that assumptions {\bf (A1)} -  {\bf (A4)} hold. Let $\epsilon > 0.$ For $K > 0$, we have conditional on $X_1$, as $n \to \infty,$ almost surely,
\begin{align*}
\sqrt{n} \big(\wh q_{1,K}(\delta) - {\rm E}\big(q_{12}(\delta){\bf 1}(\|m_{12}\| \le K)\big|X_1\big)\big) \to G_{X_1}(\delta)\;\;\text{weakly, in }\,D[\epsilon, 1 - \epsilon],
\end{align*}
where conditional on $X_1$, $\{G_{X_1}(\delta),\delta \in [\epsilon, 1-\epsilon]\}$ is a tight Gaussian process with mean zero and ${\rm Cov}\big[G_{X_1}(\delta), G_{X_1}(\delta^\prime)\big|X_1\big] =$  ${\rm Cov}\big[ \big(q_{12}(\delta){\bf 1}(\|m_{12}\| \le K), q_{12}(\delta^\prime){\bf 1}(\|m_{12}\| \le K) \big)\big|X_1\big] $.
\end{theorem}
{\sc Remarks.}  (i) Clearly, conditional on $X_i,$ the distribution of $\sqrt{n} \big(\wh q_{i,K}(\delta) - {\rm E}\big(q_{i,K}(\delta)\big| X_i\big)\,\big)$ is the same for every $i = 1,\ldots,n$. 

(ii) Up to truncation, the conditional mean used for centering equals the quantity $E_X\big(\Psi_{x_1,F_{X_2},\delta}(X)\big)$ introduced in section~\ref{random-sets}, where  $X_1 = x_1$ is the observed value of $X_1$. As discussed in section~\ref{random-sets}, the quality of classification based on depth quantile functions as discussed in section~\ref{classification}, depends on the behavior of these quantities as a function in $\delta.$ 

(iii) Theorem~\ref{as-norm-ave-quant} concerns an average of processes conditional on $X_1$, centered at their (conditional) means ${\rm E}\big(q_{12}(\delta)|X_1\big)$. If we, as above, center the individual processes on $q_{12}(\delta)$ instead, meaning we are subtracting the conditional expectation given $X_2$ (and given $X_1$), then the resulting average process $\frac{1}{n}\sum_{j=2}^n \sqrt{n}\big(\wh q_{1j}(\delta) - q_{1j}(\delta)\big){\bf 1}\{\|m_{1j}\| \le K\}$, $ \delta \in [\epsilon, 1 - \epsilon]$ can be approximated by a degenerate $U$-process, which converges at a rate faster than $1/\sqrt{n}$. This is shown in the proof of Theorem~\ref{as-norm-ave-quant}. 

(iv) It is interesting to observe that, in contrast to Theorem~\ref{as-distr-quant}, there is no exceptional set. In fact, the exceptional sets `average out' because we are averaging over many different lines. Here is where assumptions {\bf (A3)} and {\bf (A4)} come into play, because they make sure that individual exceptional sets are `small'. (See proof of Theorem~\ref{as-norm-ave-quant} for details.)




%

\section{Proofs}
Here we present some supplementing results, as well as the proofs of the theoretical results presented in the main part of the manuscript. 

First we present a lemma that is the empirical analog to Lemma~\ref{quant}. To this end, let $\wh I_{x,u}(t) = \{s \in\R: \wh d_{x,u}(s) \le t\},\, t \ge 0$ denote the sublevel sets of the empirical depth function $\wh d_{x,u}(s).$ 

\begin{lemma}\label{prop-emp-depth-quant}
Assume that $F$ and $G$ have positive Lebesgue densities, and fix $(x,u) \in \R^d \times S^{d-1}.$
\begin{itemize}
\item[(a)] The function $s \to \wh d_{x,u}(s), s \in \R$ is a piecewise constant function that is non-increasing and left continuous for $s \le 0$, and non-decreasing and right continuous for $s \ge 0.$ Thus, almost surely, the intervals $\wh I_{x,u}(t), t \ge 0$ are open. We also have $\min_s \wh d_{x,u}(s) = \wh d_{x,u}(0) = 0$ almost surely. 
\item[(b)] Let $0 < \delta < 1.$ With $\wh I_{x,u}(\wh q_{x,u}(\delta)) = (\wh s^{\,l}_{x,u}(\delta), \wh s^{\,r}_{x,u}(\delta)),$ we have almost surely that
\begin{align}\label{depth-quant-discr}
\wh q_{x,u}(\delta) \;\le\; \min \big\{\wh d_{x,u}(\wh s^{\,l}_{x,u}(\delta)),\,\wh d_{x,u}(\wh s^{\,r}_{x,u}(\delta))\big\}\; \le \; \wh q_{x,u}(\delta)  + n^{-1}.
\end{align}
\end{itemize}
\end{lemma}
{\sc Proof.} (a) U-shapedness of $ \wh d_{x,u}(s)$ follows from the fact that the classes $\{A_{x,u}(s), s \le 0\}$, $\{A_{x,u}(s), s \ge 0\},$ $\{B_{x,u}(s), s \le 0\}$ and $\{B_{x,u}(s), s > 0\}$ consist of nested sets, increasing as $|s|$ increases. Thus $s \mapsto F_n(A_{x,u}(s))$ and $s \mapsto F_n(B_{x,u}(s))$ are both non-increasing for $s \le 0$ and non-decreasing for $ s \ge 0$. Using the nestedness again, along with the fact that all the sets $A_{x,u}(s), B_{x,u}(s)$ are closed, it follows that the functions $ F_n(A_{x,u}(s), s \le 0,$ and $F_n(B_{x,u}(s)), s \le 0,$ are both left continuous, and $F_n(A_{x,u}(s), s \ge 0,$ and $F_n(B_{x,u}(s)), s \ge 0,$ are right continuous. Finally, $\wh d_{x,u}(0) = F_n(\{x\}) = 0$ almost surely.

Part (b) follows by observing that almost surely the jump sizes of the functions $F_n(A_{x,u}(s))$ and $F_n(B_{x,u}(s))$ are $n^{-1}$, meaning that the boundaries of the level sets of $\wh d_{x,u}(s)$ are points at which $\wh d_{x,u}(s)$ jumps by $n^{-1}$. \hfill $\square$\\

{\sc Proof of Lemma~\ref{quant}.} (a) Both $A_{x,u}(s)$ and $B_{x,u}(s)$ are monotonically increasing (by inclusion) as $s$ is moving away from $0$ in either direction. Hence, it follows that $d_{x,u}(s)$ is U-shaped, and that the minimum value of $d_{x,u}(s)$ is attained at $s=0$. This, in turn, implies that the sublevel sets are intervals containing $0$. (b) If $F$ is continuous then $d_{x,u}(s)$ is continuous, because the surfaces of the sets and $A_{x,u}(s)$ and $B_{x,u}(s)$ have zero $d$-dimensional Lebesgue measure. Strict U-shapedness of $d_{x,u}(s)$ follows from the fact that $f$ has a positive density, so that both $\big| F(A_{x,u}(s)) - F(A_{x,u}(s^\prime))\big| > 0$ and $\big| F(B_{x,u}(s)) - F(B_{x,u}(s^\prime))\big| > 0$ for $s \ne s^\prime$. (c) This assertion follows from (\ref{alternative-def}) by using the fact that for any (measurable) set $C$, we have $F_{\theta}(C) = F(\sigma O C +\mu),$ and that $\sigma O A_{x,u}(s) + \mu = A_{\sigma O x + \mu, O u}(\sigma s)$ (and similarly for $B_{x,u}(s))$.
\hfill $\square$\\

{\sc Proof of Lemma~\ref{multiscale}}  
Part (a): $F$ and $G$ having a strictly positive Lebesgue densities on $\R^d$ and $\R$, respectively, implies that the sublevel sets $I_{x,u}(t)$ are {\em strictly} monotone (w.r.t. inclusion), and that, for each $\delta > 0,$ there exists a unique $t(\delta)$ with  $G\big([s^l_{x,u}(t(\delta)), s^r_{x,u}(t(\delta))]\big) = \delta.$ That $q_{x,u}(0) = 0$ follows from $d_{x,u}(0) = F(A_{x,u}(0)) = F\{x\} = 0.$

Part (b): Suppose $s \to -\infty$. Then, by construction, $F(A_{x,u}(s)) \to F_{x,u}(x)$ and $F(B_{x,u}(s)) \to 1 - F_{x,u}(x),$ and thus $d_{x,u}(s) = \min\big(F(A_{x,u}(s)), F(B_{x,u}(s))\big)$ $\to \min\big(F_{x,u}(x), 1 - F_{x,u}(x)\big)$. For $s \to \infty$, we have $F(A_{x,u}(s)) \to 1 - F_{x,u}(x)$ and $F(B_{x,u}(s)) \to F_{x,u}(x)$. This implies the assertion.

Part (c):  Here we use Lemma~\ref{quant}. For small enough $\delta$ we have $q_{x,u}(\delta) = F\big(A_{x,u}(s_{x,u}^l(\delta))\big) = F\big(A_{x,u}(s_{x,u}^r(\delta))\big) $. We know that 
$${\rm vol}_d\big(A_{x,u}(s_{x,u}^l(\delta))\big)= \frac{1}{d} h_l V_{d-1} \big(\tan(\alpha)h_l\big)^{d-1} = \big(\tan(\alpha)\big)^{d-1} \frac{1}{d} V_{d-1} h_l^d,$$ 
where $V_{d-1}$ is the volume of the unit ball in ${\mathbb R}^{d-1}$, and $h_l = \|s_x^l(\delta) - x\|$ is the height of the cone $A_{x,u}(s_{x,u}^l(\delta))$. A similar expression holds for $A_{x,u}(s_{x,u}^r(\delta))$ with $h_r$ denoting the height of $A_{x,u}(s_{x,u}^r(\delta)).$ Then, $\delta = G\big[s_{x,u}^l(\delta)), s_{x,u}^r(\delta)] = \big(g(0)+O(h_l + h_r)\big)\big(h_l + h_r\big),$ and, as $\delta \to 0$, 
$$F\big(A_{x,u}(s_{x,u}^l(\delta))\big) = (f(x) + o(h_l)){\rm vol}_d\big(A_{x,u}(s_{x,u}^l(\delta))\big) = c_d f(x) h_l^d \big(1 + O(h_l)\big),$$
where $c_d$ only depends on $d.$ A similar expansion, with $h_l$ replaced by $h_r$, holds for $F\big(A_{x,u}(s_{x,u}^r(\delta))\big).$ 
Since $F\big(A_{x,u}(s_x^r(\delta))\big) = F\big(A_{x,u}(s_{x,u}^l(\delta))\big),$ we obtain $h_r - h_l = o(\delta)$ as $\delta \to 0.$ This implies $2 g(0) h_l = \delta (1 + o(1))$, and thus
$$F\big(A_{x,u}(s_{x,u}^l(\delta))\big) = c_d 2^d f(x) \Big(\tfrac{\delta}{g(0)}\Big)^d (1 + o(1)).$$
This is the assertion.

Part (d): We have \begin{align*}
q^{\theta}_{x,u}(\delta) &= \inf\big\{ t \ge 0: G_\theta(s: d^{\theta}_{x,u}(s) \le t) \ge \delta \big\}\\[3pt]
&=\inf\big\{ t \ge 0: G_\theta(s: d_{\sigma Ox + \mu,Ou}(\sigma s) \le t) \ge \delta \big\} \quad\text{(Lemma 2.1 (c))}\\[3pt]
&=\inf\big\{ t \ge 0: G_\theta\big(\tfrac{s}{\sigma}: d_{\sigma O x + \mu,Ou}(s) \le t\big) \ge \delta \big\}\\[3pt]
&=\inf\big\{ t \ge 0: G\big(s: d_{\sigma x + \mu,u}(s) \le t\big) \ge \delta \big\}\qquad\quad\text{(definition of $G_\theta$)}\\[3pt]
&= q_{\sigma O x + \mu, O u}(\delta).
\end{align*}

\hfill $\square$

\subsection{Proofs of main results}
\label{proofs}
Before presenting the proofs, we give a heuristic high level discussion. One of the main ideas is to approximate $\wh d_{x,u}(s) - d_{x,u}(s) = F_n(\wh C^*_{x,u}(s)) - F(C^*_{x,u}(s))$ by $(F_n - F)(C^*_{x,u}(s))$. In other words, replace  $\wh C^*_{x,u}(s),$ the (random) minimizer of the empirical depth function, by its population version $C^*_{x,u}(s).$ Tools from empirical processes now become available. Replacing $\wh C^*_{x,u}(s)$ by $C^*_{x,u}(s)$ is based on the following observation. Both $C^*_{x,u}(s)$ and $\wh C^*_{x,u}(s)$ lie in $\{A_{x,u}(s), B_{x,u}(s)\}$. Suppose that $s$ is fixed, and $F(A_{x,u}(s)) < F(B_{x,u}(s))$ so that $d_{x,u}(s) = F(A_{x,u}(s)).$ Since $F_n(A_{x,u}(s)) = F(A_{x,u}(s)) + O_P(1/\sqrt{n})$ (and the same holds for $B_{x,u}(s)$), we have for $n$ large enough that with high probability also $F_n(A_{x,u}(s)) < F_n(B_{x,u}(s))$, meaning that $\wh C^*_{x,u}(s) = C^*_{x,u}(s),$ as desired. If $F(A_{x,u}(s)) - F(B_{x,u}(s)) = 0,$ or, more generally, if this difference is of the order $O(1/\sqrt{n})$, then we have problems with this argument. This is the reason for the exceptional sets $S_{x,u}(c)$ and $D_{x,u}(\eta,\epsilon)$ being introduced in the formulation of the theorems. The main challenges in the proofs of weak convergences is to control remainder terms {\em uniformly} over the appropriate ranges of the parameters. For showing equicontinuity, we need to make sure that if $|s - s^\prime| \le \eta$ for $\eta$ small enough, then  $C^*_{x,u}(s)$ and $C^*_{x,u}(s^\prime)$ are both either `$A$-sets' or `$B$-sets, because then $|F(C^*_{x,u}(s)) - F(C^*_{x,u}(s^\prime))| = F\big(C^*_{x,u}(s) \Delta C^*_{x,u}(s^\prime)\big)$, because the $A$-sets are nested, and so are the $B$-sets (as long as $(x,u)$ is fixed). This is where assumption {\bf (A3)} comes into play. In the proof of Theorem~\ref{as-norm-ave-quant}, we also use theory for $U$-processes.

\subsection{Proofs of Theorems~\ref{depth} and \ref{as-distr-quant}} 

{\sc Proof of Theorem~\ref{depth}.}  For $c > 0$, let $E_n(c) = \{\sup_{C \in \D_{x,u}}|\sqrt{n}(F_n - F)(C)| \le c\}$. As $\D_{x,u}$ is a VC-class, empirical processes theory tells us (e.g. see van der Waart and Wellner, 1996), that, for any $\epsilon > 0,$ there exists a $c > 0,$ such that $P(E_n(c)) > 1-\epsilon$. Consequently, by definition of $\D_{x,u}$, for any given $\epsilon > 0,$ there exists $c_\epsilon > 0$, such that, for $c \ge c_\epsilon$,\\[5pt]
\scalebox{0.93}{$P\Big(\max\big(\sup\limits_{s \in \R}|(F_n - F)(A_{x,u}(s)|, \sup\limits_{s \in \R}|(F_n - F)(B_{x,u}(s)|\big) > \tfrac{c}{\sqrt{n}}\Big) \le P\big(E_n(c)\big) \le \epsilon.$}
\vspace*{0.2cm}

On this event, we have for $s \in S_{x,u}(2c/\sqrt{n})$ with $F(A_{x,u}(s)) > F(B_{x,u}(s)) + \frac{2c}{\sqrt{n}}$, say, that  $F_n(A_{x,u}(s)) > F_n(B_{x,u}(s)).$ In other words, for $s \in S_{x,u}(2c/\sqrt{n})$, the unique minimizer $\Cstar_{x,u}(s)$ in the definition of $d_{x,u}(s)$ equals the unique minimizer $\Cstarhat_{x,u}(s)$ in the definition of $\wh d_{x,u}(s)$. It follows that, on $E_n(c)$, we have, for $s \in S_{x,u}\big(2c/\sqrt{n}\big),$
\begin{align*}
\sqrt{n\,}(\wh d_{x,u}(s) - d_{x,u}(s)) = \sqrt{n\,} (F_n - F)(\Cstar_{x,u}(s)) = \min_{C \in \T_{x,u}(s)}\sqrt{n\,} (F_n - F)(C),
\end{align*}
where the last equality holds trivially, since here $\T_{x,u}(s) = \{\Cstar_{x,u}(s)\}.$ Recalling that $\B^{(n)}_{x,u}(s) = $ $\displaystyle{\min_{C \in \T_{x,u}(s)}} \B^{(n)}_F(C)$, we obtain that, on $E_n(c),${\small
\begin{align*}
\sup_{s \in S_{x,u}(\frac{2c}{\sqrt{n}})}\big|\sqrt{n\,}(\wh d_{x,u}(s) - d_{x,u}(s)) - \B^{(n)}_{x,u}(s)\big| \le \max_{C \in \D_{x,u}}\big| \sqrt{n\,} (F_n - F)(C) - \B^{(n)}_F(C) \big|.
\end{align*}
}
If $s$ is such that $\{s\in \ell: \Delta_{x,u}(s) = 0\},$ so that $\T_{x,u}(s) = \{A_{x,u}(s), B_{x,u}(s)\},$ $\wh d_{x,u}(s) = F_n(A_{x,u}(s)) = F_n(B_{x,u}(s))$ and $d_{x,u}(s) = F(A_{x,u}(s)) = F(B_{x,u}(s)),$ then we trivially also have
$$\sqrt{n}\big(\widehat{d}_{x,u}(s) - d_{x,u}(s)\big) = \min_{C \in \T_{x,u}(s)} \sqrt{n}\big(F_n-F\big)(C).$$
Consequently, on $E_n(c)$,
\begin{align}\label{bound}
\sup_{s \in \Sstarsmall_{x,u}(\frac{2c}{\sqrt{n}})}\big| \sqrt{n}\big(\widehat{d}_{x,u}(s) &- d_{x,u}(s)\big) -  \B^{(n)}_{x,u}(s) \big|\nonumber \\[-5pt]
&\le \sup_{C \in \D_{x,u}}\big| \sqrt{n}\big(F_n- F\big)(C) - \B^{(n)}_F(C) \big|. 
\end{align}
Choose $c = \eta_n \sqrt{n}.$ It follows from strong approximation results for empirical processes indexed by VC classes by Brownian $F$-bridges (e.g. see Berthet and Mason, 2006) that there exists a sequence of process $\B^{(n)}_{F}(C)$, such that the term on the right hand-side is $o_P(1).$ The rate of approximation only depends on the VC dimension. We already have seen that the VC dimension of ${\cal D}_{x,u}$ is bounded by $4$, and the first assertion of the theorem follows. 

The second assertion of the theorem follows similarly: Replacing the class ${\cal D}_{x,u}$ in the definition of the event $E_n(c)$ and in (\ref{bound}) by the class $\D$ (defined at the beginning of section~\ref{large sample})  similar arguments as above show that everything also holds uniformly over $(x,u) \in \R \times S^{d-1}.$ It remains to recall that $\D$ also is a VC-class, but with VC dimension depending on $d$. \hfill$\square$\\

{\sc Proof of Theorem~\ref{as-distr-quant}.} Let $d_{x,u}^{\,l}(s)$ and $d_{x,u}^{\,r}(s)$ denote the function $d_{x,u}(s)$ restricted to $s \le 0$ and $s > 0$, respectively. Since $d_{x,u}(s)$ is strictly U-shaped with minimum value in $s=0$ (Lemma~\ref{quant}), both $d_{x,u}^{\,l}(s)$ and $d_{x,u}^{\,r}(s)$ have a unique inverse. It follows that, for every $\delta \in (0,1),$ $\{s \in \R: d_{x,u}(s) \le q_{x,u}(\delta)\} = [s^l_{x,u}(\delta), s^r_{x,u}(\delta)],$ where $s^l_{x,u}(\delta) = \big(d^{\,l}_{x,u}\big)^{-1}(q_{x,u}(\delta))$, and $s^r_{x,u}(\delta) = \big(d^{\,r}_{x,u}\big)^{-1}(q_{x,u}(\delta))$ is a compact interval on $\R$. \\[5pt]
Let $c_n = \eta_n \sqrt{n}$ and set  
%
$\cA_n(c_n) =  \big\{  \sup_{\|x\| \le K\atop u \in S^{d-1}}\sup_{s \in \R}\big|(\wh F_n - F)(A_{x,u}(s)) + (F_n - F)(B_{x,u}(s))\big| \le \tfrac{c_n}{2\sqrt{n}}\big\}.$
%
Using similar arguments as in Proposition~\ref{stoch-bound} (and its proof) and in Proposition~\ref{depthity-quant-bound}, we obtain, that there exists $c^\prime > 0$ with 
%
$P(\mycomp{A}_n(c)) \le c^\prime c^3 e^{-c^2/2}\quad\text{for all }\,c > 0.$
%
In particular, $P\big(\mycomp{\cA}_n(c_n)\big) \to 0$ as $c_n \to \infty.$ Also, it follows from arguments given in the proof of Proposition~\ref{depthity-quant-bound}, that on $\cA_n(c_n)$ we have $ \sup\limits_{\|x\| \le K; u \in S^{d-1}}\sup_{\delta \in [0,1]} \big|\wh q_{x,u}(\delta) - q_{x,u}(\delta) \big| \le c_n / 2\sqrt{n}.$ This will be used below. In the following, we assume that we are on $\cA_n(c_n).$   We will show that there exists a process $\wh{\mathbb Y}_{x,u}$, such that, for $\delta \in D_{x,u}(c_n),$ we have (on $\cA_n(c_n)$),
\begin{align*}
\wh q_{x,u}(\delta) &= \inf\big\{ t: G(s: \wh d_{x,u}(s) \le t) \ge \delta \big\}\\
&\hspace*{-0.7cm}= \inf\big\{ t: G(s: d_{x,u}(s) \le t) \ge \delta - G(s: \wh d_{x,u}(s) \le t)  + G(s: d_{x,u}(s) \le t) \big\}\\
&\hspace*{-0.7cm}\approx  \inf\big\{ t: G(s: d_{x,u}(s) \le t) \ge \delta + \wh{\mathbb Y}_{x,u}(\delta) + r_n) \big\},
\end{align*}
so that\\[-25pt]
\begin{align*}
\wh q_{x,u}(\delta) - q_{x,u}(\delta) \approx q^\prime_{x,u}(\delta)\,\wh {\mathbb Y}_{x,u}(\delta) + q^\prime_{x,u}(\delta)\,r_n.
\end{align*}
It will turn out that $\sqrt{n\,}\,r_n = \omega_n\big(c_n/\sqrt{n} \big),$ where $\omega_n$ denotes the modulus of continuity (with respect to $d_F(C,D)$) of the set-indexed empirical process $\sqrt{n}\big(F_n - F)(C), C \in {\cal D}.$ Empirical process theory for VC-classes implies that $r_n = o_p\big(n^{-1/2}\big),$ so that the asymptotic behavior of  $\sqrt{n}\,\big(\wh q_{x,u}(\delta) - q_{x,u}(\delta)\big)$ is  determined by $\sqrt{n} \wh {\mathbb Y}_{x,u}(\delta)$. In the following we will make this precise. \\[5pt]
To shorten the notation we drop in what follows the indices $x,u$ from\\ $\wh q_{x,u}, q_{x,u}, \wh d_{x,u}, d_{x,u}, s^{\,l}_{x,u}, s^{\,r}_{x,u}, S_{x,u}, A_{x,u}, B_{x,u}$ and $\alpha_{x,u}.$ We begin by noting that, on $\cA_n(c_n)$, and with $\sqrt{n}\,r_n = \omega_n\big(c_n/\sqrt{n} \big),$ we have
\begin{align}
&\inf\big\{ t:  G\big(s: d(s) \le t - (\wh d - d)(s^l(\delta)) + r_n\big)\big) \ge \delta\,\big\} \nonumber\\[2pt]
&\hspace*{1cm} \le \wh q(\delta) \label{upper-lower}\\
& \hspace*{2cm}\le \inf\big\{ t:  G\big(s: d(s) \le t - (\wh d - d)(s^l(\delta)) - r_n\big)\big)\ge \delta\,\big\}.\nonumber
\end{align}
Assume for now that (\ref{upper-lower}) holds (the proof is given below). Empirical process theory for VC-classes tells us, that $\omega_n\big(\frac{c_n}{\sqrt{n}}\big)$ converges to zero as $n \to \infty.$ With $s^l(\delta)$ replaced by $s^r(\delta),$ the same argument holds for $s > 0.$ Combining these two results then implies that on $\cA_n(c_n),$ for $t > 0$ as above, 
 %
%
%
%
%
\begin{align*}
G(s:\,&\wh d(s) \le t) - G(s:\,d(s) \le t) \\
&\le G\big(s \le 0: d(s) \le t - (\wh d - d)(s^l(\delta)) - r_n\big) - G\big(s \le 0: d^{\,l}(s) \le t\big)\\
&\; \;\;+ G\big(s > 0: d(s) \le t- (\wh d - d)(s^r(\delta)) - r_n\big) - G\big(s > 0: d^{\,r}(s) \le t\big).
\end{align*}
%
Since $d(s)$ is strictly U-shaped with minimum value in $x$ (Lemma~\ref{quant}), both $d^{\,l}(s)$ and $d^{\,r}(s)$ have unique inverses. Thus (with slight abuse of notation)
\begin{align*}
G\big(s \le 0: d(s) \le t\big) = G\big(s: \big(d^{\,l}\big)^{-1}(t) \le s \le 0\big) = G(0) - G\big(\big(d^{\,l}\big)^{-1}(t)\big).
\end{align*}
Now, on ${\cal A}_n(c_n)$ and by definition of $S(2\eta_n),$ for all $s \in S(\eta_n), s \le 0$ either $d^{\,l}(s) = F(A(s))$, or $d^{\,l}(s) = F(B(s)).$ Thus it follows from Lemma~\ref{derivative} that $d^{\,l}(s)$ is differentiable at $s \in S(2c_n)$ with derivative $\big(d^{\,l}\big)^{\prime}(s)$ being non-zero for $s \ne 0$. Consequently, for $0 < t < q(1)$, such that $\big(d^{\,l}\big)^{-1}(t) \in S(2c_n)$,
\begin{align*}
\tfrac{d}{dt} G\big(s \le 0: d(s) \le t\big) = - \tfrac{d}{dt} G\big((d^{\,l})^{-1}(t)\big) = - \tfrac{g( (d^{\,l})^{-1}(t) )}{(d^{\,l})^{\,\prime}((d^{\,l})^{-1}(t) ) }.
\end{align*}
Similarly, 
\begin{align*}
\tfrac{d}{dt} G\big(s > 0: d(s) \le t\big) =  \tfrac{g( (d^{\,r})^{-1}(t) )}{(d^{\,r})^{\,\prime}((d^{\,r})^{-1}(t) ) }.
\end{align*}
Recall that, by definition, $(d^{\,l})^{-1}(t) = s^{\,l}(t)$ and $(d^{\,r})^{-1}(t) = s^{\,r}(t)$. It follows from Lemmas~\ref{derivative} and \ref{compact} that also the second derivative of $G\big(s \le 0: d(s) \le t\big)$ exists for values of $t$ in a neighborhood of $q(\delta)$ with $\delta \in D(c_n,\epsilon)$, and that it is uniformly bounded for $\|x\| \le K,$ $u \in S^{d-1}$  (for any choice of $K,c_n, \epsilon > 0$). Using the fact that $\sup_{\|x\|\le K, u \in S^{d-1}}\big|(\wh d - d)(s^r(\delta)) + r_n\big| = O_P(n^{-1/2})$, it follows that, on $\cA_n(c_n)$, and for $ 0 < q(\delta) - \frac{c_n}{2\sqrt{n}} < t < q(\delta) + \frac{c_n}{2\sqrt{n}} < 1$,
\begin{align*}
G\big(s:\,\wh d(s) \le t\big) &- G\big(s:\,d(s) \le t\big)\nonumber \\
&\hspace*{-1.5cm}\le {\textstyle\big(\frac{g(s^{\,r}(\delta))}{d^{\,\prime}(s^{\,r}(\delta))} + O_P\big(\frac{1}{\sqrt{n\,}}\big)\big)\,\big[ (\wh d - d)(s^{\,r}(\delta)) -r_n\big] }\nonumber \\[5pt]
&\hspace*{-0.5cm} + {\textstyle\big(\frac{g(s^{\,l}(\delta))}{d^{\,\prime}(s^{\,l}(\delta))} + O_P\big(\frac{1}{\sqrt{n\,}}\big)\big) \,\big[ (\wh d - d)(s^{\,l}(\delta)) -r_n \big]}\nonumber \\[5pt]
& = :\wh{\mathbb Y}(\delta) + R_n,
\end{align*}
%

where 
\begin{align*}
\wh {\mathbb Y}(\delta) = \frac{g(s^{\,r}(\delta))}{d^{\,\prime}(s^{\,r}(\delta))}\,(\wh d - d)(s^r(\delta)) + \frac{g(s^{\,l}(\delta))}{d^{\,\prime}(s^{\,l}(\delta))} \,(\wh d - d)(s^{\,l}(\delta)),
\end{align*}
and $R_n = o_P(n^{-1/2})$ uniformly in $\|x\| \le K$ and $u \in S^{d-1}$. A similar argument gives us a lower bound with the same first order term and a slightly different remainder term that, however, is of the same order as $R_n$. In the following we thus assume that the inequality is an equality.

 From the above calculations we have that 
%
$q^\prime(\delta) = \big(\,\frac{g(s^r(\delta))}{d^{\,\prime}(s^r(\delta))} - \frac{g(s^l(\delta))}{d^{\,\prime}(s^l(\delta))}\big)^{-1}.$ 
%
Note that $q^\prime(\delta) > 0$ for $\delta > 0$, because $d(s)$ is strictly de(in)creasing for $s < 0$ ($s > 0$). Thus,
\begin{align*}
\sqrt{n}\big(\wh q(\delta) - q(\delta)\big) &= q^\prime(\delta) \sqrt{n}\,\wh{\mathbb Y}(\delta) + o_P(1)\\
&\hspace*{-2.5cm}= \alpha(\delta)\,\sqrt{n}\,(\wh d - d)(s^r(\delta))  - \beta(\delta) \,\sqrt{n}(\wh d - d)(s^l(\delta)) + o_P(1) \\
&\hspace*{-2.5cm} = \alpha(\delta)\,\sqrt{n}\,(F_n - F)\big(\Cstar(s^r(\delta)))  - \beta(\delta) \,\sqrt{n}\,(F_n - F)\big(\Cstar(s^l(\delta))) + o_P(1),
\end{align*}
where $\beta(\delta) = 1 - \alpha(\delta)$ with
\begin{align*}
\alpha(\delta) =  q^{\,\prime}(\delta)\,\tfrac{g(s^r(\delta))}{d^{\,\prime}(s^r(\delta))} = \big(1 - \tfrac{g(s^{\,l}(\delta))}{g(s^{\,r}(\delta))}\cdot \tfrac{d^{\,\prime}(s^{\,r}(\delta))}{d^{\,\prime}(s^{\,l}(\delta))}  \big)^{-1}, 
\end{align*}
satisfying $ 0 \le \alpha(\delta) \le 1.$ This is the first assertion of the theorem.

In order to see the second assertion of the theorem, we use the fact that $\{\sqrt{n} (F_n - F)(C), C \in \D\}$ can be strongly approximated by a Gaussian process with almost surely continuous sample paths. Here we are using Massart (1989), which can be applied here because the VC-classes considered here are subclasses of convex sets, and thus they fulfill the uniform Minkowski condition assumed in Massart's paper (see Definition 2 in Massart, 1989). 

\subsection{Proof of Theorem~\ref{as-norm-ave-quant}} 
First, we write $\sqrt{n} \big(\wh q_{1,K}(\delta) - {\rm E}\big(q_{1,K}(\delta)\big|X_1\big)\big) = \sqrt{n}\big(\wh q_{1,K}(\delta) - q_{1,K}(\delta)\big) + \sqrt{n}\big(q_{1,K}(\delta) - {\rm E}\big(q_{1,K}(\delta)\big|X_1\big)\big).$  We will see that this (approximately) equals a Hoeffding decomposition of a $U$-process of rank $2$, where the first summand on the right is a degenerate $U$-process (converging to zero at rate $n^{-{1/2}}$), while the second summand converges to the asserted limit distribution, which we derive first. By definition, conditional on $X_1$,
\begin{align*}
T_{1,K}(\delta):= &\sqrt{n}\big(q_{1,K}(\delta) - {\rm E}\big(q_{1,K}(\delta)\big|X_1\big)\big)\\
&= \frac{1}{\sqrt{n\,}}\sum_{j = 2}^n \Big[q_{1j}(\delta){\bf 1}(\|m_{1j}\| \le K) - {\rm E}\big(q_{1j}(\delta){\bf 1}(\|m_{1j}\| \le K)|X_1\big)\Big]
\end{align*}
is a sum of iid random variables, and $0 \le \sup_{\delta \in [0,1]}q_{1j}(\delta){\bf 1}(\|m_{1j}\| \le K) \le 1.$ Given $X_1$, the central limit theorem thus immediately gives the convergence of the finite dimensional distributions of the process $T_{1,K}(\delta);\,\delta \in [\epsilon, 1 - \epsilon]$ to the asserted limiting normal distribution. It remains to show asymptotic equicontinuity. To this end, we consider $T_{1,K}(\delta)$ as an empirical process indexed by the class of functions $\{h^K_{\delta}(x) = q_{m_{1,x},u_{1,x}}(\delta){\bf 1}\big(|m_1(x)\| \le K\big),\,\delta \in [\epsilon, 1 - \epsilon]\}$ with $m_{1,x} = \frac{X_1 + x}{2}$, $u_{1,x} = \frac{X_1 - x}{\|X_1 - x\|},$  
and show that the map $[\epsilon, 1 - \epsilon] \ni \delta \to h^K_{\delta}(\cdot)$ is Lipschitz continuous with respect to $L_2(F)$. This then implies that the metric entropy with bracketing w.r.t. $L_2(F)$ of the class of functions $h_{\delta,K}(\cdot), \delta \in [\epsilon, 1 - \epsilon]\}$ behaves like the metric entropy with bracketing of the interval $[\epsilon, 1 - \epsilon]$, which in turn is bounded by $C \epsilon^{-1}$, for some constant $C >0.$ Asymptotic equicontinuity of the process $T_{1,K}(\delta) = \nu_n(h^K_{\delta})$ follows.
 
It remains to show the Lipschitz continuity of $\delta \to h^K_{\delta}(\cdot).$ From Lemma~\ref{continuity}, we have that $\{\delta \to q_{x,u}(\delta), \delta \in [\epsilon, 1 - \epsilon], \|x\| \le K, u \in S^{d-1}\}$ is equi-Lipschitz continuous. Thus, for $\delta,\delta^\prime \in [\epsilon, 1 - \epsilon]$,
\begin{align*}
{\rm E}_F\big( h^K_{\delta}(X) - &h^K_{\delta^\prime}(X)\big)^2 \\
&= \int \big(q_{m_{1,x},u_{1,x}}(\delta) - q_{m_{1,x},u_{1,x}}(\delta^\prime)\big)^2 {\bf 1}\big(\|m_{1,x}\| \le K\big) dF(x)\\
& \le \int \sup_{\|x\| \le K, u \in S^{d-1}}\big(q_{x,u}(\delta) - q_{x,u}(\delta^\prime)\big)^2dF(x)\\
& \le\, C\,\big(\delta - \delta^\prime\big)^2
\end{align*}
for some $C > 0.$ This is what we wanted to show.

Next we show that $\sqrt{n}\big(\wh q_{1,K}(\delta) - q_{1,K}(\delta)\big)$ is asymptotically negligible. To this end, and with $\mycomp{D}_{1j} (\eta, \epsilon) = [\epsilon, 1 - \epsilon] \setminus D_{1j} (\eta, \epsilon)$, we split up this process in two terms,
\begin{align*}
\sqrt{n}\big(\wh q_{1,K}(\delta) &- q_{1,K}(\delta)\big) \\
& = \frac{1}{n}\sum_{j=2}^n \sqrt{n}\big(\wh q_{1j}(\delta) - q_{1j}(\delta)\big) {\bf 1}(\| m_{ij}\| \le K;\,\delta \in \mycomp{D}_{1j} (\eta, \epsilon))\\
& + \frac{1}{n}\sum_{j=2}^n \sqrt{n}\big(\wh q_{1j}(\delta) - q_{1j}(\delta)\big){\bf 1}(\| m_{ij}\| \le K;\,\delta \in D_{1j}(\eta, \epsilon))\nonumber\\
&= \qquad(I)_\eta\quad+ \quad (II)_\eta,
\end{align*}
where $\eta > 0$ is arbitrary. 
We show that, for each $\epsilon, \gamma > 0,$ there exists $\eta > 0$ such that both 
\begin{align}\label{prop}
P\big(|(I)_\eta| > \gamma\big) \to 0,\;\;\text{and}\;\;P\big(|(II)_\eta| > \gamma\big) \to 0\quad\text{as }\; n\to \infty.
\end{align} 
Since $\sqrt{n}\big(\wh q_{1,K}(\delta) - q_{1,K}(\delta)\big)$ does not depend on $\eta$, this implies the result.

First consider $(I)_\eta$. Since $\max\limits_{j}\sup\limits_{\delta \in [\epsilon, 1 - \epsilon]}\big|\sqrt{n} \big(\wh q_{1j}(\delta) - q_{1j}(\delta)\big)\big| = O_P(1),$ which follows from Proposition~\ref{depthity-quant-bound}, we obtain that
\begin{align*}
\Big| \frac{1}{n}\sum_{j=2}^n \sqrt{n}\big(&\wh q_{1j}(\delta) - q_{1j}(\delta)\big){\bf 1}(\delta \in \mycomp{D}_{1j}(\eta, \epsilon)) {\bf 1}(\| m_{1j}\| \le K)\Big| \\[-5pt]
& \le O_P(1)\; \sup_{\delta \in [\epsilon, 1 - \epsilon]}\frac{1}{n}\sum_{j=2}^n {\bf 1}(\delta \in  \mycomp{D}_{1j}(\eta,\epsilon))){\bf 1}(\| m_{1j}\| \le K).
\end{align*}
The last sum is a sum of iid random variables with (conditional) mean bounded by $P(\delta \in  \mycomp{D}_{1j}(\eta,\epsilon)| X_1)$, where $\sup_{\delta \in [\epsilon, 1 - \epsilon]}P(\delta \in \mycomp{D}_{1j}(\eta,\epsilon)| X_1) = o(1)$ a.s., as $\eta \to 0.$ To see the latter, observe that by assumption {\bf (A4)}, for any $u \in S^{d-1}$ and $\delta \in [\epsilon, 1 - \epsilon]$, ${\cal H}^1\big\{x \in \ell_{0,u}: \big|\Psi_{\delta,u}(x)\big| \le \eta\big\} \lesssim \eta,$ where ${\cal H}^1$ denotes $1$-dimensional Hausdorff measure, and `$\lesssim$' means `$\le$' up to an universal (multiplicative) constant. Thus, with $m_1(x) = \frac{X_1 + x}{2}, u_1(x) = \frac{X_1 - x}{\|X_1 - x\|}$ and $\delta \in [\epsilon, 1 - \epsilon]$, for each $x_1 \in \R^d$, 
\begin{align}
&P\big(\delta \in \mycomp{D}_{1j}(\eta,\epsilon),\,\big\| {\textstyle \frac{X_j + X_1}{2}}\big\| \le K\big| X_1= x_1\big) \nonumber\\
&= P\big( \delta \in \mycomp{D}_{m_1(X_j),u_1(X_j)}(\eta,\epsilon),\,\big\| {\textstyle \frac{X_j + X_1}{2}}\big\| \le K \Big| X_1= x_1\big)\nonumber\\
&\le F\big(x: \Big| \Psi_{\delta, \frac{x - x_1}{\|x - x_1\|}}\big({\textstyle \frac{x + x_1}{2}}\big)\Big| \le \eta,\,\big\| {\textstyle \frac{x + x_1}{2}}\big\| \le K \big)\nonumber \\
&\le F\big( x: \Big| \Psi_{\delta, u_{x,x_1}}\big(x_1 + \tfrac{u_{x,x_1} r_{x,x_1}}{2}\big) \Big| \le \eta, \| x\| \le \|x_1\| + 2K\big)\nonumber\\[-7pt]
\intertext{\hspace*{4cm} (where $u_{x,x_1} = \frac{x - x_1}{\|x -x_1\|}$ and $r_{x,x_1} = \|x - x_1\|$) \vspace*{-0.2cm} } 
& \le F\Big(x = r u: \Big| \Psi_{\delta, u}\big({\textstyle x_1 +\frac{r u}{2}}\big)\Big| \le \eta,\,\|r\| \le 2(\|x_1| + K), u \in S^{d-1}\Big) \nonumber\\[3pt]
%
%
& =  \int_{S^{d-1}}\int_0^{2(\|x_1\| + K)} r^{d-1} {\bf 1}\big( \big| \Psi_{\delta, u}\big(\tfrac{x_1 + ru}{2}\big)\big| \le \eta \big)f(ru)\,dr\,d{\cal H}^{d-1}(u) \label{ooo}\\[3pt]
&\hspace*{1cm}\lesssim\; \eta,\nonumber
\end{align}
where we have been using polar coordinates in $\R^d$. Thus, in order to conclude that $(I)_\eta$ satisfies (\ref{prop}), it remains to show that, for each $\eta > 0$,{\small
\begin{align}\label{unif-except}
\sup_{\delta \in [\epsilon,1-\epsilon]}\Big| \frac{1}{n}\sum_{j=1}^n{\bf 1}\big(\delta \in  \mycomp{D}_{1j}\big(\eta,\epsilon),\,\|m_{1j}\| \le K\big) - P\big(\delta \in \mycomp{D}_{1j}(\eta,\epsilon),&\,\|m_{1j}\| \le K| X_1\big)\Big| \nonumber\\
&= o_P(1).
\end{align}}
It is shown below that this uniform law of large numbers holds. 

In order to show that also $(II)_\eta$ satisfies (\ref{prop}), we use the approximation derived in Theorem~\ref{as-distr-quant}. Recall the notation $\alpha_{x,u}(\delta), s_{x,u}^{\,l}(\delta)$ and $s_{x,u}^{\,r}(\delta)$ (cf. Theorem~\ref{as-distr-quant} and the discussion immediately before the formulation of this theorem), and let $\alpha_{1j}(\delta), s_{1i}^{\,l}(\delta)$ and $s_{1j}^{\,r}(\delta),$ respectively, denote these expressions with $(x,u)$ replaced by $(m_{1j},u_{1j})$. With 
$$\nu^{\,\wh{ij}}_n(s^l_{x,u}(\delta))  = \frac{1}{\sqrt{n}}\sum_{k=1\atop k \ne i,j}^n \Big[{\bf 1}(X_k \in \Cstar_{x,u}(s^l_{x,u}(\delta))) - F\big(\Cstar_{x,u}(s^l_{x,u}(\delta))\big)\Big]$$ 
let $M_{n,\wh j}(\delta) = \alpha_{1j}(\delta) \nu^{\,\wh{1j}}_n(s^l_{1j}(\delta))- (1-\alpha_{1j}(\delta))\,\nu^{\,\wh {1j}}_n(s^r_{1j}(\delta)),$ and
$$ \wh K_{1j}(\delta)  = \big[\sqrt{n\,}\big(\wh q_{1j}(\delta) - q_{1j}(\delta)\big) -M_{n,\wh j}(\delta)\big] {\bf 1}(\delta \in D_{1j}(\eta, \epsilon)).$$
Now write $(II)_\eta = A_{1,\eta} + A_{2,\eta}$, where
\begin{align*}
A_{1,\eta} = &\frac{1}{n}\sum_{j=2}^n \wh K_{1j}(\delta) {\bf 1}(\| m_{1j}\| \le K)\qquad\text{and}\\
A_{2,\eta}=  &\frac{1}{n}\sum_{j=2}^n M_{n,\wh j}(\delta) {\bf 1}(\delta \in D_{1j}(\eta, \epsilon)){\bf 1}(\| m_{1j}\| \le K).
\end{align*}
We have  $A_{1,\eta} = o_P(1)$, because 
by definition of $D_{1j}(\eta,\epsilon)$, for any $\epsilon, \eta > 0$,
\begin{align*} 
&\sup_{\delta \in [\epsilon, 1 - \epsilon]}\sup_{j = 1,\ldots,n}\big| \wh K_{1j}(\delta){\bf 1}\big(\|m_{1j} \| \le K\big)\, \big| \\
%
%
&\hspace*{0.5cm}\le \sup_{u \in S^{d-1} \atop \|x\| \le K} \sup_{\delta \in D_{x,u}(\eta,\epsilon)}\Big| \sqrt{n}\big(\wh q_{x,u}(\delta) - q_{x,u}(\delta) \big) -M_n(\delta)\Big| + o_P(1)\\
&\hspace*{0.5cm}= o_P\big(1\,\big),
\end{align*}
where $M_n(\delta)$ is as $M_{n,\wh j}(\delta)$, but with $\nu_n^{\wh{1j}}$ replaced by $\nu_n$. In the last step we use Theorem~\ref{as-distr-quant}. 
Now we show that $A_{2,\eta} = o_P(1).$ Recall that all the arguments here are conditional on $X_1$, which, for simplicity, is not indicated explicitly in the following. Plugging in the definitions of $\nu^{\,\wh{1j}}_n(s^{\,l}_{1j}(\delta))$ and $\nu^{\,\wh{1j}}_n(s^{\,r}_{1j}(\delta)),$  and using the notation
{\small
\begin{align*}
w_\delta(X_j,X_k) =\big[\alpha_{1j}(\delta) {\mathbf 1}\big(X_k \in C^*_{1j}(s_{1j}^{\,l}(\delta))\big) &- \big(1 - \alpha_{1j}(\delta)\big){\mathbf 1}\big(X_k \in C^*_{1j}(s_{1j}^{\,r}(\delta))\big)\big]\, \\ 
& \times\;{\bf 1}\big(\|m_{1j} \| \le K\big) {\bf 1}\big(\delta \in D_{1j}(\eta,\epsilon)\big),
\end{align*}}
we have
\begin{align*}
\frac{1}{n}\sum_{j=2}^n &M_{n,\wh j}(\delta)\,{\bf 1}\big(\|m_{1j} \| \le K\big) {\bf 1}\big(\delta \in D_{1j}(\eta,\epsilon)\big)\\
&= \frac{\sqrt{n}}{n^2}\sum_{j=2}^n\sum_{k=2\atop k \ne j}^n \big[w_{\delta}(X_j,X_k) - {\rm E}(w_{\delta}(X_j,X_k)|X_j)\big]\\
& = c_n\,\frac{\sqrt{n}}{\binom{n-1}{2}}\sum_{2 \le j< k \le n} \widetilde w_{\delta}(X_j,X_k) 
\; := \; c_n\,\sqrt{n}\, U_{n}(\delta),
\end{align*}
where $c_n \to \frac{1}{2}$ as $n\to\infty$. Cconditional on $X_1$, the term $U_{n}(\delta)$ is a mean zero $U$-process of rank $2$ (based on $(n-1)$ observations) with $\widetilde w_\delta(x_j,x_k) = \frac{1}{2} [w_{\delta}(x_j,x_k) + w_{\delta}(x_k,x_j)- ({\rm E}(w_{\delta}(x_j,X_k) + w_{\delta}(X_j,x_k)) ]$ as symmetric kernel. Since we are subtracting ${\rm E}(w_{\delta}(X_j,X_k)|X_j)$ from each summand, $U_n(\delta)$ is a degenerate $U$-process. It follows from well-known theory for $U$-statistics that, for each fixed $\delta \in (0,1)$, $\sqrt{n}\,U_n(\delta) \to 0$ in probability. We show below that this also holds uniformly in $\delta \in [\epsilon, 1 - \epsilon]$. 

First we show that there exists a constant $C > 0$ such that, for $\delta, \delta^\prime \in [\epsilon, 1 - \epsilon],$ we have $|\delta - \delta^\prime| \le \eta\;\Rightarrow\;{\rm E}\big(w_\delta(X_j,X_k) - w_{\delta^\prime}(X_j,X_k) \big)^2 \le C \eta$, for $\eta > 0$ small enough. In the following, to shorten notation, we let `${\rm E}$` and '$P$' denoting conditional expectation and probability, respectively, given $X_1$.  For $2 \le j < k \le n$, 
\begin{align}
&{\rm E}\big(w_\delta(X_j,X_k) - w_{\delta^\prime}(X_j,X_k) \big)^2\nonumber \\
&\hspace*{0.5cm}\le 4{\rm E}\big[\big[{\mathbf 1}\big(X_k \in C^*_{1j}(s_{1j}^{\,l}(\delta))\big) - {\mathbf 1}\big(X_k \in C^*_{1j}(s_{1j}^{\,l}(\delta^\prime))\big)\big] {\bf 1}\big(\|m_{1j} \| \le K\big) \big]^2\label{term1} \\
&\hspace*{1cm}+ 4{\rm E}\big[\big[{\mathbf 1}\big(X_k \in C^*_{1j}(s_{1j}^{\,r}(\delta))\big) - {\mathbf 1}\big(X_k \in C^*_{1j}(s_{1j}^{\,r}(\delta^\prime))\big)\big] {\bf 1}\big(\|m_{1j} \| \le K\big) \big]^2\label{term2}\\
&\hspace*{2cm} +\; 8 {\rm E}\big[\big({\bf 1}\big(\delta \in D_{1j}(\eta,\epsilon)\big) - {\bf 1}\big(\delta^\prime \in D_{1j}(\eta,\epsilon)\big)\big) {\bf 1}\big(\|m_{1j} \| \le K\big)\big]^2.\label{term3}\\
&\hspace*{5cm}+\;8{\rm E}\big[ \alpha_{1j}(\delta) - \alpha_{1j}(\delta^\prime)\big]^2.\label{term4}
\end{align} 
Equi-Lipschitz continuity of the maps $[\epsilon, 1 - \epsilon] \in \delta \to \alpha_{x,u}(\delta)$ (\ref{term4}), follows directly from the definition of these functions along with (\ref{q-prop}). As for (\ref{term3}), observe that, by using equi-Lipschitz continuity of  $\delta \to \Psi_{\delta,u}(x) = F\big(A_{x, u}(s_{x, u}(\delta))\big) - F\big(B_{x, u}(s_{x, u}(\delta))\big), \|x\| \le K, u \in S^{d-1}$, for some $C>0$, and $\delta, \delta^\prime \in [\epsilon, 1 - \epsilon]$,
%
\begin{align}
{\rm E}&\big[\big({\bf 1}\big(\delta \in D_{1j}(\eta,\epsilon)\big) - {\bf 1}\big(\delta^\prime \in D_{1j}(\eta,\epsilon)\big)\big)^2 {\bf 1}\big(\|m_{1j} \| \le K\big)\Big| X_1 = x_1\big]\nonumber\\
%
%
& = F\Big(x: \Big|\Psi_{\delta,u_{x,x_1}}\big(x_1 + \tfrac{r_{x,x_1} u_{x,x_1}}{2}\big)\Big| > \eta; \nonumber\\
&\hspace*{2cm}\,\Big|\Psi_{\delta^\prime,u_{x,x_1}}\big( x_1 + \tfrac{r_{x,x_1} u_{x,x_1}}{2}\big)\Big| \le \eta;\,\|x\| \le \|x_1\| + 2K \Big)\nonumber\\
&\hspace*{0.8cm} + F\Big(x: \Big|\Psi_{\delta^\prime,u_{x,x_1}}\big(x_1 + \tfrac{r_{x,x_1} u_{x,x_1}}{2}\big)\Big| > \eta; \nonumber\\
&\hspace*{2cm}\,\Big|\Psi_{\delta,u_{x,x_1}}\big( x_1 + \tfrac{r_{x,x_1} u_{x,x_1}}{2}\big)\Big| \le \eta;\,\|x\| \le \|x_1\| + 2K \Big)\nonumber\\
%
%
&\le F\big(x: \big|\Psi_{\delta,u_{x,x_1}}\big(x_1 + \tfrac{r_{x,x_1} u_{x,x_1}}{2}\big)\big| \in \big[\eta, \eta + C|\delta - \delta^\prime|\big];\,\|x\| \le \|x_1\| + 2K\big)\nonumber\\[2pt]
&\hspace*{0.2cm} + F\big(x: \big|\Psi_{\delta^\prime,u_{x,x_1}}\big(x_1 + \tfrac{r_{x,x_1} u_{x,x_1}}{2}\big)\big| \in \big[\eta, \eta + C|\delta - \delta^\prime|\big];\,\|x\| \le \|x_1\| + 2K\big)\nonumber\\
%
%
&\le \,2\int\limits_{S^{d-1}}\int\limits_0^{2(\|x_1\| + K)} r^{d-1}\,{\bf 1}\big( |\Psi_{\delta^\prime,u}(ru)| \in[\eta, \eta + C|\delta - \delta^\prime|\big]\big) f(ru) dr d{\cal H}^{d-1}(u)\nonumber\\
&\hspace*{2cm}\lesssim \big|\delta - \delta^\prime\big|,\label{Lip}
\end{align}
where the second to last inequality uses similar arguments as in (\ref{ooo}). Next  consider (\ref{term1}). (The term in (\ref{term2}) can be treated similarly.) First observe that (\ref{Lip}) says that the probability of two values $\delta,\delta^\prime$ not both lying in $D_{1j}(\eta,\epsilon)$ is bounded by a (universal) constant times the distance between $\delta$ and $\delta^\prime$. By definition of $D_{1j}(\eta,\epsilon),$ $\delta,\delta^\prime$ both lying in this set means that $s^{\,l}_{1j}(\delta), s^{\,l}_{1j}(\delta^\prime)$ and $s^{\,r}_{1j}(\delta), s^{\,r}_{1j}(\delta^\prime)$ all lie in $S_{1j}(\eta,\epsilon)$. 

As is shown in Lemma~\ref{short-intervals}, for each $x\in \R^n, u \in S^{d-1}$, the set $\mycomp{S}_{x,u}(\eta,\epsilon)$ consists of a union of finitely many intervals $I_j, j = 1,\ldots,J$, with lengths $|I_j|$ satisfying $c\eta \le |I_j| \le C\eta,$ and constants $c,C,J$ not depending on $x,u$. Thus, $S_{x,u}(\eta,\epsilon)$ consists of a union of intervals $J_i$, separated by `small' intervals of length at least $c\eta.$ It is shown in Lemma~\ref{continuity} that the (monotonic) functions $\delta \to s_{x,u}^{\,l}(\delta)$ with $\delta$ ranging over $[\epsilon, 1 - \epsilon]$ are equi-Lipschitz continuous, i.e. there exists a constant $C_1 > 0$ (not depending on $(x,u)$), such that $\big|s_{x,u}^{\,l}(\delta) - s_{x,u}^{\,l}(\delta^\prime)\big| \le C_1 |\delta - \delta^\prime|$ for all $\|x \| \le K$ and $u \in S^{d-1}$. It follows that, if $|\delta - \delta^\prime| < \frac{c}{C_1}\eta$, then $\big|s_{x,u}^{\,l}(\delta) - s_{x,u}^{\,l}(\delta^\prime)\big| \le C_1 |\delta - \delta^\prime\big| \le c \eta.$  The same holds for the functions $s_{x,u}^{\,r}(\delta)$.

As a consequence of all this, (\ref{Lip}) says that, provided $|\delta - \delta^\prime| < \frac{c}{C_1}\eta$, the probability of $s_{1j}^{\,l}(\delta), s_{1j}^{\,l}(\delta^\prime)$ not both lying in the same interval $J_i$ is bounded by a constant times their distance (the same holds for $s_{1j}^{\,r}(\delta), s_{1j}^{\,r}(\delta^\prime)$). Within each interval $J_i$, all the minimizing sets $C^*_{1j}(s), s \in J_i$ are either all of the form $A_{1j}(s)$ or they are all of the form $B_{1j}(s)$, and (recall), for each $(x,u)$, the class $\{A_{x,u}(s), s \in \R\}$ is a nested class of sets, and so is the class $\{B_{x,u}(s), s \in \R\}$. Using this, we obtain for $|\delta - \delta^\prime| \le \frac{c}{C_1}\eta$, and some $C^* > 0$, that
%
\begin{align*}
{\rm E}&{\rm E}\big[\big[{\mathbf 1}\big(X_k \in C^*_{1j}(s_{1j}^{\,l}(\delta))\big) - {\mathbf 1}\big(X_k \in C^*_{1j}(s_{1j}^{\,l}(\delta^\prime))\big)\big] {\bf 1}\big(\|m_{1j} \| \le K\big) \big]^2 \big| X_j\big)\\
&\le {\rm E}\big[F\big(C^*_{1j}(s^{\,l}_{1j}(\delta))\,\Delta\,C^*_{1j}(s^{\,l}_{1j}(\delta^\prime))\big){\bf 1}\big(\|m_{1j} \| \le K,\,\delta, \delta^\prime \in  D_{1j}(\eta, \epsilon)\big)\big] \\
&\hspace*{7cm} + \big(1 - P\big(\delta, \delta^\prime \in  D_{1j}(\eta, \epsilon\big)\big)\\
& \le {\rm E}\big[\max \big\{F\big( A_{1j}(s^{\,l}_{1j}(\delta))\,\Delta\,A_{1j}(s^{\,l}_{1j}(\delta^\prime) )\big), F\big( B_{1j}(s^{\,l}_{1j}(\delta))\,\Delta\,B_{1j}(s^{\,l}_{1j}(\delta^\prime) )\big)\big\}\\
&\hspace*{1.5cm}\times {\bf 1}\big(\|m_{1j} \| \le K,\,\delta, \delta^\prime \in  D_{1j}(\eta, \epsilon\big)\big]  
%
+ \big(1 - P\big(\delta, \delta^\prime \in  D_{1j}(\eta, \epsilon\big)\big)\\
&\le\max\big\{  \sup_{\|x\| \le K; u \in S^{d-1}} \big|\,F\big(A_{x,u}(s^{\,l}_{x,u}(\delta))\big) - F\big(A_{x,u}(s^{\,l}_{x,u}(\delta^\prime))\big)\,\big|,\; \\
&\hspace*{2.5cm}\sup_{\|x\| \le K; u \in S^{d-1}} \big|\,F\big(B_{x,u}(s^{\,l}_{x,u}(\delta))\big) - F\big(B_{x,u}(s^{\,l}_{x,u}(\delta^\prime))\big)\,\big|\big\} \\
&\hspace*{7cm} +  \big(1 - P\big(\delta, \delta^\prime \in  D_{1j}(\eta, \epsilon)\big)\big)\\
& \le C^* |\delta - \delta^\prime|.
\end{align*}
To complete the proof of the theorem, it remains to show (\ref{unif-except}), and that $\sup_{\delta \in [\epsilon, 1 - \epsilon]}\sqrt{n}\,U_n(\delta) = o_P(1).$  The proofs of both of these properties involve arguments from empirical process theory. For this we need to estimate $L_2(F \otimes F)$-bracketing numbers $N_{2,B}(\eta, \tilde{\cal W} )$ of the class of functions $\tilde{\cal W} = \{\widetilde w_\delta(\cdot,\cdot),\,\delta \in [\epsilon, 1 - \epsilon]\},$ and of the $L_2(F)$-bracketing numbers of the class of functions $\big\{\delta \to {\bf 1}(\delta \in D_{1j}(\eta,\epsilon\big),\,\delta \in [\epsilon, 1 - \epsilon]\big\}$, respectively. 

Lipschitz-continuity of $\delta \to \widetilde w_\delta(\cdot,\cdot)$, as just shown, implies that $N_{2,B}(\eta, \tilde{\cal W} )$ $\le N_{2,B}(\eta^2/C^*), {\cal I}_\epsilon ),$ where the latter denotes the bracketing number of the class of intervals ${\cal I}_\epsilon = \{[\epsilon, t],\,\epsilon \le t \le 1 - \epsilon\}$. It is well known that $N_B(\eta, {\cal I}_\epsilon ) = O(1/\eta).$ Thus,  the covering integral $\int_0^1 \sqrt{\log N_B(\eta, \tilde{\cal W})\,}d\eta$ is finite. The fact that $\sup_{\delta \in [\epsilon, 1 - \epsilon]}|\sqrt{n\,}U_n(\delta)| = o_P(1)$ now follows from results in Arcones and Gin\'{e} (1993). 

It remains to show (\ref{unif-except}). This follows by using $L_2(F)$-Lipschitz continuity of the class $\big\{x \to {\bf 1}\big(\delta \in D_{m_1(x),u_1(x)}(\eta,\epsilon);\, \|x\| \le K \big), \delta \in [\epsilon, 1 - \epsilon]  \big\}$ that is shown in (\ref{Lip}). This implies finiteness of the $L_2(F)$-bracketing numbers of this class. The fact that finiteness of the bracketing numbers implies a Glivenko-Cantelli result, is a basic result in empirical process theory. 
\hfill $\square$\\

\subsection{Proofs of Propositions~\ref{stoch-bound}, \ref{uniform-bound} and \ref{depthity-quant-bound}}

{\sc Proof of Proposition~\ref{stoch-bound}.} We have
\begin{align*}
&|\wh d_{x,u_0}(s) - d_{x,u_0}(s)| \\
&= \big|\min\big( F_n(A_{x,u_0}(s)), F_n(B_{x,u_0}(s))\big) - \min\big( F(A_{x,u_0}(s)), F(B_{x,u_0}(s))\big)\big| \\
&\le \max\big(\big|(F_n - F)(A_{x,u_0}(s))\big|, \big|(F_n - F)(B_{x,u_0}(s))\big|\big)\\
& \le 2 \sup_{C \in {\cal D}_{x_0,u_0}}\big|(F_n - F)(C)\big|.
\end{align*}
As discussed above, ${\cal D}_{x_0,u_0}$ is a VC-class with VC-dimension bounded by 4. The asserted bound thus follows from known result from empirical process theory (see Theorem 2.14.13 in van der Vaart and Wellner, 1996).\hfill $\square$\\

{\sc Proof of Proposition~\ref{uniform-bound}.} The idea of the proof is to condition on $X_i,X_j$, and to then apply Proposition~\ref{stoch-bound}. Since $\wh d_{ij}(s)$ depends on $(X_i,X_j)$ not only though $m_{ij}$ and $u_{ij}$, but also through $F_n$, we consider $F_n^{-ij},$ the empirical distribution with $X_i,X_j$ being dropped. With 
$$\wh d^-_{ij}(s) = \min\big(F_n^{-ij}(A_{ij}(s)), F_n^{-ij}(B_{ij}(s))\big)$$
we clearly have $\sup_{s \in \R}\big|\wh d_{ij}(s) - \wh d_{ij}^-(s)\big| \le 2/n$, almost surely. Thus, applying Proposition~\ref{stoch-bound} conditional on $X_i,X_j$, and using a simple union bound, we immediately obtain that, for $C> 0$ sufficiently large,
\begin{align*}
P\big( \max_{i,j = 1,\ldots,n \atop i \ne j} \sup_{s \in \R}\big|\wh d_{i,j}(s) - d_{i,j}(s)\big|& \ge C \sqrt{\tfrac{\log n}{n}}\; \big) \\
& \le {\textstyle\binom{n}{2}}c \log^{3/2} n\,C^3 \exp\big\{-\tfrac{C \log n}{2}\big\},
\end{align*}
and the right-hand side converges to zero for $C > 4.$ This gives the first part of the asserted rate. Moreover, it is known from empirical process theory that for a VC-class of sets ${\cal C}$ with VC-dimension $O(d)$, we have $\sup_{C \in {\cal C}}\big|(F_n - F)(C)| = O_P\big(\sqrt{\frac{d}{n}}\big)$ (e.g. see Bousquet et al. 2004). Applying this to ${\cal D}$, along with the fact that 
\begin{align*}
\max_{i,j = 1,\ldots,n \atop i \ne j} \sup_{s \in \R}\big|\wh d_{i,j}(s) - d_{i,j}(s)\big| &\le \sup_{x\in \R^d,u \in S^{d-1}}\sup_{s \in \R}\big|\wh d_{x,u}(s) - d_{x,u}(s)\big| \\
&\le 2\,\sup_{C \in \D}\big| (F_n - F)(D)\big|
\end{align*}
(see proof of Proposition~\ref{stoch-bound}), we obtain the second part of the asserted rate. This completes the proof of the asserted rate. The asserted adaptation to `sparsity', follows from the simple geometric fact that intersections of cones with hyperplanes that are passing through the axis of symmetry of the cone are again cones, lying in the hyperplane. \hfill$\square$\\

{\sc Proof of Proposition~\ref{depthity-quant-bound}} With $\wh b_{x_0,u_0} = \sup\limits_{x \in \ell_{x_0,u_0} \atop s \in \R}\big| \wh d_{x,u_0}(s)- d_{x,u_0}(s) \big|,$ and $x \in \ell_{x_0,u_0}$, we have\\[-10pt]
\begin{align*}
\wh q_{x,u_0}(\delta) = \inf\big\{ &t: G(s: \wh d_{x,u_0}(s) \le t) \ge \delta \big\}\\
& \le \inf\big\{ t: G(s: d_{x,u_0}(s) \le t - \wh b_{x_0,u_0}) \ge \delta ) \big\} \\
&= q_{x,u_0}(\delta) + \wh b_{x_0,u_0},
\end{align*}
and similarly,
\begin{align*}
\wh q_{x,u_0}(\delta) \ge \inf\big\{ t: G(s: d_{x,u_0}(s) \le t + \wh b_n) \ge \delta) \big\} = q_{x,u_0}(\delta) - \wh b_{x_0,u_0}, 
\end{align*}
so that $\sup\limits_{x \in \ell_{x_0,u_0} \atop \delta \in [0,1]} \big|\wh q_{x,u_0}(\delta) - q_{x,u_0}(\delta)\big| \le \wh b_{x_0,u_0}$. Using Proposition~\ref{stoch-bound}, the assertion follows.
\vspace*{-0.5cm}
\begin{flushright} $\square$\end{flushright} 

\subsection{Proof of miscellaneous technical results}

{\sc Proof of (\ref{upper-lower}).}  To see this, observe that,  on $A_n(c_n)$, 
\begin{align*}
\wh q_{x,u}(\delta) = \inf\big\{ &t: G(s: \wh d_{x,u}(s) \le t) \ge \delta \big\}\\
& = \inf\big\{ t:  G(s \le 0: \wh d_{x,u}(s) \le t) + G(s > 0: \wh d_{x,u}(s) \le t) \ge \delta;\\
&\hspace*{4cm} \,q_{x,u}(\delta) - \tfrac{c_n}{2\sqrt{n}} \le t \le q_{x,u}(\delta) + \tfrac{c_n}{2\sqrt{n}}\big\}.
\end{align*}
For $q_{x,u}(\delta) - \frac{c_n}{2\sqrt{n}} \le t,$ we have, by again using the definition of $A_n(c_n),$ that $q_{x,u}(\delta) - \frac{c_n}{\sqrt{n}} \le t - (\wh d_{x,u} - d_{x,u})(s)$. For such values of $t$, it follows that 
\begin{align}\label{abc}
&G(s \le 0: \wh d_{x,u}(s) \le t)  \nonumber\\
&= G\big(s \le 0: d_{x,u}(s) \le t - (\wh d_{x,u} - d_{x,u})(s)\big) \nonumber\\
&\le G\big(s \le 0: d_{x,u}(s) \le q_{x,u}(\delta) - \tfrac{c_n}{\sqrt{n}} \big) \nonumber\\
&\hspace*{1cm}+ G\big(s \le 0: q_{x,u}(\delta) - \tfrac{c_n}{\sqrt{n}} < d_{x,u}(s) \le t - (\wh d_{x,u} - d_{x,u})(s)\big).
\end{align}
Since $t \le q_{x,u}(\delta) + \tfrac{c_n}{2\sqrt{n}}$, and since on $A_n(c_n)$ we have $\sup_{x \in \ell_{x,u},s \in \R}\big|\wh d_{x,u}(s) - d_{x,u}(s)| \le \frac{c_n}{2\sqrt{n}},$ we see that $q_{x,u}(\delta) - \tfrac{c_n}{\sqrt{n}} < d_{x,u}(s) \le t - (\wh d_{x,u} - d_{x,u})(s)$ implies $q_{x,u}(\delta) - \tfrac{c_n}{\sqrt{n}} < d_{x,u}(s) \le q_{x,u}(\delta) + \tfrac{c_n}{\sqrt{n}}$. 
Let ${\rm B}^l_{x,u}(\delta) = \{s \le 0: q_{x,u}(\delta) - \tfrac{c_n}{\sqrt{n}} < d_{x,u}(s) \le q_{x,u}(\delta) + \tfrac{c_n}{\sqrt{n}}\}$. Recalling that $q_x(\delta) = d_{x,u}(s^l_{x,u}(\delta)) = F(\Cstar_{x,u}(s^l_{x,u}(\delta)),$ we can write
 $${\rm B}^l_{x,u}(\delta) = \Big\{s \le 0:\,\big|F\big(\Cstar_{x,u}(s^l_{x,u}(\delta))\big) - F\big(\Cstar_{x,u}(s))\big)\big| \le \frac{c_n}{\sqrt{n}}\Big\}.$$
 It follows that 
 \begin{align}\label{abc1}
&G(s \le 0: \wh d_{x,u}(s) \le t) \\
&\hspace*{1cm}\le G\big(s \le 0: q_{x,u}(\delta) - \tfrac{c_n}{\sqrt{n}} < d_{x,u}(s) \le t - (\wh d_{x,u} - d_{x,u})(s^l_{x,u}(\delta)) + r_n\big)
\end{align}
where $r_n$ is such that $\sup_{s \in {\rm B}^l_{x,u}(\delta)} \big|(\wh d_{x,u} - d_{x,u})(s)) - (\wh d_{x,u} - d_{x,u})(s^l_{x,u}(\delta))\big| \le r_n$.

Recall that, by assumption, $s^{\,l}_{x,u}(\delta) \in S_{x,u}(\eta^\prime_n),$ where $\eta_n^\prime = C^\prime \eta$ for $C^\prime \ge C_\epsilon > 0,$ where $C_\epsilon > 0$ will be chosen below. It is shown in Lemma~\ref{short-intervals} that $\mycomp{S}_{x,u}(\eta^\prime_n)$ is a finite union of intervals $I_j, j = 1,\ldots,J$ of length $|I_j|$ satisfying $c \eta^\prime_n \le |I_j| \le C \eta^\prime_n$, with $c,C>0$ independent of $\|x\| \le K$ and $u \in S^{d-1}.$ Each of the intervals $I_j$ contains exactly one value $s$ with $\Delta_{x,u}(s) = 0.$ Consequently, also $S_{x,u}(\eta^\prime_n)$ is a finite union of intervals $J_i$, and for each interval $J_i,$ the set of minimizers $\{C^*_{x,u}(s), s \in J_i\},$ either equals $\{A_{x,u}(s), s \in J_i\}$ or $\{B_{x,u}(s), s \in J_i\}.$ In particular this implies that the minimizers within each interval $J_i$ are nested. We now show that for appropriate choice of $C_\epsilon$, there exists an interval $J_i \subset S_{x,u}(\eta_n)$, such that both $s^{\,l}_{x,u}(\delta)$ and $B^l_{x,u}(\delta)$ lie in the same interval $J_i$. To see this, observe that $ \mycomp{S}_{x,u}(c_n) \cap \{s \le 0\}$ is a union of intervals with properties mentioned above. This is so, because $\mycomp{S}_{x,u}(c_n)$ has this property with intervals $J_i$ not containing $0$. The latter follows from $\inf_{\|x\| \le K, u \in S^{d-1}}\Delta_{x,u}(0) = \inf_{\|x\| \le K, u \in S^{d-1}}F(B_{x,u}(0)) \ge \epsilon_0 > 0$ for some $\epsilon_0 > 0.$ Also, $d^{\,l}_{x,u}(s)$, the restriction of $d_{x,u}(s)$ to $s \le 0$, is strictly monotonic. It follows by using assumption {\bf (A3)} that in order for $s, s^\prime$ to lie in different intervals $J_i$, we need to have $\big| d_{x,u}(s) - d_{x,u}(s^\prime)\big| \ge \epsilon_K c^\prime \eta^\prime_n \ge \epsilon_K c^\prime C_\epsilon \eta_n$ for some $c^\prime > 0.$ Thus, for $C_\epsilon$ such that $\epsilon_K c^\prime C_\epsilon > 1$, $s, s^\prime$ cannot both lie in $B_n(\delta)$ and in different intervals $J_i.$ All the constants used here do not depend on $x,u$.

As a consequence of all this, $s^{\,l}_{x,u}(\delta)$ and $B_{x,u}^l(\delta)$ lie in the same interval. This implies two things: (i) $B_{x,u}^l(\delta) = \big\{s \le 0: F\big(\Cstar_{x,u}(s^l_{x,u}(\delta))\Delta (\Cstar_{x,u}(s))\big) \le c_n/\sqrt{n}\big\},$ and (ii) for $s \in B_{x,u}^l(\delta)$, the minimizers in the definitions of $\wh d_{x,u}(s)$ and $d_{x,u}(s),$ respectively, are the same sets, so that  $(\wh d_{x,u} - d_{x,u})(s) = (F_n - F)(\Cstar_{x,u}(s))$ (recall here that we are on $A_n(c_n)$). 

We obtain
\begin{align*}
\sqrt{n\,}&\,\sup_{s \in {\rm B}^l_{x,u}(\delta)} \Big| (\wh d_{x,u} - d_{x,u})(s)) - (\wh d_{x,u} - d_{x,u})(s^l_{x,u}(\delta))\Big|  \\
&\le  \sup_{\big\{s \le 0: F(\Cstar_{x,u}(s^l_{x,u}(\delta)) \Delta \Cstar_{x,u}(s))| \le \frac{c_n}{\sqrt{n}}\big\}}\big|\nu_n(\Cstar_{x,u}(s^l_{x,u}(\delta))-\nu_n(\Cstar_{x,u}(s))\big|\\
&\hspace*{1cm} \le \omega_n\Big(\frac{c_n}{\sqrt{n}}\Big) = r_n.
\end{align*}
\vspace*{-1cm}

\hfill $\square$\\

The following technical lemmas are used in the proofs above. They are proven below.

\begin{lemma} \label{derivative}Suppose that $F$ has a continuously differentiable density $f$ and that {\bf (A2)} holds. For any line $\ell_{x,u} \subset \R^d$, the maps  $d^A_{x,u}(s) := F(A_{x,u}(s)), s \in \R$  and $d^B_{x,u}(s) = F(B_{x,u}(s)), s \in \R$ are twice continuously differentiable. For each $K  > 0$, the derivatives are uniformly bounded over $\|x\| \le K, s\in \R, u \in S^{d-1}.$ Moreover, for each $\epsilon^\prime, K^\prime > 0,$  there exists $\eta > 0$, such that $\big|\frac{d}{ds}d^A_{x,u}(s)\big| > \eta$ and $\big|\frac{d}{ds}d^B_{x,u}(s)\big| > \eta$ for $\epsilon^\prime < |s| < K^\prime, \|x\| \le K, u \in S^{d-1}.$
\end{lemma}

\begin{lemma} \label{short-intervals}Assume that {\bf (A1), (A2)} and {\bf(A3)} hold. For each $\epsilon, K > 0$, there exists constants $c,C, J > 0$, such that for each $(x,u) \in \R^d \times S^{d-1}$ with $\|x\| \le K,$ and each $\epsilon, \eta > 0,$ the sets $\mycomp{S}_{x,u}(\eta)$ and  $\mycomp{D}_{x,u}(\eta,\epsilon)$ consist of at most $J$ intervals $I_j$ of lengths $|I_j|, j = 1,\ldots,J$, satisfying $c\eta \le |I_j| \le C\eta.$
\end{lemma}
\begin{lemma}\label{continuity}Suppose that {\bf (A1), (A2)} and {\bf(A3)} hold. Then, for each $\epsilon, K > 0,$ the following classes of functions are equi-Lipschitz continuous:
\begin{itemize}
\item  $\big\{\delta \to s^{\,l}_{x,u}(\delta),\delta \in [\epsilon, 1 - \epsilon],\,(x,u) \in \{x \in \R^d:\|x\| \le K\} \times S^{d-1}\big\},$ 
\item $\big\{\delta \to s^{\,r}_{x,u}(\delta),\delta \in [\epsilon, 1 - \epsilon],\,(x,u) \in \{x \in \R^d:\|x\| \le K\} \times S^{d-1}\big\},$ 
\item $\big\{q_{x,u}(\delta), \delta \in [\epsilon, 1 - \epsilon], \,(x,u) \in \{x \in \R^d:\|x\| \le K\} \times S^{d-1}\big\}.$
\end{itemize}
Moreover, there exist constants $c,C > 0$, only depending on $\epsilon$ and $K$, such that, for all $\|x\| \le K, u \in S^{d-1}$, 
\begin{align}\label{q-prop}
c |\delta - \delta^\prime| \le |q_{x,u}(\delta) - q_{x,u}(\delta^\prime)| \le C |\delta - \delta^\prime|\;\; \text{for }\, \delta,\delta^\prime \in [\epsilon, 1 - \epsilon].
\end{align}
\end{lemma}
\begin{lemma}\label{compact}
The set $ \bigcup_{x: \|x\| \le K, u \in S^{d-1}}[s^l_{x,u}(\delta), s^r_{x,u}(\delta)]$ is contained in a compact set.\end{lemma}

{\sc Proof of Lemma~\ref{derivative}.} Let $u = (1,0,\ldots,0)$ and $x = (0,\ldots,0)^\prime \in {\mathbb R}^d$, and consider $F(A_{x,u}(s))$ for $s \le 0.$ Recall that all our cones have opening angle $\alpha$, and write $y = (y_1,\ldots,y_d) = (y_1,z), z \in \R^{d-1}, y_1 \in \R$. Using Fubini's theorem, we have, for $s^\prime < s < 0$,  
\begin{align*}
0 \le F(A_{x,u}(s^\prime)) &- F(A_{x,u}(s)) 
= \int^{0}_{s^\prime} {\cal S}_{s,s^\prime}(y_1) dy_1,
\end{align*}
where
\begin{align*}
{\cal S}_{s,s^\prime}(y_1) = \int_{\R^{d-1}}{\bf 1}\big\{ (\tan(\alpha) ( y_1 - s) \vee 0) \le \|z\| \le  \tan(\alpha) \big(y_1 -s^\prime \big) \big\} f(y_1, z)\, dz.
\end{align*}
Letting ${\cal B}_0(r)$ denote a $(d-1)$-dimensional ball with midpoint zero and radius $r,$ and $F_{y_1}$ the measure on $\R^{d-1}$ with Radon-Nikodym derivative $f_{y_1}(z) = f(y_1,z)$, we see that, for $s \le y_1$, the quantity ${\cal S}_{s,s^\prime}(y_1)$ equals the $F_{y_1}$-measure of the shell ${\cal B}_0(\tan(\alpha)((|s^\prime - y_1|)) \setminus {\cal B}_0(\tan(\alpha)( |s - y_1| ))$, and for $s^\prime  < y_1 < s$, it becomes the $F_{y_1}$-measure of the ball ${\cal B}_0(\tan(\alpha)((|y_1 - s^\prime|) ))$. Using polar coordinates (in $\R^{d-1}$), we have
\begin{align*}
{\cal S}_{s,s^\prime}(y_1) = \int_{\max(\tan(\alpha) (y_1 - s),0)}^{\tan(\alpha) (y_1 -  s^\prime )} \int_{S^{d-2}} r^{d-2}f_{y_1}(r v)\,d\H^{d-2}(v)\,dr,
\end{align*}
with $\H^{d-2}$ denoting $(d-2)$-dimensional Hausdorff measure on $\R^{d-1}$. We obtain, for each fixed $s \le y_1 \le 0$, 
\begin{align*}
\lim_{s^\prime \to s}\;\tfrac{{\cal S}_{s,s^\prime}(y_1)}{\tan(\alpha)|s - s^\prime|} =  \int_{S^{d-2}} \big(\tan(\alpha) (y_1 - s )\big)^{d-2}f_{y_1}(\tan(\alpha) ( y_1 - s)\, v)\,d\H^{d-2}(v),
\end{align*}
and for $ s^\prime < y_1 < s$, the limit equals zero because the volume of ${\cal B}_0\big(\tan \alpha (s - y_1)\big)$ is of smaller order than $|s - s^\prime|.$ 

Clearly, for each fixed $s$ and $y,$ the limit is finite and we obtain
%
%

\begin{figure}[h]
\begin{center}
\includegraphics[height=2.1in, width= 3.9in]{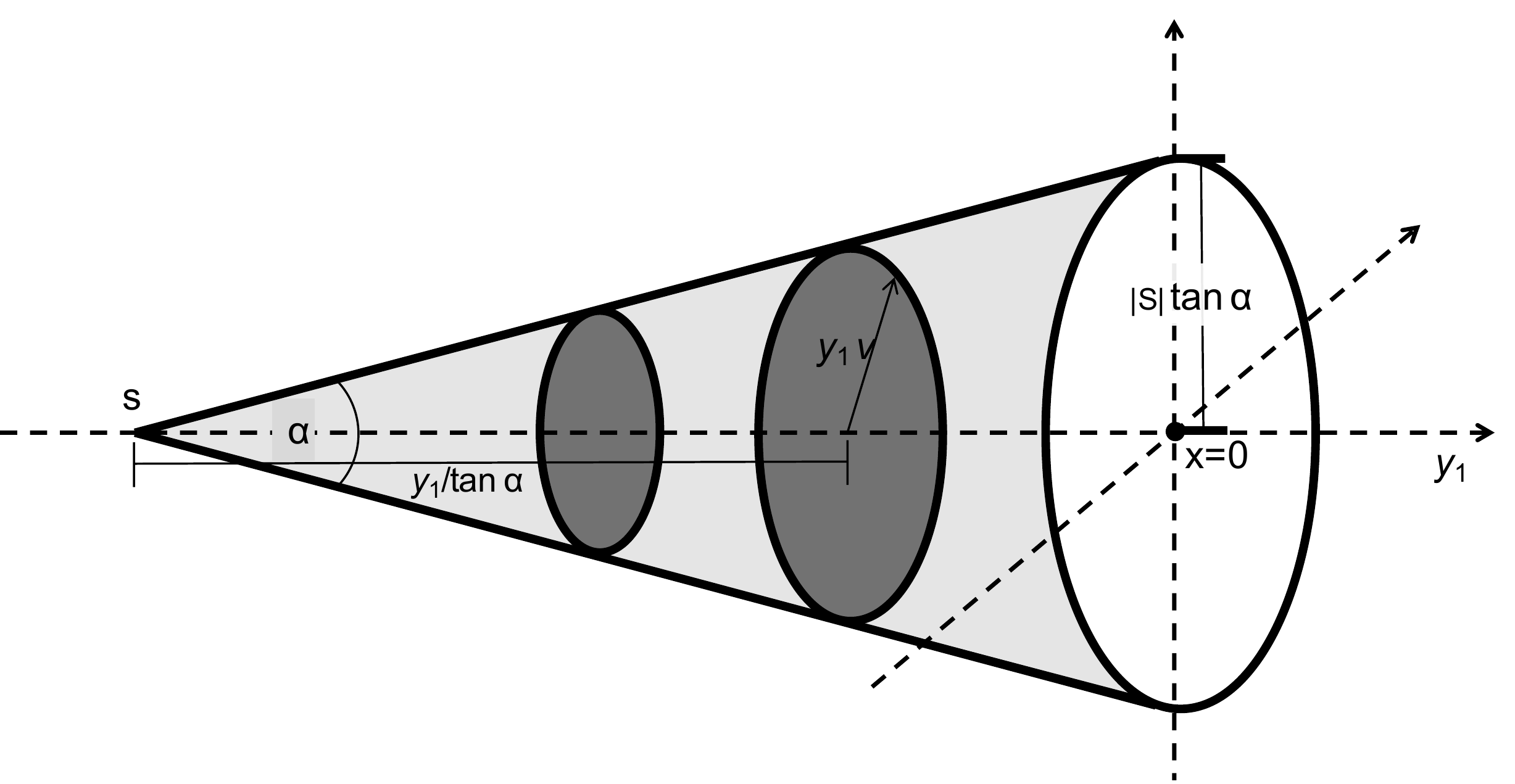}
\end{center}
\caption{\small Schematic illustration of the derivative formula of the function $s \to A_{x,u}(s)$ given in (\ref{deriv}). }
\end{figure}

%
\begin{align}
&\lim_{s^\prime \nearrow s}\tfrac{F(A_{x,u}(s)) - F(A_{x,u}(s^\prime))}{s- s^\prime} \nonumber \\
&= - \tan(\alpha)\,\int\limits^0_{s} \int\limits_{S^{d-2}} \big(\tan(\alpha) (y_1 - s)\big)^{d-2}f_{y_1}(\tan(\alpha) ( y_1 - s)\, v)\,d\H^{d-2}(v)\,dy_1\quad\nonumber \\[3pt]
&\hspace*{1cm} = -\int^{|s|\, \tan(\alpha)}_0 \int_{S^{d-2}}r^{d-2} f\big(\tfrac{r}{\tan(\alpha)} - |s|,\,rv\big)\,d\H^{d-2}(v)\, dr\label{deriv}\\[5pt]
&\hspace*{1cm} = - \int_{\partial \overset{\scriptscriptstyle \circ}{\rm A}_{0,u}(s)} f(w)\,d{\cal H}^{d-1}(w),\nonumber
\end{align}
and $\partial\overset{\scriptscriptstyle \circ}{\rm A}_{0,u}(s)$ denotes the surface of the cone $A_{0,u}(s)$ minus its base. 

For general lines $\ell_{x,u},$ a similar argument as above gives, for $s \le 0$, 
\begin{align}
&\tfrac{d}{ds}F(A_{x,u}(s)) = - \int_{\partial \overset{\scriptscriptstyle \circ}{\rm A}_{x,u}(s)} f(w)\,d{\cal H}^{d-1}(w)\nonumber\\
&= -\int\limits^{|s|\, \tan(\alpha)}_0 \int\limits_{S^{d-2}}r^{d-2} f\big(x + R_u^{-1}(\tfrac{r}{\tan(\alpha)} - |s|,\,rv)\big)\,d{\cal H}^{d-2}(v)\, dr\nonumber\\
&= {\rm sign}(s) \int\limits^{|s|\,\tan(\alpha)}_0 \int\limits_{S^{d-2}}r^{d-2} f\big(x + (\tfrac{r}{\tan(\alpha)} - |s|) u + rR_u^{-1}(0,v)^\prime\big)\,d{\cal H}^{d-2}(v)\, dr, \label{deriv-gen}
\end{align}
where $0 \in \R$, $v \in S^{d-2} \subset \R^{d-1}$, so that $(0,v)^\prime \in \R^d$, and $R_u$ denotes the $(d \times d)$ rotation matrix with $u = R_u^{-1} e_1$. For $s > 0$, the derivative has a similar form.

Differentiability of the function $s \to F(B_{x,u}(s))$ follows by a similar argument. Here the form of the derivative is, 
\begin{align*}
&\tfrac{d}{ds}F(B_{x,u}(s)) \\
&= {\rm sign}(s) \int\limits_{|s|\,\tan(\alpha)}^\infty \int\limits_{S^{d-2}}r^{d-2} f\big(x + \big(\tfrac{r}{\tan(\alpha)} - |s|\big) u + rR_u^{-1}(0,v)^\prime\big)\,d{\cal H}^{d-2}(v)\, dr.
\end{align*}
Notice that the outer integral is over an unbounded range. Assumption {\bf (A2)(i)} implies that this integral still is bounded uniformly in $u \in S^{d-1}$ for each given $s$ and $x$. The argument is as follows: With $C$ from assumption {\bf(A2)(i)}, we have
\begin{align*}
&\int_{|s|\,\tan(\alpha)}^\infty \int_{S^{d-2}}r^{d-2} f\big(x + \big(\tfrac{r}{\tan(\alpha)} - |s|\big) u + rR_u^{-1}(0,v)^\prime\big)\,d{\cal H}^{d-2}(v)\, dr\\
\le &\; C \int_{|s|\,\tan(\alpha)}^\infty \int_{S^{d-2}}\frac{r^{d-2}}{\|x + \big(\tfrac{r}{\tan(\alpha)} - |s|\big) u + rR_u^{-1}(0,v)^\prime\|^{d-1+\epsilon}}\,d{\cal H}^{d-2}(v)\, dr.
\end{align*}
By using orthogonality of $u$ and $R_u^{-1}(0,v)^\prime$,
$$\big\|\big(\tfrac{r}{\tan(\alpha)}- |s|\big) u + rR_u^{-1}(0,u)^\prime\big\| = \sqrt{|\big(\tfrac{r}{\tan(\alpha)}- |s|\big)\big|^2 + r^2\,} \ge r,$$
and thus, for $r \ge 2\|x\|,$ we have  $\|x + (s + \tfrac{r}{\tan(\alpha)}) u + rR_u^{-1}(0,u)^\prime\| \ge \big|\,\|x\| - \|(s + \tfrac{r}{\tan(\alpha)}) u + rR_u^{-1}(0,v)^\prime\|\big| \ge r/2.$ Consequently,
\begin{align*}
&\int_{|s|\,\tan(\alpha)}^\infty \int_{S^{d-2}}r^{d-2} f\big(x + \big(\tfrac{r}{\tan(\alpha)} - |s|\big) u + rR_u^{-1}(0,v)^\prime\big)\,d{\cal H}^{d-2}(v)\, dr \\
&=  \Big( \int\limits_{|s|\, \tan(\alpha)}^{2\|x\|} + \int\limits_{2\|x\|}^\infty\Big) \int\limits_{S^{d-2}}r^{d-2} f\big(x + \big(\tfrac{r}{\tan(\alpha)}- |s| \big) u + rR_u^{-1}(0,v)^\prime\big)\,d{\cal H}^{d-2}(v)\, dr\\
&\le C + 2^{d-1+\epsilon}\,A_{d-2} \int_{2\|x\|}^\infty r^{-(1 +\epsilon)} \,dr < \infty,
\end{align*}
with $A_{d-2}$ denoting the surface area of $S^{d-2}$. Note that, if $\|x\| \le K$, we can simply split the first integral in these calculations into $\int_{|s| \tan \alpha}^{2K} \cdots$ and $\int_{2K}^\infty\cdots$, and we see that the resulting bound only depends on $K,$ and not on $x$. 

In order to derive the second order derivative of these functions, we use chain rule and assumption {\bf (A2)(ii)}. From (\ref{deriv-gen}), we obtain, for $s < 0$, by using chain rule 
\begin{align*}
\tfrac{d^2}{ds^2} F(A_{x,u}(s)) &= -|s|^{d-2}\tan(\alpha)^{d-1} \int\limits_{S^{d-2}}f(x + |s| \tan(\alpha) R_u^{-1}(0,v)^\prime)\,d{\cal H}^{d-2}(v) \\
&- \int_{|s| \tan(\alpha)}^0 \int_{S^{d-2}}r^{d-2}\big({\rm grad}f(x + r R_u^{-1}(0,v)^\prime)\big)^\prime u\,d{\cal H}^{d-2}(v) dr.
\end{align*}
Observe that the integrand of the second term depends on $f$ only though its directional derivatives in direction $u.$ A similar formula holds for the second derivative of $F(B_{x,u}(s))$, and as for the first derivative, the corresponding expression involves integration over unbounded domains. Finiteness of these integrals follows from assumptions {\bf (A2)(i)} and {\bf (A2)(ii)}, respectively, by using arguments similar to the ones that have been used above to show finiteness of the first derivative. Continuity of the second derivative follows by using the assumed continuity of the gradient of $f$.

The fact that the first derivatives are uniformly bounded, and that they are uniformly bounded away from zero for $s$ bounded away from zero and infinity, follows directly form the form of the derivative. \hfill $\square$\\

{\sc Proof of Lemma~\ref{compact}:} Let $\epsilon > 0$ and let $R_\epsilon \ge 1$ be such that $F\big({\cal B}_{R_\epsilon}(0) \big) \ge 1-\epsilon.$  Fix $(x_0,u_0),$ and let $(x,u)$ be such that $\|x - x_0\| \le \epsilon/R_{\epsilon}$ and $\big|1 - \langle u_0,u\rangle\big| \le \epsilon/R_{\epsilon}$. For a given $s \in \R$, let $\ell$ and $\ell_0$ denote the axes of symmetry of the cones $C_{x,u}(s)$ and $C_{x_0,u_0}(s)$, respectively. These axes are given by $y(t) = x + tu$ and $y_0(t) = x_0 + tu_0$, $ t \in \R$, and since
\begin{align*}
\|y(t) - y_0(t)\| &\le \|x_0 - x\| + \|u_0 - u\||t| \le \tfrac{\epsilon}{R_{\epsilon}} + 2(1 - (1-\tfrac{\epsilon}{R_{\epsilon}}))|t| \\
&= \tfrac{\epsilon}{R_{\epsilon}} + 2\tfrac{\epsilon}{R_{\epsilon}} |t|,
\end{align*}
we can conclude that the Hausdorff distance $H$ of the two axes restricted to the ball ${\cal B}_{R_\epsilon}(0)$ satisfies $H(\ell \cap {\cal B}_{R_\epsilon}(0),\ell_0 \cap {\cal B}_{R_\epsilon}(0)) \le \epsilon + 4\epsilon = 5\epsilon.$ Simple geometric considerations now show that $C_{x,u}\big(s - \frac{5\epsilon}{\cos \alpha}\big) \cap {\cal B}_{R_\epsilon}(0) \subset C_{x_0,u_0}(s) \cap {\cal B}_{R_\epsilon}(0) \subset C_{x,u}\big(s+\frac{5\epsilon}{\cos \alpha}\big) \cap {\cal B}_{R_\epsilon}(0).$ Using the assumed boundedness of the density $f$, we can thus find $c, \epsilon > 0$ with 
\begin{align*}
\sup_{\|x - x_0\| \le \epsilon/R_\epsilon\atop |1 - \langle u,u_0\rangle| \le \epsilon/R_\epsilon}&\big|F\big(A_{x,u}(s)\big) - F\big(A_{x_0,u_0}(s)\big)| \\
&\le \big|F\big(\big(A_{x,u}(s) \Delta A_{x_0,u_0}(s)\big) \cap {\cal B}_{R_\epsilon}(0)\big) \big| +  \epsilon \le c\,\epsilon
\end{align*}
and the same holds for the sets $B_{x,u}(s)$. Notice that this shows equi-Lipschitz-continuity of the class $\{(x,u) \to d_{x,u}(s), s \in \R,\, (x,u) \in \R^d \times S^{d-1}\}$. We obtain that, for $(x,u)$ with $\|x - x_0\| \le \epsilon/R_{\epsilon}$ and $\big|1 - \langle u_0,u\rangle\big| \le  \epsilon/R_{\epsilon},$ there exists $c > 0,$ only depending on $\epsilon$, such that, for any $\delta > 0$,
\begin{align*}
q_{x,u}(\delta) &= \inf\big\{t: G\big(s: d_{x,u}(s) \le t\big) \ge \delta\big)\\
& \le \inf\big\{t: G\big(s: d_{x_0,u_0}(s) \le t - \hspace*{-0.2cm}\sup_{{\|x - x_0\| \le \epsilon/R_\epsilon\atop |1 - \langle u,u_0\rangle| \le \epsilon/R_\epsilon}\atop s\in \R} \big|d_{x,u}(s) - d_{x_0,u_0}(s)\big|\big) \ge \delta\big)\\
& \le \inf\big\{t: G\big(s: d_{x_0,u_0}(s) \le t - c\,\epsilon\big) \ge \delta\big) = q_{x_0,u_0}(\delta) + c\,\epsilon,
\end{align*}
and similarly, $q_{x,u}(\delta) \ge q_{x_0,u_0}(\delta) - c\,\epsilon$. It follows that, with this choice of $\epsilon$ and $c$,
\begin{align*}
[s^l_{x,u}(\delta), &s^r_{x,u}(\delta)] 
= \{s \in \R:\, d_{x,u}(s) \le q_{x,u}(\delta)\}  \\
&\subset \Big\{s \in \R: d_{x_0,u_0}(s) \le q_{x_0,u_0}(\delta) \\
&\hspace*{1.5cm}+ \sup_{x,u,s}\Big| d_{x,u}(s) - d_{x_0,u_0}(s)\Big| + \sup_{x,u}\Big|  q_{x,u}(\delta) -  q_{x_0,u_0}(\delta) \Big|\Big\} \\
& \subset \Big\{s \in \R: d_{x_0,u_0}(s) \le q_{x_0,u_0}(\delta) + 2c\,\epsilon\Big\},
\end{align*}
where, as above, suprema over $x, u$ and $s$ are extended over $\|x - x_0\| \le \epsilon/R_\epsilon, |1 - \langle u,u_0\rangle| \le \epsilon/R_\epsilon, s\in \R.$ As indicated at the beginning of the proof of Theorem~\ref{as-distr-quant}, for $\epsilon > 0$ small enough, so that $q_{x_0,u_0}(\delta) + 2c\, \epsilon < q_{x_0,u_0}(1)$, the set in the last displayed formula is a finite interval. Finally, for every given $\epsilon > 0$, there exist finitely many points $(x_i,u_i), i = 1,\ldots,L$ such that the balls ${\cal B}_{\epsilon}(x_i), i = 1,\ldots,L$ cover $\{\|x\| \le K\},$ and the sets $\{u \in S^{d-1}: |\langle u,u_i\rangle| \le \epsilon\}, i = 1,\ldots,L$  cover $S^{d-1}$. Using the above, we have
\begin{align*}
\bigcup_{i=1}^L  \Big\{s \in \R: d_{x_i,u_i}(s) \le q_{x_0,u_0}(\delta) + 2c\,\epsilon\Big\} \supset \bigcup_{\{x: \|x\| \le K, u \in S^{d-1}\}}[s^l_{x,u}(\delta), s^r_{x,u}(\delta)].
\end{align*}
As a finite union of compact sets, the left hand side is compact. This is what we wanted to show. \hfill$\square$\\

{\sc Proof of Lemma~\ref{short-intervals}.} Assumption {\bf (A3)} says that the set $\{s: \Delta_{x,u}(s) = 0\}$ consist of isolated points. More precisely, since the second derivatives $\frac{d^2}{ds^2}F(A_{x,u}(s))$ and $\frac{d^2}{ds^2}F(B_{x,u}(s))$ are uniformly bounded over $u \in S^{d-1}, \|x\| \le K$, we can conclude from assumption {\bf (A3)} that, for $\eta > 0$ sufficiently small, 
$$\sup_{\|x\| \le K \atop u \in S^{d-1}} \sup_{\{s: \Delta_{x,u}(s)\le \eta\}} \Big|\frac{d}{ds}F(A_{x,u}(s)) - \frac{d}{ds}F(B_{x,u}(s))\Big| > \epsilon_K/2,$$ 
where $\epsilon_K$ is from assumption {\bf (A3)}. Combined with the fact that the first derivatives are bounded, this also implies that there exist constants $c,C > 0$ such that, for each $s_j \in \{s:\Delta_{x,u}(s) = 0\}$, and $\eta > 0$ sufficiently small, there exists an interval $\tilde I_j$ with $s_j \in \tilde I_j$, and
$$\mycomp{S}_{x,u}(\eta) = \{s: | \Delta_{x,u}(s)| \le \eta \big\}  = \bigcup_{j} \tilde I_j$$
with the union consisting of disjoint intervals, and
$$c\eta \le \text{vol}(\tilde I_j) \le C \eta,\qquad j = 1,2,\ldots$$ 
Note that the constants $c,C$ do not depend on $(x,u).$ Since the functions $\delta \to s_{x,u}^{\,l}(\delta)$ and $\delta \to s_{x,u}^{\,r}(\delta)$ with $\delta$ ranging over $[\epsilon, 1 - \epsilon]$ are monotonically increasing and differentiable with derivatives bounded away from $0$ and infinity, uniformly over $\|x\| \le K$ and $u \in S^{d-1},$ a similar property holds for
the set $D_{x,u}(\eta)$, meaning that $D_{x,u}(\eta)$ can be written as a disjoint union of intervals, the number of which is bounded uniformly over $\|x\| \le K$ and $u \in S^{d-1}$.\hfill $\square$\\

{\sc Proof of Lemma~\ref{continuity}.}  Recall that $[s^{\,l}_{x,u}(\delta), s^{\,r}_{x,u}(\delta) ] = \{s \in \R: d_{x,u}(s) \le q_{x,u}(\delta)\},$ where $q_{x,u}(\delta)$ is the inverse of the function $\gamma \to G(s: d_{x,u}(s) \le \gamma), \gamma \ge 0.$  Also recall that $d^{\,l}_{x,u}(s)$ and $d^{\,r}_{x,u}(s),$ are the restrictions of $d_{x,u}(s)$ to $s \le 0$ and $s > 0$, respectively.

We first note that both $d^{\,l}_{x,u}(s)$ and $d^{\,r}_{x,u}(s),$ $(x,u) \in \{\|x\| \le K\} \times  S^{d-1}$ are equi-Lipschitz continuous classes of functions. This follows from the fact that $s\to F(A_{x,u}(s))$ and $s \to F(B_{x,u}(s))$ $(x,u) \in \{\|x\| \le K\} \times  S^{d-1}$ are differentiable functions with derivatives uniformly bounded over $s \in \R, \|x\| \le K, u \in S^{d-1}$ (see Lemma~\ref{derivative}). Moreover, both $d^{\,l}_{x,u}(s)$ and $d^{\,r}_{x,u}(s),$ are piecewise differentiable. For each $\epsilon^\prime, K^\prime > 0$, the derivatives of these functions are also bounded away from zero on $\epsilon^\prime \le |s| \le K^\prime$ (uniformly over $\|x\| \le K$ and $u \in S^{d-1}$), it follows that, for each $0 < \epsilon^{\prime\prime} \le 1/2,$ the inverses $\big(d^{\,r}_{x,u}\big)^{-1}(\gamma)$ and $\big(d^{\,l}_{x,u}\big)^{-1}(\gamma), \gamma \in [\epsilon^{\prime\prime}, 1 - \epsilon^{\prime\prime}]$ are differentiable with derivatives bounded away from zero and infinity, uniformly over $\|x\| \le K$ and $u \in S^{d-1}$. This, along with the fact that by assumption $g$ is continuous and strictly positive, implies that the class $\{[0,1]) \ni \gamma \mapsto q_{x,u}^{-1}(\gamma) := G(s: d_{x,u}(s) \le \gamma), \|x\| \le K, u \in S^{d-1}\}$ is also an equi-Lipschitz continuous class of functions. Moreover, these functions are also piecewise differentiable, and, for each $\epsilon > 0$, their derivatives are uniformly bounded away from zero and infinity for $\gamma \in [\epsilon, q_{x,u}(1) - \epsilon]$. Similar properties holds for the class of inverse functions $\{q_{x,u}(\delta), \delta \in [\epsilon, 1 - \epsilon], \|x\| \le K, u \in S^{d-1}\}$, where $\epsilon > 0$. This follows by using that $q_{x,u}(\delta)$ is piecewise differentiable, with $\frac{d}{d\delta}q_{x,u}(\delta) = 1/\big(q_{x,u}^{-1}\big)^\prime(q_{x,u}(\delta)),$ and that, for each $\epsilon > 0,$ there exist $\gamma_0 > 0$ such that%
\begin{align}\label{lb}
\inf_{\|x\| \le K, u \in S^{d-1}}\inf_{\delta \ge \epsilon} q_{x,u}(\delta) =  \inf_{\|x\| \le K, u \in S^{d-1}}q_{x,u}(\epsilon) > \gamma_0.
\end{align}
To see (\ref{lb}), first observe that our assumptions imply the existence of $s_0 > 0$, such that for all $\|x\| \le K$ and $u \in S^{d-1}$, $d_{x,u}(s) = F(A_{x,u}(s))$ for all $|s| < s_0.$ Since ${\rm vol}(A_{x,u}(s)) = C_d |s|^{d}$, for some constant $C_d > 0$ (not depending on $x,u$), we can conclude that there exists constants, $c,C > 0$ such that, for all $\|x\| \le K, u \in S^{d-1}$,  $c |s|^d \le d_{x,u}(s) \le C |s|^d$, for $|s| \le s_0$. Using the fact that $d^{\,l}_{x,u}(s^{\,l}_{x,u}(\delta)) =  d^{\,r}_{x,u}(s^{\,r}_{x,u}(\delta)) = q_{x,u}(\delta)$, we obtain that, for $q_{x,u}(\delta) \le d_{x,u}(s_0)$, 
\begin{align*}
\frac{1}{C^{1/d}}\, q^{1/d}_{x,u}(\delta) \le& s^{\,r}_{x,u}(\delta) \le \frac{1}{c^{1/d}} \,q^{1/d}_{x,u}(\delta)\quad\text{and}\\
 -\frac{1}{c^{1/d}}\, q^{1/d}_{x,u}(\delta) \le &s^{\,l}_{x,u}(\delta) \le -\frac{1}{C^{1/d}} \,q^{1/d}_{x,u}(\delta).
\end{align*}
Recalling that by definition of $s^{\,l(r)}_{x,u}(\delta)$, $\delta = G(s^{\,r}_{x,u}(\delta)) - G(s^{\,l}_{x,u}(\delta))$, we obtain
\begin{align*}
\delta &\le G\Big(\frac{1}{c^{1/d}} \,q^{1/d}_{x,u}(\delta)\Big) - G\Big(-\frac{1}{c^{1/d}} \,q^{1/d}_{x,u}(\delta)\Big) 
 \le \max_{s\in \R}g(s)\,\frac{2}{c^{1/d}}q^{1/d}_{x,u}(\delta).
\end{align*}
or
\begin{align}\label{bb}
q_{x,u}(\delta) \ge c_d \delta^{d}\quad \text{for all $\delta$ with }\, \delta \le \big(q_{x,u}\big)^{-1} \big( d_{x,u}(s_0)\big),
\end{align}
where $c_d > 0$ is a constant not depending on $(x,u) \in \{\|x\| \le K\} \times S^{d-1}.$ Let $\delta_{x,u} = \big(q_{x,u}\big)^{-1}(d_{x,u}(s_0))$. Our assumptions imply that $\delta_{x,u}$ is a continuous function of $(x,u)$, and thus there exists $\widetilde \delta_0 > 0$  such that $\inf_{\|x\| \le K, u \in S^{d-1}} \delta_{x,u} \ge \widetilde\delta_0,$ so that $q_{x,u}(\delta) \ge c_d\delta,$ for all $\delta \le \widetilde\delta_0,$ with $\widetilde\delta_0$ not depending on $(x,u) \in \{\|x\| \le K\} \times S^{d-1}$. W.l.o.g. assume that $\epsilon \le \widetilde \delta_0,$ then $\inf_{\|x\| \le K, u \in S^{d-1}}\inf_{\delta \ge \epsilon} q_{x,u}(\delta) =  \inf_{\|x\| \le K, u \in S^{d-1}}q_{x,u}(\epsilon) > c_d \epsilon,$ i.e. (\ref{lb}) holds.
 
Finally, since $s_{x,u}^{\,l}(\delta) = \big(d_{x,u}^{\,l}\big)^{-1}(q_{x,u}(\delta))$, equi-Lipschitz continuity of this class of functions follows, and the same holds for the functions $s_{x,u}^{\,r}(\delta).$ \hfill $\square$
\section{References}
\begin{description}
\item Aeberhard, S., Coomans, D. and de Vel, O. (1992): 
Comparison of Classifiers in High Dimensional Settings, 
Tech. Rep. no. 92-02, Dept. of Computer Science and Dept. of 
Mathematics and Statistics, James Cook University of North Queensland. 
\item Agostinelli, C. (2018): Local half-region depth for functional data. {\em J. Multivariate Anal.}, {\bf 163}, 67 - 79.
\item Agostinelli, C. and Romanazzi, M. (2011): Local depth. {\em J. Statist. Plan. Inference} {\bf 141}, 817 - 830.
\item Alon, U., Barkai, N., Notterman, D.A., Gish, K., Ybarra, S., Mack, D. and Levine, AJ. (1999): Broad patterns of gene expression revealed by clustering analysis of tumor and normal colon tissues probed by oligonucleotide arrays. {\it Proc. Natl. Acad. Sci. USA} {\bf 96}(12), 6745-6750.
\item  Arcones, M.A., and Gin\'{e}, E. (1993): Limit theorems for U-processes. {\it Ann Probab.}, {\bf 21}, 1494 - 1542.
\item Berthet, P. and Mason, D. M. (2006):  Revisiting two strong approximation results of Dudley and Philipp, In {\it IMS Lecture Notes - Monograph Series, High Dimensional Probability,} Gin\'e, E., Koltchinskii, V., Li, W., and Zinn, J. eds., {\bf 51}, 155 - 172.
\item Bubenik, P. (2015). Statistical topological data analysis using persistence landscapes. {\em J. Mach. Learn. Res.}, {\bf 16}, 77-102.
\item Bubenik, P. and Kim, P. T. (2007). A statistical approach to persistent homology. {\em Homology Homotopy Appl.}, {\bf 9}, 337-362.
\item Bousquet O., Boucheron S., and Lugosi G. (2004): Introduction to Statistical Learning Theory. In: Bousquet O., von Luxburg U., R\"{a}tsch G. (eds) {\em Advanced Lectures on Machine Learning. ML 2003. Lecture Notes in Computer Science, vol. 3176.} Springer, Berlin, Heidelberg.
\item Campos, G. O., Zimek, A.,  Sander, J., Campello, R. J. G. B. , Micenková, B., Schubert, E.,  Assent, I., and Houle, M. E. (2016): On the Evaluation of Unsupervised Outlier Detection: Measures, Datasets, and an Empirical Study. {\it 
Data Mining and Knowledge Discovery} {\bf 30}, 891-927. 
\item Chaudhuri, P. and Marron, J.S. (1999): SiZer for exploration of structure in curves. {\it J. Amer. Stat. Assoc.} {\bf 94}, 807-823.
%
%
\item Cuturi, M. (2010): Positive definite kernels in machine learning. {\it Technical report.}
\item Dutta, S., Sarkar, S., and Ghosh, A.K. (2016): Multi-scale classification using localized spatial depth. {\it J. Machine Learn. Res.}, {\bf 17}, 1 - 30.
\item Einmahl, J.H.J. and Mason, D. (1992): Generalized quantile processes. {\it Ann. Stat.} {\bf  20}, 1062-1078.
\item Elmore, R.T., Hettmansperger, T.P. and Xuan, F. (2006): Spherical data depth and a multivariate median, In: {\em Data Depth: Robust Multivariate Analysis, Computational
Geometry and Applications}, DIMACS Ser. Discrete Math. Theoret. Comput. Sci., {\bf 72}, 87-101.
\item Fraiman, R., Gamboa, F. and Moreno, L. (2019): Connecting pairwise geodesic spheres by depth: DCOPS. {\em J. Multivariate Anal.}, {\bf 169}, 81-94.
\item Genton, M. (2001): Classes of kernels for machine learning: A statistics perspective. {\it J. Machine Learn. Res.} {\bf 2}, 299-312.
%
%
%
\item Kotik, L. and Hlubinka, D. (2017): A weighted localization of halfspace depth
and its properties. {\em J. Multivariate Anal.}, {\bf 157}, 53-69.
%
\item Leng, X. and M\"{u}ller, H.G. (2006): Classification using functional data analysis for temporal gene expression data. {\em Bioinformatics} {\bf 22}, 68-76.
\item Liu, R.Y. (1988): On a notion of simplicial depth, {\em Proceedings of the National Academy of Sciences}, {\bf 85}, 1732-1734.
\item Liu, R.Y. (1990): On a notion of data depth based on random simplices, {\em Ann. Statist.}, {\bf 18}, 405-414.
\item Massart, P. (1989): Strong approximation for multivariate empirical and related processes, via KMO constructions. {\it Ann. Probab.} {\bf 17}, 266 - 291.
\item Mass\'{e}, J-C. (2004): Asymptotics for the Tukey depth process with an application to a multivariate trimmed mean. {\it Bernoulli} {\bf 10}, 397-419.
\item Minotte, M. and Scott, D. (1993): The mode tree: A tool for visualization of nonparametric density features, {\it J. Comput. Graph. Statist.} {\bf 2}, 51-68.
\item Paindaveine, D. and van Bever, G. (2012): From depth to local depth: A focus on centrality. {\it JASA}, {\bf 108:503}, 1105 - 1119.
\item Pham, N. (2018): L1-Depth Revisited : A Robust Angle-based Outlier Factor in High-dimensional Space. {\it Machine Learning and Knowledge Discovery in Databases}, pp 105-121.
\item R\'{e}nyi, A. and Sulanke, R. (1963): \"{U}ber die konvexe H\"{u}lle von n zuf\"{a}llig gerw\"{a}hten Punkten I. Z. Wahrsch. Verw. Gebiete, {\bf 2}, 75?84.
\item Serfling, R. (2019): Depth Functions on General Data Spaces, I. Perspectives, with Consideration of ?Density? and ?Local? Depths. {\em Available at} \url{www.utdallas.edu/~serfling/papers/I.Perspectives.pdf}
\item Shawe-Taylor, J. and Cristianini, N. (2004): {\it Kernel methods for pattern analysis.} Cambridge University Press.
\item Spellman, P.T., Sherlock, G., Zhang, M.Q., Iyer, V.R., Anders, K., Eisen, M.B.,  Brown, P.O., Botstein, D., and Futcher, B. (1998): Comprehensive identification of cell cycle-regulated genes of the yeast saccharomyces cerevisiae by microarray hybridization.  {\it Molecular Biology of the Cell} {\bf 9}, 3273 - 3297. 
%
%
\item Ting, K.M., Zhou, G-T., Liu, F.T., and Tan, S.C. (2013): Mass estimation. {\it  Mach. Learn.}, {\bf 90}, 127-160.
%
%
\item van der Vaart, A. W. and Wellner, J. A. (1996): {\it Weak convergence of empirical processes: With applications to Statistics.} Springer
\item Zuo, Y. and Serfling, R. (2000): Notions of statistical depth functions. {\em Ann. Statist.}, {\bf 28}, 461-482.
\end{description}
\end{document}